\documentclass[12pt]{report}   
\usepackage{graphicx}  
\usepackage[letterpaper, margin=3cm]{geometry}
\usepackage{setspace}  
\usepackage{times}  
\usepackage[explicit]{titlesec}  
\usepackage[titles]{tocloft}  
\usepackage[backend=bibtex, sorting=none, bibstyle=alphabetic, citestyle=alphabetic, sorting=nyt]{biblatex}  
\usepackage[bookmarks=true, hidelinks]{hyperref}
\usepackage[page]{appendix}  
\usepackage{rotating}  
\usepackage[normalem]{ulem}  

\usepackage{mypackage} 

\usepackage{afterpage}

\bibliography{references}

\AtEveryBibitem{\clearfield{issn}}
\AtEveryBibitem{\clearlist{issn}}

\AtEveryBibitem{\clearfield{language}}
\AtEveryBibitem{\clearlist{language}}

\AtEveryBibitem{\clearfield{doi}}
\AtEveryBibitem{\clearlist{doi}}

\AtEveryBibitem{\clearfield{url}}
\AtEveryBibitem{\clearlist{url}}

\AtEveryBibitem{%
  \ifentrytype{online}
    {}
    {\clearfield{urlyear}\clearfield{urlmonth}\clearfield{urlday}}}

\begin{document}

\onehalfspacing



\newcommand{\thesisTitle}{\Large Twisted Equivariant Tate K-Theory}
\newcommand{\yourName}{Thomas Dove}
\newcommand{\yourSchool}{Mathematics and Statistics}
\newcommand{\yourMonth}{October}
\newcommand{\yourYear}{2019}

\newcommand{\HRule}[1]{\rule{\linewidth}{#1}}


\begin{titlepage}

\begin{center}

\begin{singlespacing}
\phantom{XX}
\vspace{4\baselineskip}
\textbf{\MakeUppercase{\thesisTitle}}\\
\vspace{3\baselineskip}
By\\
\vspace{3\baselineskip}
{\Large \yourName}\\
\vspace{3\baselineskip}
Supervised by \\
Nora Ganter \\
\vspace{3\baselineskip}
Submitted In Partial Fulfillment\\
of the Requirements for the Degree\\
Master of Science in the\\
School of \yourSchool\\
\vspace{3\baselineskip}
{\Large University of Melbourne}\\
\vspace{\baselineskip}
\yourMonth{} \yourYear{}
\vfill

\end{singlespacing}

\end{center}

\end{titlepage}

\currentpdfbookmark{Title Page}{titlePage}  

\pagenumbering{roman}

\setcounter{page}{2} 

\phantom{X}\clearpage

\clearpage

\begin{centering}
\textbf{ABSTRACT}\\
\vspace{\baselineskip}
\end{centering}

Starting with a $\C^*$-valued cocycle on the global quotient orbifold $X \mmod G$, we apply transgression techniques for 2-gerbes, as developed by Lupercio and Uribe, to construct a gerbe on the orbifold loop space $\L(X \mmod G)$. This gives a loop based definition of twisted equivariant Tate K-theory for finite groups that conforms to the definition that Luecke provides for compact, connected Lie groups. We relate our construction to Ganter's work on Moonshine and Huan's quasi-elliptic cohomology.

\vspace{2cm}

\begin{centering}
\textbf{ACKNOWLEDGEMENTS}\\
\vspace{\baselineskip}
\end{centering}

To my supervisor, Nora Ganter: you have given me a warm welcome to the world of research, and for this I thank you. To all of my friends in Peter Hall: you have made coming into university each day a true joy. To my dear friends and family who exist outside of my mathematics life: your support means the world to me, even if you don't quite know what it is that I do all day. To the several people who have proofread this thesis: it has been a pleasure fighting by your side in the endless war against typos.

\clearpage

\phantom{X}\clearpage 

\addtocontents{toc}{\cftpagenumbersoff{chapter}} 

\addtocontents{toc}{\cftpagenumberson{chapter}}


\renewcommand{\cftchapdotsep}{\cftdotsep}  
\renewcommand{\cftchapfont}{\bfseries}  
\renewcommand{\cftchappagefont}{}  
\renewcommand{\cftchappresnum}{Chapter }
\renewcommand{\cftchapaftersnum}{:}
\renewcommand{\cftchapnumwidth}{5em}
\renewcommand{\cftchapafterpnum}{\vskip\baselineskip} 
\renewcommand{\cftsecafterpnum}{\vskip\baselineskip}  
\renewcommand{\cftsubsecafterpnum}{\vskip\baselineskip} 
\renewcommand{\cftsubsubsecafterpnum}{\vskip\baselineskip} 

\titleformat{\chapter}[display]
{\normalfont\bfseries\filcenter}{\chaptertitlename\ \thechapter}{0pt}{\MakeUppercase{#1}}

\renewcommand\contentsname{Table of Contents}

\setcounter{tocdepth}{1}

\begin{singlespace}
\tableofcontents
\end{singlespace}

\currentpdfbookmark{Table of Contents}{TOC}

\clearpage


\clearpage
\pagenumbering{arabic}
\setcounter{page}{1} 

\titleformat{\chapter}[display]
{\normalfont\bfseries\filcenter}{\MakeUppercase\chaptertitlename\ \thechapter}{0pt}{\MakeUppercase{#1}}  
\titlespacing*{\chapter}
  {0pt}{0pt}{30pt}	
  
\titleformat{\section}{\large\bfseries}{\thesection}{1em}{#1}

\titleformat{\subsection}{\normalfont\bfseries}{\thesubsection}{1em}{#1}

\titleformat{\subsubsection}{\normalfont\itshape}{\thesubsection}{1em}{#1}


\chapter{Introduction}

Tate K-theory is an elliptic cohomology theory built from the $S^1$-equivariant K-theory of the loop space $LX$. In the equivariant setting this becomes more subtle, as $LX$ is replaced with the orbifold loop groupoid $\L(X \mmod G)$. When $G$ is finite, a definition of $G$-equivariant Tate K-theory is provided by Ganter in \cite[\textsection 3.1]{ganterstringy}, where it is attributed to the ideas of Devoto \cite{devoto1996}. In the case where $G$ is a Lie group, there is a construction by Huan in terms of quasi-elliptic cohomology \cite{Huan:thesis}. In this thesis, we build on these ideas to define twisted $G$-equivariant Tate K-theory, where $G$ is a finite group.

The twistings of equivariant K-theory correspond to elements of $H^3(X_G;\Z)$, where $X_G$ is the ``Borel construction'' of $X$, see \cite{AStwisted}. In \cite{ando2010}, the authors make the case that twistings of equivariant elliptic cohomology are associated to elements of $H^4(X_G;\Z)$. Being an equivariant elliptic cohomology theory, it is expected that equivariant Tate K-theory is twisted by $H^4(X_G;\Z)$ as well.

This prediction holds true in the case that $X$ is a point; the formalism appears, for instance, in the context of Moonshine. Let $M$ be the Monster group. Norton's generalised Moonshine conjecture, published in the appendix of \cite{Mason:Norton}, is concerned with the existence of a family of modular functions that can be interpreted as elements of the $M$-equivariant Tate K-theory of a point. As this conjecture was explored further, the existence of representations that are twisted by a particular $\alpha \in H^4(BM;\Z)$ became evident \cite{Mason:unpublished}. The properties of $\alpha$ are slowly being uncovered, for instance in \cite{JohnsonFreyd:MoonshineAnomaly}, Johnson-Freyd proves that $\alpha$ is of order 24 and is not a Chern class. Due to the presence of these twists, the modular functions in the generalised Moonshine conjecture can be interpreted as elements of the $\alpha$-twisted $M$-equivariant Tate K-theory of a point, see \cite{ganter2009}. We expand on this in Chapter 6.

Elements of $H^4(X_G;\Z)$ can be interpreted as 2-gerbes on the global quotient orbifold $X \mmod G$, up to a suitable notion of isomorphism. The holonomy of a 2-gerbe with connective structure on $X \mmod G$ produces a gerbe on the loop groupoid $\L(X \mmod G)$. This process can be described as a transgression map on Deligne cohomology - this is discussed in Chapter 5. In the orbifold setting, the theory of gerbes and the holonomy of gerbes was developed by Lupercio and Uribe in \cite{LU04} and \cite{LUhol}. Twisted equivariant Tate K-theory is modelled on the K-theory of $\L(X \mmod G)$, twisted by the gerbe we obtain via transgression.

In his recent paper \cite{Kiranspaper}, Luecke has formulated twisted $G$-equivariant Tate K-theory in the case that $G$ is a compact connected Lie group. Luecke twists by an $S^1$-central extension on the appropriate loop groupoid, which is equivalent to an $S^1$-gerbe \cite{Behren:Xu:DifferentiableStacksandGerbes}. As we are working with finite groups, Luecke's picture simplifies. In particular, Luecke's construction overcomes two issues which are not present in the finite case: choosing a loop groupoid with a strict $S^1$-action and choosing the appropriate notion of completed $S^1$-equivariant K-theory. 

It is well known that there is a strong relationship between equivariant elliptic cohomology, the K-theory of loop spaces, and the ordinary cohomology of double loop spaces \cite{spongphd}. In this context, it is expected that the introduction of transgression techniques, as discussed in this thesis, will play a natural role in connecting the twists of these theories.

\chapter{Orbifolds and Groupoids}\label{Groupoids}

\section{Groupoids}

In this section, we introduce the main definitions surrounding groupoids. In particular, we establish the notation that will be used throughout the thesis. The reader is referred to \cite{Morbifolds} for more details, as this is our main reference for this section.

\begin{definition}\label{def:groupoid}
A \emph{groupoid} $\G$ is a small category in which every morphism is an isomorphism. We write $\G_0$ for the set of objects and $\G_1$ for the set of morphisms, and we have the following structure maps:
\begin{itemize}
    \item The source and target maps $s,t\colon \G_1 \to \G_0$, which sends an arrow to its source/domain and target/codomain respectively.
    \item The unit map $u\colon \G_0 \to \G_1$, which maps an object to its identity arrow.
    \item The inverse map $i\colon \G_1 \to \G_1$, which sends an arrow to its inverse.
    \item The multiplication map $m\colon \G_1 \mathbin{_t \times_s} \G_1 \to \G_1$, which sends a pair of composable arrows $(g,h)$ to their composition $g\cdot h$.
\end{itemize}
These structure maps satisfy the expected conditions, for example $s \circ i = t$ and $i \circ u = u$. We will often write a groupoid as 
\[
\G = (\G_1 \rightrightarrows \G_0)
\quad \text{or} \quad
\G = 
\begin{array}{c}
\G_1 \\
\downdownarrows \\
\G_0
\end{array}
\]
where the two parallel arrows represent the source and target maps.
\end{definition}

\begin{remark}
From this point forward we will always use the letters $s,t,u,i,m$ for the structure maps of the groupoid, even when considering different groupoids. It will typically be clear in context which groupoid the structure maps belong to. 
\end{remark}

\begin{remark} 
Even though groupoids are categories, we rarely think of them in this way. We will see that groupoids are most often talked about as if they were groups or topological spaces. One can actually view groupoids as a generalisation of (i) groups, (ii) manifolds or (iii) equivalence relations. This will be seen in the coming examples. It is for this reason that if $g\colon x \to y$ and $h\colon y \to z$ are two arrows, then we write the composition as $g \cdot h \colon x \to z$ instead of $h \circ g\colon x \to z$. In other words, we treat the composition of arrows as a group multiplication rather than a composition of functions. 
\end{remark}

We will write $\G_n$ \label{Gn} for the space containing sequences of $n$ composable arrows in $\G$, that is,
\[
\G_n = \G_1 \mathbin{_t \times_s} \G_1 \mathbin{_t \times_s} ... \mathbin{_t \times_s} \G_1
\]
where the fibre product is repeated $n$ times. With this notation, the multiplication map is $m\colon \G_2 \to \G_1$. Define the maps $d_j\colon \G_n \to \G_{n-1}$ for $j \in \{0, ..., n\}$ which ``delete'' the $j$th object in the string of arrows. For instance, the maps $d_j\colon \G_2 \to \G_1$ are
\begin{align*}
   & d_0(x  \xrightarrow{g} y \xrightarrow{h} z) =  y \xrightarrow {h} z, \\
   & d_1(x  \xrightarrow{g} y \xrightarrow{h} z) =  x \xrightarrow {g \cdot h} z, \\
   & d_2(x  \xrightarrow{g} y \xrightarrow{h} z) = x \xrightarrow {g} y . 
 \end{align*}
These maps turn our groupoid into a simplicial set, called the \emph{nerve} of $\G$:
\begin{center}
\begin{tikzpicture}[scale=0.5]

    \foreach \n in {3, 2, 1}{
      \node at (12-4*\n, 0) {$\G_{\n}$};
     }
    \node at (12, 0) {$\G_{0}$.}; 
    
    \def \l {0.35}

    \foreach \n in {-1, 1} {
        \draw [->] (9, \n*\l/2) -- (11,\n*\l/2);
    }
    
    \foreach \n in {-1, 0, 1} {
        \draw [->] (5, \n*\l) -- (7,\n*\l);
    }

    \foreach \n in {0,1,2,3} {
        \draw [->] (1, \n*\l- 3*\l / 2) -- (3,\n*\l-3*\l/2);
    }

    \foreach \n in {-2, -1, 0, 1, 2} {
        \draw [->] (-3, \n*\l) -- (-1,\n*\l);
    }

\end{tikzpicture}
\end{center}
The geometric realisation of the nerve of $\G$ is called the classifying space $B\G$. 

We now introduce the groupoids that are most important to us in this thesis.

\begin{example}\label{ch2:actiongroupoid}
Let $X$ be a manifold acted on the right by a group $G$. The groupoid $X \mmod G$, called the \emph{action groupoid}, is the groupoid $X \times G \rightrightarrows X$ where 
\[
s(x,g) = x, \quad t(x,g) = x \cdot g, \quad \text{and} \quad m\bigl((x,g),(x\cdot g, h)\bigr) = x \cdot (gh).
\]
The classifying space of $X \mmod G$ is a model for the Borel construction $X_G$, see for instance \cite{Alejandro:Michele:lectures}. The importance of $X \mmod G$ will become very apparent throughout the thesis.
\end{example}

\begin{example}\label{ch2:group}
A group $G$ can be thought of as a groupoid with one object, which we write as $* \mmod G$ or $\mathbb{B}G$. The set of arrows is $G$ itself and the multiplication map is the group operation. This matches the perspective of groups as being symmetries of an object. One could then think of groupoids as being symmetries of \emph{many} objects.
\end{example}

\begin{example}\label{ch2:manifold}
Any topological space $M$ can be viewed as a groupoid $M \mmod * = (M \rightrightarrows M)$. In this case the only arrows are the identity arrows and all of the structure maps produce the relevant identity morphisms.
\end{example}

\begin{example}\label{ch2:opencovermanifold}
Let $M$ be a topological space and let $\U = \{U_i\}$ be an open cover of $M$. The open cover groupoid $M[\U]$ is defined as
\[
M[\U]_0 = \bigsqcup_{i} U_i = \{ (x, i) \mid x \in U_i \}, 
\]
\[
M[\U]_1 = \bigsqcup_{i,j} U_i \cap U_j = \{ (x, i, j) \mid x \in U_i \cap U_j \}.
\]
Here $(x,i,j)$ is a morphism from $(x,i)$ to $(x, j)$. The identity arrow of $(x,i)$ is $(x,i,i)$.
\end{example}

\begin{example}\label{ch2:opencovergroupoid}
Let $\G$ be a groupoid such that $\G_0$ is a topological space and let $\U = \{U_i\}$ be an open cover of $\G_0$. The open cover groupoid $\G[\U]$ is similar to the groupoid in the previous example:
\[
\G[\U]_0 = \bigsqcup_{i} U_i = \{ (x, i) \mid x \in U_i \},
\]
\[
\G[\U]_1 = \bigsqcup_{i,j} s^{-1}(U_i) \cap t^{-1}(U_{j}) = \{ (g,i, j) \mid s(g) \in U_i, \, t(g) \in U_j \}.
\]
This groupoid inherits the structure maps from $\G$, so that $s(g, i, j) = (s(g), j)$ and $t(g,i, j) = (t(g), j)$. In other words, a morphism from $(x,i)$ to $(y,j)$ is just a morphism $x \to y$ in $\G$.
\end{example}

\begin{example}\label{ch2:inertiagroupoid}
Let $\G$ be a groupoid. The \emph{inertia groupoid} $\Lambda \G$ is the groupoid
\[
(\Lambda \G)_0 = \{ g \in \G_1 \mid s(g) = t(g) \},
\]
\[
(\Lambda \G)_1 = \{ (g,h) \in \G_1 \times \G_1 \mid s(g) = t(g) = s(h) \},
\]
with $s(g,h) = g$ and $t(g,h) = h^{-1} g h$. So an object of $\Lambda \G$ is an automorphism of an object in $\G$ and an arrow is as in the following picture:


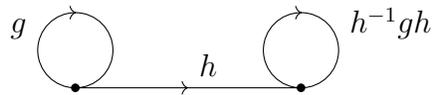
\begin{center}
\begin{tikzpicture}

\draw [fill] (-2,0) circle [radius=0.05];
\draw [fill] (1,0) circle [radius=0.05];

\draw [->] (-2,0) -- (-0.5,0);
\draw (-0.5,0) -- (1,0);

\draw (-2, 0.5) circle [radius = 0.5];
\draw (1, 0.5) circle [radius = 0.5];

\node [above right] at (-0.5,0) {$h$};
\node [above left] at (-2.5,0.5) {$g$};
\node [above right] at (1.5,0.5) {$h^{-1}gh$};

\draw [->] (-2.01, 1) -- (-1.99,1);
\draw [->] (0.99, 1) -- (1.01,1);

\end{tikzpicture}
\end{center}

\end{example}

This example perhaps highlights that it is not always obvious from $\G_1$ what the source and target maps are. One might expect on first glance that $(g,h) \in \Lambda\G_1$ is a morphism from $g$ to $h$, yet the structure maps indicate that this is not the case. In category theory, one thinks of a morphism as belonging to a ``hom-set'' of a particular source and target. In groupoids, we instead look at the space of all morphisms between all objects and rely on the structure maps to tell us how they behave as arrows. 

\begin{example}\label{ch2:equivalencerelation}
Let $X$ be a set with an equivalence relation $\sim$. The groupoid $X \mmod \!\sim$ has elements of $X$ as objects and there is an arrow from $x$ to $y$ if $x \sim y$. Conversely, given any groupoid $\G$ we can define an equivalence relation on $\G_0$ by setting $x \sim y$ if there is an arrow from $x$ to $y$ in $\G$. 
\end{example}

This last example leads us to define the quotient $\G_0 /\!\sim$ where $\sim$ is the equivalence relation determined by the arrows $\G_1$. This is called the \emph{orbit space} of $\G$ and will be denoted by $\lvert \G \rvert$ or $\G_0 / \G_1$. The reader can confirm the following:
\begin{itemize}
\item The orbit space of the action groupoid of $X \mmod G$ is the quotient space $X / G$.
\item The orbit space of $M[\U]$ is $M$.
\item The orbit space of $\G[\U]$ is the same as the orbit space of $\G$ itself.
\end{itemize}
Let $X$ be a topological space and $\G$ a groupoid such that $\G_0$ is a topological space and $\lvert\G\rvert$ is homeomorphic to $X$. In this case we will say that $\G$ \emph{represents} the space $X$. For example, both of the groupoids $S^1 \mmod *$ and $\R \mmod \Z$ represent the circle $S^1$.

The center of a group $G$ is the set of all $g \in G$ such that $gh = hg$ for all $h \in G$. This is generalised to groupoids in the following definition. We will make use of the center of a groupoid in Chapter 6. For now, we consider this as an example of how ideas in group theory can be generalised to groupoids.

\begin{definition}
The \emph{center} of a groupoid $\G$ is the group of natural transformations from the identity functor $\id_{\G}\colon \G \to \G$ to itself,
\[
\textnormal{Center}(\G) := \mathcal{N}\!at(\id_{\G}, \id_{\G}).
\]
Explicitly, an element $\xi \in \textnormal{Center}(\G)$ is, for each object $x$, an automorphism $\xi_x \colon x \to x$ such that if $g \colon x \to y$ is a morphism in $\G$, then $\xi_x \cdot g = g \cdot \xi_y$. 
\end{definition}

\begin{example}
The center of $\mathbb{B}G$ is the center of the group $G$.
\end{example}

\begin{example}\label{ex:centerofXgCg}
Let $X$ be acted on by a group $G$ and let $C_g$ be the centraliser of $g \in G$. Let $X^g$ denote the points of $X$ that are fixed by $g$,
\[
X^g = \{ x \in X \mid x\cdot g = x \}.
\]
Then $\xi_x = (x,g)$ defines an element in the center of $X^g \mmod C_g$. Indeed, if $(x,h) \in X^g \times C_g$ then 
\[
\xi_x \cdot (x,h) = (x,gh) = (x, hg) = (x,h) \cdot \xi_{x\cdot h}.
\]
This will play an important role in Chapter 6.
\end{example}


\section{Orbifold Groupoids}

In this section, we consider proper, \'etale Lie groupoids and homomorphisms between them.

\begin{definition}
A \emph{topological groupoid} is a groupoid $\G$ in which $\G_0$ and $\G_1$ are topological spaces and all of the structure maps are continuous.

A \emph{Lie groupoid} is a groupoid $\G$ in which $\G_0$ and $\G_1$ are smooth manifolds, all of the structure maps are smooth and the source and target maps are surjective submersions.

A Lie groupoid is \emph{\'etale} if the source and target maps are local diffeomorphisms.

A Lie groupoid is \emph{proper} if the map $(s,t)\colon \G_1 \to \G_0 \times \G_0$ is proper, that is, the pre-image of compact sets are compact.

An \emph{orbifold groupoid} is a proper, \'etale Lie groupoid.
\end{definition}

An orbifold groupoid is called such because it is the proper, \'etale Lie groupoids that represent orbifolds. This is made precise in Theorem \ref{thm:equivorbifoldsandgroupoids} at the end of the chapter. 

Since groupoids are categories we expect that morphisms between them will be functors. For Lie groupoids we also want to preserve the smooth structure:

\begin{definition}
A  \emph{homomorphism} $\phi\colon \G \to \H$ between Lie groupoids is a smooth functor, that is, a pair of smooth maps $\phi_1\colon \G_1 \to \H_1$ and $\phi_0\colon \G_0 \to \H_0$ that commute with the structure maps.
\end{definition}

There is also the notion of 2-morphisms between Lie groupoids, which are just smooth natural transformations. Explicitly, a natural transformation $\eta\colon \phi \Rightarrow \psi$ between smooth functors $\phi, \psi\colon \G \to \H$ is a smooth map $\eta\colon \G_0 \to \H_1$ such that $\eta(x)$ is a morphism from  $\phi_0(x)$ to $\psi_0(x)$ and if $g\colon x \to y$ is an arrow in $\G$, then the following diagram in $\H$ commutes:
\[
\begin{tikzcd}
\phi_0(x) \arrow[r, "\eta(x)"] \arrow[d, "\phi_1(g)", swap] & \psi_0(x) \arrow[d, "\psi_1(g)"] \\
\phi_0(y) \arrow[r, "\eta(y)"]           & \psi_0(y)         
\end{tikzcd}
\]

Between categories there is the notion of an equivalence of categories. We want to have the same notion for Lie groupoids. This leads to a smooth formulation of an equivalence of categories:

\begin{definition} 
An \emph{equivalence} $\phi\colon \G \to \H$ of Lie groupoids is a homomorphism satisfying the following conditions:
\begin{enumerate}[label=(\roman*)]
    \item The composition
    \[
    \H_1 \mathbin{_s\times_{\phi_0}} \G_0 \xrightarrow{\pr_1} \H_1 \xrightarrow{\,\, t \,\,} \H_0
    \]
    is a surjective submersion.
    \item The following diagram is a pullback square:
    \begin{center}
    \begin{tikzcd}
            \G_1 \arrow[r, "\phi_1"] \arrow[d, "{(s,t)}", swap] & \H_1 \arrow[d, "{(s,t)}"] \\
            \G_0 \times \G_0 \arrow[r, "{\phi_0\times\phi_0}"]           & \H_0 \times \H_0          
    \end{tikzcd}
\end{center}
\end{enumerate}
\end{definition}

\noindent One may think of these conditions as the smooth formulation of essential surjectivity and fully faithfulness, respectively. 

\begin{remark}
Some authors require that the composition
\[
\H_1 \mathbin{_t \times_\phi} \G_0 \xrightarrow{\text{pr}_1} \H_1 \xrightarrow{\,\,s\,\,} \H_0
\]
is a surjective submersion instead of the map in condition (i). These are equivalent because if $(h,x) \in \H_1 \mathbin{_t\times_\phi} \G_0$ then $(h^{-1},x) \in \H_1 \mathbin{_s \times_\phi}\G_0$,
\[
t \circ \text{pr}_1(h,x) = t(h) = s(h^{-1}) = s \circ \text{pr}_1(h^{-1}, x),
\]
and the inverse map is a diffeomorphism. For condition (ii) we can equivalently require that the map
\[
\G_1 \to (\G_0 \times \G_0) \mathbin{_{(\phi,\phi)}\times_{(s,t)}} \H_1, 
\quad
g \mapsto (s(g), t(g), \phi_1(g))
\]
is a diffeomorphism. 
\end{remark}

\begin{example}\label{ch2:circleequiv}
The homomorphism $\phi \colon \R \mmod \Z \to S^1 \mmod *$ defined on objects by $\phi_0(x) = e^{2\pi i t}$ and on arrows by $\phi_1(t,n) = e^{2\pi i n t}$ is an equivalence. 
\end{example}

\begin{example}\label{ch2:opencoverequivalence}
Let $M$ be a smooth manifold with an open cover $\,\U = \{U_i\}$. The homomorphism
\[
\varepsilon\colon M[\U] \to M \mmod *
\]
where $\varepsilon_0(x,i) = x$ and $\varepsilon_1(x, i,j) = x$ is an equivalence.
\end{example}

\begin{example}\label{ch2:groupoidopencoverequivalence}
Let $\G$ be a Lie groupoid with $\U = \{U_i\}$ an open cover of $\G_0$. The homomorphism
\[
\varepsilon\colon \G[\U] \to \G
\]
where $\varepsilon_0(x, i) = x$ and $\varepsilon_1(g, i,j) = g$ is an equivalence.
\end{example}

\begin{lemma}
An equivalence $\phi\colon \G \to \H$ induces a homeomorphism
\[
\lvert \phi \rvert \colon \lvert \G \rvert \to \vert \H \rvert, \quad [x] \mapsto [\phi_0(x)].
\]
\end{lemma}
\begin{proof}
We construct an inverse function. Using the first condition we can find an open cover $\{U_i\}$ of $\H_0$ and local sections of $t \circ \pr_1$,
\[
\sigma^i\colon U_i \to \H_1 \fp{s}{\phi_0} \G_0,
\]
which we denote by $\sigma^i(y) = (\sigma^i_1(y), \sigma^i_2(y))$. Define the map
\[
\psi\colon \vert \H \rvert \to \vert\G\rvert, \quad [y] \mapsto [\sigma^i_2(y)] \text{ where } y \in U_i.
\]
This map is well-defined: if $y \in U_i \cap U_j$ then $\sigma_1^i(y) \cdot \sigma_1^j(y)^{-1}$ is a morphism from $\phi_0 (\sigma_2^i(y))$ to $\phi_0 (\sigma_2^j(y))$ in $\H$. Then the second condition on $\phi$ implies that there is a morphism from $\sigma_2^i(y)$ to $\sigma_2^j(y)$ in $\G$, making them equal in the orbit space. A similar argument shows that $\psi$ is a left inverse to $\lvert \phi \rvert$. To show that this is a right inverse, note that for $y \in U_i$ we have
\[
[y] \xmapsto{\, \psi\,} [\sigma^i_2(y)] \xmapsto{\, \lvert \phi \rvert\,} [\phi_0(\sigma^i_2(y))] = [s\sigma^i_1(y)] = [t \sigma_1^i(y)] = [y].
\]
Therefore, $\psi$ is a right inverse to $\lvert \phi \rvert$ as well.
\end{proof}

\begin{lemma}
If $\phi\colon \G \to \H$ and $\psi\colon \H \to \K$ are equivalences then the composition $\psi \circ \phi$ is an equivalence.
\end{lemma}
\begin{proof}
Consider the following commutative diagram:
\[
\begin{tikzcd}
(\K_1 \fp{s}{\psi_0} \H_0 ) \fp{\text{pr}_2}{\phi_0} \G_0 \arrow[d, "\cong"]  \arrow[r, ] &  \K_1 \fp{s}{\psi_0} \H_0 \arrow[d, "t \text{pr}_1"]  \\
\K_1 \fp{s}{\psi_0 \phi_0} \G_0   \arrow[r, "t \text{pr}_1"]   & \K_0          
\end{tikzcd}
\]
The left vertical map is defined by $((k, x), y) \mapsto (k, y)$ and is a diffeomorphism with inverse $(k,y) \mapsto ((k, \phi_0(y)), y)$. The top horizontal map is a projection from a fibre product, which is a surjective submersion. The right vertical map is a surjective submersion since $\psi$ is an equivalence. Therefore the bottom map is a surjective submersion, as required. 

For the second condition we consider the following diagram:
\[
\begin{tikzcd}
\G_1 \arrow[d, "{(s,t)}"] \arrow[r, "\phi_1"] & \H_1 \arrow[r, "\psi_1"] \arrow[d, "{(s,t)}"] & \K_1 \arrow[d, "{(s,t)}"] \\
\G_0 \times \G_0  \arrow[r, "\phi_0 \times \phi_0"] & \H_0 \times \H_0 \arrow[r, "\psi_0 \times \psi_0"]           & \K_0 \times \K_0       
\end{tikzcd}
\]
By hypothesis the left and right squares are pullback diagrams, and so the outer square is a pullback diagram as well.
\end{proof}

When we have an equivalence of categories, there always exists an ``inverse'' functor that is also an equivalence of categories. The problem is that the inverse of an equivalence may not be smooth. For example, the equivalence in Example \ref{ch2:opencoverequivalence} does not have a smooth inverse unless the open cover is trivial, because there is no smooth map $M \to \bigsqcup_{\alpha} U_{\alpha}$. Ideally, we would like a category of Lie groupoids in which morphisms are homomorphisms and isomorphisms are equivalences. However, the inverse of an equivalence may not even be a homomorphism. The solution is to instead define generalised maps between Lie groupoids:

\begin{definition}
A \emph{generalised map} $\G \to \H$ between Lie groupoids $\G$ and $\H$ is an equivalence class of triples $(\G', \varepsilon, \phi)$ fitting into the diagram,
\[
\G \xleftarrow{\,\, \varepsilon \,\,} \G' \xrightarrow{\,\,\phi\,\,} \H
\]
where $\varepsilon$ is an equivalence. The triple $(\G', \varepsilon, \phi)$ is equivalent to $(\G'', \varepsilon', \phi')$ if there exists a homomorphism $\gamma\colon \G' \to \G''$ making the following diagram commutes up to natural transformation:
\begin{center}
\begin{tikzpicture}
    \node (L) at (0,1)  {$\G$};
    \node (R) at (5,1) {$\H$};
    \node (U) at (2.5, 1.75) {$\G'$};
    \node (D) at (2.5, 0.25) {$\G''$};
    \draw [<-] (L) -- (U) node [midway, above] {\small$\varepsilon$};
    \draw [->] (U) -- (R) node [midway, above] {\small$\phi$};
    \draw [<-] (L) -- (D) node [midway, below] {\small$\varepsilon'$};
    \draw [->] (D) -- (R) node [midway, below] {\small$\phi'$};
    \draw [->] (U) -- (D) node [midway, right] {\small$\gamma$};
\end{tikzpicture}
\end{center}
\end{definition}

\begin{remark}
An alternative notion of maps between orbifold groupoids is that of \emph{bibundles} and \emph{Hilsum-Skandalis maps}. In this approach, a map between Lie groupoids $\G$ and $\H$ is a principal $\H$-bundle on $\G_0$ with an additional action by $\G$. See Appendix \ref{appendix:bibundles} for a discussion of these ideas.
\end{remark}

\noindent The reader should think of the groupoid $\G'$ as an open cover of $\G$, in fact we have the following:

\begin{proposition}
Any generalised map $\G \to \H$ can be represented by an open cover $\mathcal{U}$ of $\G_0$ and a homomorphism $\phi$ so that we have
\[
\G \xleftarrow{\,\,\varepsilon\,\,} \G[\U] \xrightarrow{\,\,\phi\,\,} \H.
\]
where $\varepsilon$ is defined in Example \ref{ch2:groupoidopencoverequivalence}. Furthermore, another such triple $(\G[\U'], \varepsilon', \phi')$ represents the same generalised map if and only if $\phi$ and $\phi'$ are related by a natural transformation when restricted to a common refinement of $\mathcal{U}$ and $\mathcal{U}'$.
\end{proposition}

This result is best proven using the language of bibundles, where there is a natural notion of local triviality. A proof is provided in Appendix \ref{appendix:bibundles}.

The suitable notion of equivalence between Lie groupoids is that of Morita equivalence:

\begin{definition}
A \emph{Morita equivalence} between $\G$ and $\H$ is a generalised map $(\H', \varepsilon, \phi)$ such that $\phi$ is an equivalence. If a Morita equivalence exists, then $\G$ and $\H$ are said to be \emph{Morita equivalent}.
\end{definition}

If one were to discuss orbifolds without the language of groupoids, they would speak of topological spaces equipped with an equivalence class of orbifold atlases. This is analogous to smooth manifolds equipped with smooth atlases. In a similar method to how we constructed a groupoid from an open cover of a manifold, one can build a groupoid using an orbifold atlas. Such a groupoid will be a proper, \'etale Lie groupoid. If a different, yet compatible, orbifold atlas was used, then the groupoid obtained would be Morita equivalent to the original. In this sense, when we think of a groupoid $\G$ representing an orbifold, $\G$ plays the role of an atlas on the topological space $\vert \G \rvert$. For more details on this, the original reference is \cite{MP97} but the reader should also see \cite{Morbifolds} and \cite{amenta}. The important theorem here is the following \cite[Thm 4.1 (4) $\Rightarrow$ (1)]{MP97}: 

\begin{theorem}\label{thm:equivorbifoldsandgroupoids}
There is an equivalence of categories between the category of orbifolds and the category of orbifold groupoids with generalised maps.
\end{theorem}

\begin{remark}
There is a formal way of forcing a certain class of maps in a category to be invertible. Let $\CC$ be a category and $W$ a subclass of morphisms in $\CC$. A (the) localisation of $\CC$ with respect to $W$ is a category $\CC[W^{-1}]$ with a functor $L\colon \CC \to \CC[W^{-1}]$ such that:
\begin{enumerate}[label=(\roman*)]
    \item For any $w \in W$ the arrow $L(w)$ is invertible in $\CC[W^{-1}]$.
    \item If $\phi\colon \CC \to D$ is any functor such that $\phi(w)$ is invertible for every $w \in W$, then there is a unique functor $\psi\colon \CC[W^{-1}] \to D$ such that $\psi \circ L = \phi$.
\end{enumerate}
This definition is from \cite[\textsection 3.1]{Lerman}. The category of orbifold groupoids with generalised morphisms is actually the category $\operatorname{Gp}[W^{-1}]$ where $\operatorname{Gp}$ is the category of orbifold groupoids with isomorphism classes of smooth functors and $W$ is the collection of arrows represented by equivalences.
\end{remark}

\chapter{Loop Groupoids}\label{LoopGroupoids}


\section{Loop Groupoids}

The free loop space $LM$ of a manifold $M$ is the space of all loops $S^1 \to M$. We want to generalise this construction to orbifolds, defining a loop groupoid that represents the loop space of the orbifold under consideration. The reference for these ideas is \cite{LU02}, in which Lupercio and Uribe formulate the loop groupoid and investigate its properties. This section aims to provide a complementary account of the material presented in their paper

Given a Lie groupoid $\G$, one expects the objects of the loop groupoid $\L \G$ to be orbifold maps $\S^1 \to \G$, where $\S^1$ is a groupoid representing the circle. The canonical choice of $\S^1$ is 
\[
\S^1 := S^1 \mmod * = \bigl( S^1 \rightrightarrows S^1\bigr).
\]
An orbifold map from $\S^1$ to $\G$ is an equivalence class of diagrams of the form
\[
\S^1 \leftarrow \S^1[\U] \to \G
\]
where $\U$ is an open cover of $S^1$.  Pulling $\U$ back along the exponential map $t \mapsto e^{2\pi i t}$, we can produce an open cover $\V$ of $\R$. The groupoid $\R \mmod \Z [\V]$ fits into the diagram
\begin{center}
\begin{tikzpicture}
    \node (L) at (0,1)  {$\S^1$};
    \node (R) at (5,1) {$\G$};
    \node (U) at (2.5, 1.75) {$\S^1 [\U]$};
    \node (D) at (2.5, 0.25) {$\R \mmod \Z [\V]$};
    \draw [<-] (L) -- (U);
    \draw [->] (U) -- (R);
    \draw [<-] (L) -- (D);
    \draw [->] (D) -- (R);
    \draw [->] (D) -- (U);
\end{tikzpicture}
\end{center}
 where the vertical map is an equivalence induced by the exponential map and the bottom two maps are chosen to make the diagram commute. This means that the generalised maps 
\[
\S^1 \leftarrow \S^1 [\U] \to \G
\quad \text{and} \quad
\S^1 \leftarrow \R \mmod \Z [\V] \to \G
\]
are equivalent. Introducing the notation
\[\label{S1U}
\S^1_{\U} := \R \mmod \Z [\V],
\]
where $\U$ is an open cover of $S^1$ and $\V$ is the corresponding open cover of $\R$, we are now justified in restricting our attention to generalised maps of the form
\[
\S^1 \leftarrow \S^1_{\U} \to \G.
\]
The reason for switching from $\S^1[\U]$ to $\S^1_{\U}$ is because it is more convenient to map out of $\R$ than to map out of $S^1$. It is analogous to working with $\R / \Z$ instead of $S^1$ or by considering loops that are maps out of $[0,1]$ rather than $S^1$.

We now ask which open covers of the circle it is sufficient to consider. Immediately, we can restrict to finite covers: compactness of $S^1$ implies that any open cover $\U$ has a finite sub-cover $\U'$. The inclusion $\S^1_{\U'} \hookrightarrow \S^1_{\U}$ is an equivalence and fits into the diagram
\begin{center}
\begin{tikzpicture}
    \node (L) at (0,1)  {$\S^1$};
    \node (R) at (5,1) {$\G$};
    \node (U) at (2.5, 1.75) {$\S^1_{\U}$};
    \node (D) at (2.5, 0.25) {$\S^1_{\U'}$};
    \draw [<-] (L) -- (U);
    \draw [->] (U) -- (R);
    \draw [<-] (L) -- (D);
    \draw [->] (D) -- (R);
    \draw [->] (D) -- (U);
\end{tikzpicture}
\end{center}
where the bottom arrows are chosen to make the diagram commute. The top and bottom rows hence give the same generalised map, which means that we can always represent a generalised map with a finite cover. This argument, which was applied to sub-covers of $\U$, will also work for refinements of $\U$. The complexity of the orbifold maps $\S^1 \to \G$ will depend on the complexity of the open covers we choose, so we want these to be as nice as possible. The question is: what is a suitably ``nice'' class of open covers of $S^1$ that any finite cover can be refined into?

\begin{figure}
    \centering
    \begin{tikzpicture}

        \def \ptsize {0.05};
        \def \e {0.4};
        \def \angle {10};

        \begin{scope}[shift={(3,2)}]

            \def \a {-2}; \def \b {-0.7}; \def \c {0.6}; \def \d {2};
            
            \draw [ultra thick, myblue, (-)] (\a-\e, 0) -- (\b+\e, 0) node [black, midway, above] {$I_1$};
            \draw [ultra thick, orange, (-)] (\b-\e, 0) -- (\c+\e, 0) node [black, midway, above] {$I_2$};
            \draw [ultra thick, mygreen, (-)] (\c-\e, 0) -- (\d+\e, 0) node [black, midway, above] {$I_3$};

            \draw [fill] (\a, 0) circle [radius=\ptsize] node [below, yshift=-2.5pt] {$\alpha_0$};
            \draw [fill] (\b, 0) circle [radius=\ptsize] node [below, yshift=-2.5pt] {$\alpha_1$};
            \draw [fill] (\c, 0) circle [radius=\ptsize] node [below, yshift=-2.5pt] {$\alpha_2$};
            \draw [fill] (\d, 0) circle [radius=\ptsize] node [below] {$\alpha_0+1$};
            
            \def \down {-0.7}
            \draw [(-)] (\b - \e, \down) -- (\b+\e, \down) node [midway, below] {$2\epsilon$};
            \draw [(-)] (\c - \e, \down) -- (\c+\e, \down) node [midway, below] {$2\epsilon$};

        \end{scope}
        
        \begin{scope}[shift={(11,2)}]
        
            \draw (0,0) circle [radius=2];
            
            \draw [ultra thick, mygreen, -)] (0:2) arc (0:40:2);
            \draw [ultra thick, myblue, (-)] (20:2) arc (20:180:2);
            \draw [ultra thick, orange, (-)] (160:2) arc (160:280:2);
            \draw [ultra thick, mygreen, (-] (260:2) arc (260:360:2);
            
            \draw [fill] (30:2) circle [radius = \ptsize] node [above right, yshift=-2pt] {$\alpha_0$};
            \draw [fill] (170:2) circle [radius = \ptsize] node [above left, yshift=-3.5pt] {$\alpha_1$};
            \draw [fill] (270:2) circle [radius = \ptsize] node [below] {$\alpha_2$};
        
        \end{scope}
        
        \draw [thick, ->] (6.5, 2) -- (7.5,2) node [midway, above] {$\exp$};

\end{tikzpicture}
    \caption{An admissible cover of $S^1$.}
    \label{fig:admissiblecover}
\end{figure}
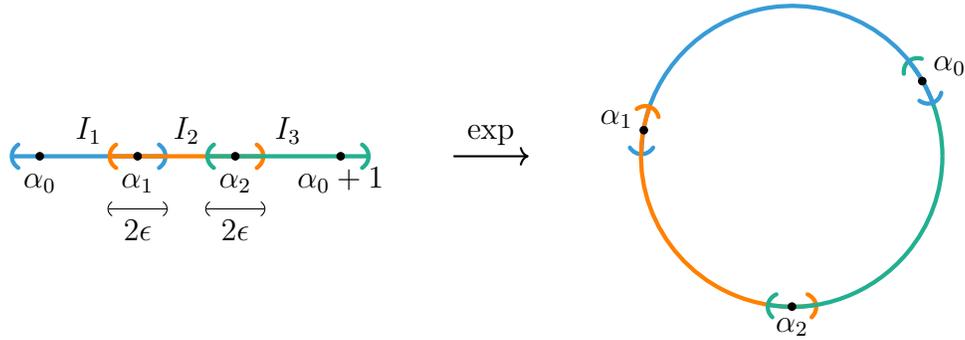

Lupercio and Uribe's notion of \emph{admissible} covers provide an answer to this question. Consider a finite collection of points
\[
0 < \alpha_0 < \alpha_1 < ... < \alpha_n  \leq  1
\]
and define the intervals
\[
I_i = (\alpha_{i-1} - \epsilon, \alpha_i + \epsilon), \quad i \in \{1, ..., n+1\},
\]
where:
\begin{itemize}
    \item $\alpha_{n+1} := \alpha_0+1$.
    \item $\epsilon > 0$ is chosen small enough that triple intersections of the $I_i$ are empty after mapping to the circle via the exponential map $t \mapsto e^{2\pi i t}$.
\end{itemize}
Under the exponential map, each interval maps to an open set in $S^1$ and together form an open cover of $S^1$, see Figure \ref{fig:admissiblecover}. A cover of $S^1$ obtained in this way is called an \emph{admissible cover}. We note a few properties:
\begin{itemize}
\item Choosing a smaller $\epsilon$ results in a refinement of the original admissible cover, and hence for our purposes is the same. We will therefore define an admissible cover using its ``vertices'' $\alpha_i$ (see Remark \ref{rem:admissiblecover}) and assume that $\epsilon$ and the intervals $I_i$ are implicitly defined. 
\item Given any two admissible covers, there is a common refinement which is also an admissible cover. This is obtained, for instance, by taking the union of all the chosen vertices with a sufficiently small $\epsilon$. 
\item Any finite cover of $S^1$ can be refined to an admissible cover. This amounts to removing any ``redundant'' open sets, breaking the remaining sets up so that they're connected, and then shrinking them until they're the correct length. Each of these steps can be achieved by a suitable refinement of the initial open cover.
\end{itemize}

\begin{remark}\label{rem:admissiblecover}
One may think of an admissible cover as a triangulation of $S^1$ that has been converted into an open cover by extending the end of each edge by $\epsilon$. A common refinement of two admissible covers is then a common triangulation with a sufficiently small $\epsilon$. As the $\epsilon$ may be chosen arbitrarily small, in application we will often think of admissible covers and triangulations as being one and the same. This is the reason we refer to the $\alpha_i$ as ``vertices'' of the admissible cover.
\end{remark}

Let $\U$ be an admissible cover of $S^1$. If $\U$ is the trivial cover of $S^1$, which occurs when $\U$ is obtained from a single vertex, then $\S^1_{\U} = \R \mmod \Z$.  Otherwise, an explicit description of the groupoid $\S^1_{\U}$ is:
\[
\S^1_{\U} =
\left(
\begin{array}{c}
    \bigsqcup_{i,j} (W_i \cap W_j)  \times \Z \\
     \downdownarrows \\
     \bigsqcup_i W_i
\end{array}
\right),
\quad \text{where} \quad
W_i = \bigsqcup_{n \in \Z} (I_i + n).
\] 
The source and target maps are
\[
s(x, n, i,j) = (x,i) \quad \text{and} \quad t(x,n,i,j) = (x+n,j).
\]
If $\U$ and $\V$ are two admissible covers, then their common refinement will be denoted by $\U \cap \V$. We have that
\[
\S^1_{\U \cap \V} = 
\left(
\begin{array}{c}
    \bigsqcup_{i,j, k, l} (W_i \cap W_j' \cap W_k \cap W_l' )  \times \Z \\
     \downdownarrows \\
     \bigsqcup_{i,j} W_i \cap W_j'
\end{array}
\right)
\]
where the $W_j'$ are obtained from $\V$ in the same way that the $W_i$ are obtained from $\U$. This may appear to be a large disjoint union but in fact almost all of the $W_i \cap W_j'$ are empty.

With a suitably nice family of open covers established, we turn to defining the loop groupoid. As discussed, the objects of the loop groupoid of $\G$ should be orbifold maps
\begin{equation}\label{genmap}
\S^1 \leftarrow \S^1_{\U} \to \G
\end{equation}
where $\U$ ranges across all admissible covers and their refinements. This suggests a colimit construction; the loop groupoid $\L\G(\U)$ will be defined for a particular admissible cover $\U$ and the complete loop groupoid $\L\G$ will be a colimit of $\L\G(\U)$ over all of the admissible covers $\U$.

Fix an admissible cover $\U$. We construct a groupoid in which the objects are orbifold maps of the form \eqref{genmap}. The left homomorphism $\S^1_{\U} \to \S^1$ is fixed by the choice of $\U$, so the objects will be smooth functors $\S^1_{\U} \to \G$. If the objects are to be functors, the morphisms should be natural transformations between these functors. Given two smooth functors $\Phi, \Psi \colon \S^1_{\U} \to \G$, a natural transformation $\Lambda \colon \Phi \to \Psi$ is a map
\[
\Lambda \colon (\S^1_{\U})_0 \to \G_1
\]
such that if $x \in (\S^1_{\U})_0$, then $\Lambda(x)$ is a morphism from $\Phi_0(x)$ to $\Psi_0(x)$. This must satisfy the following naturality condition: if $f \colon x \to y$ is a morphism in $\S^1_{\U}$, then the following diagram in $\G$ commutes:
\[ 
\begin{tikzcd} \Phi(x) \arrow[r, "{\Lambda(x)}"] \arrow[d, "{\Phi_1(f)}", swap] &       \Psi(x) \arrow[d, "{\Psi_1(f)}"] \\
    \Phi(y) \arrow[r, "{\Lambda(y)}"]           &  \Psi(y). 
\end{tikzcd}
\]  
This discussion is summarised in the following definition \cite[Def 3.2.1]{LU02}:
\begin{definition}\label{def:LG(U)}
Let $\G$ be a Lie groupoid and $\U$ an admissible cover. The \emph{loop groupoid of $\G$ associated to the open cover $\U$}, denoted by $\L\G(\U)$, is defined as follows:
\begin{itemize}
    \item Objects are smooth functors $\Phi \colon \S^1_{\U} \to \G$.
    \item Morphisms are smooth natural transformations $\Lambda \colon \Phi \to \Psi$.
    \item Multiplication of morphisms is given by the composition of arrows in $\G$. Concretely, if $\Lambda \colon \Phi \to \Psi$ and $\Omega \colon \Psi \to \Theta$ are morphisms, then $\Lambda \cdot \Omega$ is given by the map
    \[
    \Lambda \cdot \Omega  \colon (\S^1_{\U})_0 \to \G_1, \quad (x,i) \mapsto \Lambda(x,i) \cdot \Omega(x,i). 
    \]
\end{itemize}
Objects of $\L\G(\U)$ are called loops associated to the open cover $\U$.
\end{definition}

\begin{remark}
Lupercio and Uribe define morphisms in $\L\G(\U)$ as maps $\Lambda \colon (\S^1_{\U})_1 \to \G_1$ such that $\Lambda(x,n,i,j)$ is a morphism from $\Phi(x,i)$ to $\Psi(x+n,j)$ and
\[
\Lambda(x,n,i,j) = \Phi_1(x,n,i,i) \cdot \Lambda(x+n,0,i,j) = \Lambda(x,0,i,j) \cdot \Psi(x,n,j,j). 
\]
In other words, the naturality condition is used to extend from a map on $(\S^1_{\U})_0$ to a map on $(\S^1_{\U})_1$. The two definitions are equivalent.
\end{remark}

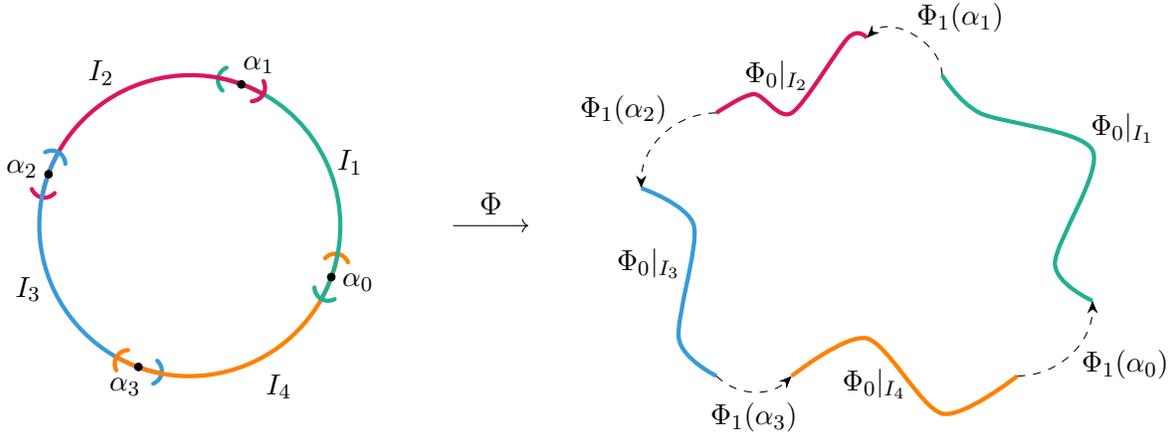
\begin{figure}
    \centering
    \begin{tikzpicture}
    
    \def \ptsize {0.05};
    
    \begin{scope}[shift={(3,2)}]
    
        \draw (0,0) circle [radius=2];
        
        \draw[ultra thick, mygreen, -)] (0:2) arc (0:80:2);
        \draw [ultra thick, mypink, (-)] (60:2) arc (60:170:2) node [black, midway, above left, yshift=-3pt] {\small$I_2$};
        \draw [ultra thick, myblue, (-)] (150:2) arc (150:260:2) node [black, midway, left] {\small$I_3$};
        \draw [ultra thick, orange, (-)] (240:2) arc (240:350:2) node [black, midway, below right] {\small$I_4$};
        \draw [ultra thick, mygreen, (-] (330:2) arc (330:360:2);
        \node [right] at (25:2) {\small$I_1$};
        
        \draw [fill] (70:2) circle [radius=\ptsize] node [above right, xshift=-3.5pt] {\small$\alpha_1$};
        \draw [fill] (160:2) circle [radius=\ptsize] node [left, yshift=2pt] {\small$\alpha_2$};
        \draw [fill] (250:2) circle [radius=\ptsize] node [below left, xshift=5pt] {\small$\alpha_3$};
        \draw [fill] (340:2) circle [radius=\ptsize] node [right, yshift=-2pt]  {\small$\alpha_0$};
    
    \end{scope}
    
    \begin{scope}[shift={(12,2)}]
    
        \draw [ultra thick, myblue] plot [smooth] coordinates {(-3, 0.5) (-2.3,0)  (-2.5,-1.5) (-2,-2)}; 
        \draw [ultra thick, orange] plot [smooth] coordinates {(-1,-2) (0,-1.5) (1,-2.5) (2,-2)};
        \draw [ultra thick, mygreen] plot [smooth] coordinates {(3,-1) (2.5,-0.5) (3,1) (1.5,1.5) (1,2)};
        \draw [ultra thick, mypink] plot [smooth] coordinates {(0,2.5) (-0.25,2.5) (-1,1.5) (-1.5, 1.75) (-2,1.5)};
        
        \draw [dashed, -{Stealth[scale=1,angle'=45]}] (-2,1.5) to[out=180,in=90] node [midway, above left, xshift=2mm] {\small$\Phi_1(\alpha_2)$} (-3,0.5) ;
        \draw [dashed, -{Stealth[scale=1,angle'=45]}] (-2,-2) to[out=-45, in=-135] node [midway, below] {\small$\Phi_1(\alpha_3)$} (-1,-2);
        \draw [dashed, -{Stealth[scale=1,angle'=45]}] (2,-2) to[out=0, in=-90] node [midway, below right, yshift=2mm] {\small$\Phi_1(\alpha_0)$} (3,-1);
        \draw [dashed, -{Stealth[scale=1,angle'=45]}] (1,2) to[out=90, in=45] node [midway, above right, yshift=-1mm, xshift=-1mm] {\small$\Phi_1(\alpha_1)$} (0,2.5);     
    
        \node at (-2.9,-0.5) {\small $\Phi_0|_{I_3}$};
        \node at (0.1,-2.1) {\small $\Phi_0|_{I_4}$};
        \node at (3.4, 1.3) {\small $\Phi_0|_{I_1}$};
        \node at (-1.2,2) {\small $\Phi_0|_{I_2}$};
    
    \end{scope}
    
    \draw [->] (6.5, 2) -- (7.5, 2) node [above, midway] {$\Phi$};
    
\end{tikzpicture}
    \caption{A loop in $\L\G(\U)$ is a disjoint collection of paths in $\G_0$ with morphisms connecting each path to the next.}
    \label{fig:orbifoldloop}
\end{figure}

An object $\Phi \in \L\G(\U)_0$ is depicted schematically in Figure \ref{fig:orbifoldloop}. Since $(\S^1_{\U})_0$ is a disjoint union of intervals, $\Phi_0$ is a disjoint collection of paths in $\G_0$. The images $\Phi_0(I_i)$ and $\Phi_0(I_i + n)$ are identified by morphisms of the form $\Phi_1(x, x+n, i,i)$. This means that each $\Phi_0(W_i)$ is a collection of ``copies'' of a single path $\Phi_0(I_i)$. The remaining morphisms, those between $W_i$ and $W_j$ for $i \neq j$, are mapped by $\Phi_1$ to morphisms in $\G$ that identify the ends of the $\Phi_0(I_i)$ and $\Phi_0(I_j)$ to form a complete loop when isomorphic objects are identified. Note that in Figure \ref{fig:orbifoldloop} we only bother to draw the images of $I_1, ...., I_n$ and omit the images of their translates.

Two morphisms $\Phi, \Psi \colon \S^1_{\U} \to \G$ are equivalent if there is a common refinement $\U \cap \V$ of $\U$ and $\V$ such that the following diagram of functors commutes up to natural transformation:
\[
\begin{tikzcd}
\S^1_{\U \cap \V} \arrow{r} \arrow{d} & \S^1_{\U} \arrow[d, "\Phi"] \\
\S^1_{\V} \arrow[r, "\Psi",swap]           & \G          
\end{tikzcd}
\]
The natural way of considering all the loops associated to all the admissible open covers under this equivalence is by using a colimit construction.

\begin{definition}
Let $\G$ be a topological groupoid. The \emph{loop groupoid} $\L\G$ is defined to be
\[
\L\G = \varinjlim_{\U} \L\G(\U)
\]
where the colimit ranges over the admissible covers of $S^1$.
\end{definition}

\begin{example}\label{ex:trivialloop}
Our favourite loops will be those associated to the trivial cover. As discussed, if $\U$ is the trivial cover of $S^1$ then $S^1_{\U} = \R \mmod \Z$. Loops associated with this cover are smooth functors
\[
\Phi \colon \R \mmod \Z \to \G.
\]
These loops are fully determined by the value of $\Phi_1$ on $\R \times \{1\}$, because
\begin{gather*}
    \Phi_0(x) = s(\Phi_1(x,1)), \\
    \Phi_0(x+1) = t(\Phi_1(x,1)), \\
    \Phi_1(x,0) = u( \Phi_0(x)), \\
    \Phi_1(x,n) = \Phi_1(x,1) \Phi_1(x+1, 1) \dotsm \Phi_1(x+n-1,1),  \quad \text{and}\\
    \Phi_1(x,-n) = i(\Phi_1(x-n ,n)),
\end{gather*}
where $n$ is a positive integer. Moreover, the restriction of $\Phi_0$ to an interval $[x,x+1]$ descends to a well-defined loop in $\vert \G \rvert$, see Figure \ref{fig:loop1}. 
\begin{figure}[h]
    \centering
    \begin{tikzpicture}

\def \pt {0.05}

\draw [ultra thick, myblue] plot [smooth] coordinates {(0,0) (0.8, -0.5) (0,-1) (0.8, -1.5) (1.5, -1) (2.5, -0.3) (3.5, -1) (4, 0)};

\draw [dashed] (0,0) .. controls (0,1) and (4,1) .. (4,0);

\draw [->] (1.99,0.74) -- (2.01, 0.74);
\draw [myblue, ->] (2-0.1,-0.57-0.1) -- (2+0.1,-0.572+0.1);

\draw [fill] (0,0) circle [radius = \pt] node [left=2pt] {$\Phi_0(0)$};
\draw [fill] (4,0) circle [radius = \pt] node [right=2pt] {$\Phi_0(1)$};
\node [above=2pt] at (2, 0.74) {$\Phi_1(0)$};
\node [below right=2pt] at (2-0.5, -0.57-0.1) {$\Phi_0$};

\begin{scope}[xshift = 0.65]
    
    \draw [ultra thick, myblue] plot [smooth cycle] coordinates {(9,0) (10, -0.5) (9.6, -1) (11, -1.3) (12, -0.5) (13,0) (12.6, 0.5) (11,0.4) (10, 0.8) };

    \draw [fill] (9,0) circle [radius = \pt];
    \node [left] at (9+0.4,0-0.2) {\begin{tabular}{c}$[\Phi_0(0)]$ \\  $=[\Phi_0(1)]$ \end{tabular} };
    
\end{scope}

\draw [->] (5.5, 0.5) -- (7, 0.5) node [above, midway] {quotient};

\end{tikzpicture}
    \caption{$\Phi_0$ descends to a loop in $\lvert \G \rvert$.}
    \label{fig:loop1}
\end{figure}
\end{example}


\section{Rotation and the Inertia Groupoid}

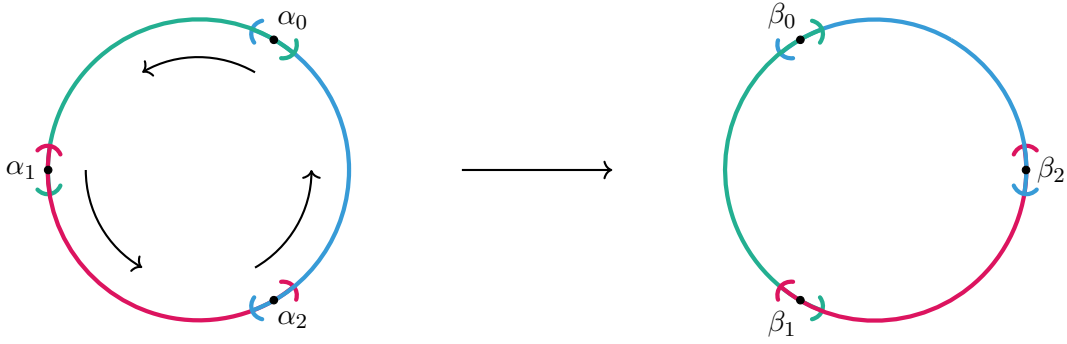
\begin{figure}
    \centering
    \begin{tikzpicture}
    
    \def \pt {0.05};

    \begin{scope}[shift={(4,2)}]
    
        \draw (0,0) circle [radius=2];
    
        \draw [ultra thick, myblue, -)] (0:2) arc (0:70:2);
        \draw [ultra thick, mygreen, (-)] (50:2) arc (50:190:2);
        \draw [ultra thick, mypink, (-)] (170:2) arc (170:310:2);
        \draw [ultra thick, myblue, (-] (290:2) arc (290:360:2);
        
        \draw [fill] (60:2) circle [radius=\pt] node [above right, xshift=-1mm] {\small $\alpha_0$};
        \draw [fill] (180:2) circle [radius=\pt] node [left] {\small $\alpha_1$};
        \draw [fill] (300:2) circle [radius=\pt] node [below right, xshift=-1mm] {\small $\alpha_2$}; 
    
        \draw [thick, ->] (60:1.5) arc (60:120:1.5); 
        \draw [thick, ->] (180:1.5) arc (180:240:1.5);
        \draw [thick, ->] (300:1.5) arc (300:360:1.5);
    
    \end{scope}
    
    \begin{scope}[shift={(13,2)}]
    
        \draw (0,0) circle [radius=2];
    
        \draw [ultra thick, myblue, -)] (60:2) arc (60:130:2);
        \draw [ultra thick, mygreen, (-)] (110:2) arc (110:250:2);
        \draw [ultra thick, mypink, (-)] (230:2) arc (230:370:2);
        \draw [ultra thick, myblue, (-] (-10:2) arc (-10:60:2);
    
        \draw [fill] (240:2) circle [radius=\pt] node [below left, xshift=1mm] {\small $\beta_1$};
        \draw [fill] (0:2) circle [radius=\pt] node [right] {\small $\beta_2$};
        \draw [fill] (120:2) circle [radius=\pt] node [above left, xshift=1mm] {\small $\beta_0$}; 

    \end{scope}
    
    \draw [thick, ->] (7.5,2) -- (9.5,2);

\end{tikzpicture}
    \caption{Rotation of an admissible cover by $z = 1/6$.}
    \label{fig:rotation}
\end{figure}

Any loop construction should have a natural $\R$-action which comes from rotating the loops. Let $\U$ be an admissible cover of $S^1$ given by $0 < \alpha_0 < ... < \alpha_n \leq 1$ and let $z \in \R$. Define 
\[
\S^1_{\U^z} = 
\left(
\begin{array}{c}
    \bigsqcup_{i,j} (W_i^z \cap W_j^z)  \times \Z \\
     \downdownarrows \\
     \bigsqcup_i W_i^z
\end{array}
\right)
\quad \text{where} \quad
W_i^z = W_i + z. 
\]
The notation is explained as follows: for each $i$, let $\beta_i \in (0,1]$ be a representative of $\alpha_i + z$ modulo $\Z$. This gives a new set of vertices $\{\beta_0, \dotsc, \beta_n\} \subset (0,1]$ and, using the same $\epsilon$, we have a new admissible cover $\U^z$ with $\S^1_{\U^z}$ as above. Formally, this is just a rotation of the circle, see Figure \ref{fig:rotation}.

There is a functor $\iota \colon \S^1_{\U^z} \to \S^1_{\U}$ defined by undoing the translation by $z$,
\[
\iota_0(x,i) = (x-z, i), \quad \iota_1(x,n,i,j) = (x-z, n, i, j).
\]
Given a loop $\Phi$ associated to $\U$ we define the rotated loop $\Phi^z$ to be  $\Phi^z  := \Phi \circ \iota$. Explicitly,
\[
\Phi^z_0(x,i) = \Phi_0(x-z, i)
\quad \text{and} \quad
\Phi^z_1(x,n,i,j) = \Phi_1(x-z,n,i,j).
\]
Therefore, the rotation action on $\L\G$ is the function
\[
(\L\G)_0 \times \R \to (\L\G)_0, \quad (\Phi, z) \mapsto \Phi^z.
\]
We now ask which loops are rotation invariant. Such loops form a sub-groupoid of $\L\G$ that is denoted $\L\G^{\R}$. Consider a simple example.

\begin{example}\label{ex:rotationinvarianttrivialloop}
Consider a loop $\Phi \colon \R \mmod \Z \to \G$ that is invariant under rotation. This means that if $z \in \R$, then $\Phi = \Phi^z$ on a common refinement of $\U$ and $\U^z$. Unpacking this, we can conclude that for all $z \in \R$,
\[
\Phi_0(x) = \Phi_0(x-z)
\quad \text{and} \quad
\Phi_1(x,n) = \Phi_1(x-z,n).
\]
Therefore, $\Phi_0$ and $\Phi_1(-,n)$ are constant and $\Phi_1(x,n) = g^n$ where $g :=\Phi_1(x,1)$. In other words, we have an automorphism of a fixed object $\Phi_0(x) \in \G_0$. This is an object in the inertia groupoid $\Lambda G$.
\end{example}

This example, together with the following result, leads to a nice description of $\L\G^{\R}$.

\begin{lemma}\label{lem:rotloopstriv}
Every object in $\L\G^{\R}$ is isomorphic to another object in $\L\G^{\R}$ that is associated to the trivial cover of $S^1$.
\end{lemma}

The idea behind the proof is that rotation invariance of $\Phi \in \L\G^{\R}_0$ implies that $\Phi_1$ and $\Phi_0$ are locally constant. The loop is then a collection of constant paths all of which are connected by morphisms in the image of $\Phi_1$, allowing the construction of an isomorphic loop associated to the trivial cover. For a complete proof, we refer the reader to \cite[Lemma 3.6.3]{LU02} and mention that this proof uses similar technique to Proposition \ref{prop:trivialcoverXmmodG}, for which a detailed proof sketch is provided.

Example \ref{ex:rotationinvarianttrivialloop} and Lemma \ref{lem:rotloopstriv} lead to the following \cite[Theorem 3.6.6]{LU02}:

\begin{proposition}
 The groupoid $\L\G^{\R}$ is Morita equivalent to the inertia groupoid $\Lambda \G$.
\end{proposition}

\noindent Of particular note is that there is an inclusion functor $\iota \colon \Lambda \G \hookrightarrow \L\G $ defined as follows:
\begin{itemize}
    \item Given an object $g$ in $\Lambda \G$, the loop $\Phi = \iota_0(g) \colon \R \mmod \Z \to \G$ is defined by
    \[
    \Phi_0(t) := s(g) = t(g)
    \quad \text{and} \quad
    \Phi_1(t,1) := g.
    \]
    This can be thought of as a constant loop at $s(g) = t(g)$.
    \item Given a morphism $(g,h)$ in $\Lambda \G$, the arrow $\Lambda = \iota_1(g, h)$ is given by 
    \[
    \Lambda \colon \R \to \G_1, \quad \Lambda(t) = h.
    \]
    This defines a natural transformation from $\Phi=\iota_0(g)$ to $\Psi = \iota_0(h^{-1}gh)$:
    \begin{gather*}
        s(\Lambda(t)) = s(h) =  \Phi_0(t), \\
        t(\Lambda(t,n)) = t(h) = \Psi_0(t), \text{ and} \\
        \Lambda(t) \cdot \Psi(t,n) = h \cdot (h^{-1} g^n h) = g^n h = \Phi(t,n) \cdot \Lambda(t),
    \end{gather*}
    The third calculation is the naturality condition.
\end{itemize}


\section{Loops of \texorpdfstring{$X \mmod G$}{X/G}}    

 In this section, we will consider the loop groupoid of $X \mmod G$ and find a  description of $\L(X \mmod G)$ in terms of certain paths in $X$. The first result which simplifies things is the following \cite[Lemma 4.1.1]{LU02}:

\begin{proposition}\label{prop:trivialcoverXmmodG}
Let $X$ be a connected space acted on by a finite group $G$. Every object in the loop groupoid $\L(X \mmod G)$ is isomorphic to a loop associated to the trivial cover. 
\end{proposition}

 A sketch of the proof, aimed to complement that of \cite{LU02}, is provided towards the end of the chapter, starting on Page \pageref{proofsketch}. Assuming the proposition, every object in $\L(X \mmod \G)$ is isomorphic to a smooth functor $\Phi \colon \R \mmod \Z \to X \mmod G$. As discussed in Example \ref{ex:trivialloop}, such a functor is fully determined by the value of $\Phi_1$ on $\R \times \{1\}$. Moreover,
\[
\Phi_1( \R \times \{1\}) \subseteq X \times \{g\}
\]
for some $g \in G$. This is because $\Phi_1(\R \times \{1\})$ lies in a connected component of $X \times G$ and $G$, being a finite group, has the discrete topology. In fact,
\[
\Phi_1(x,1) = (\Phi_0(x), g)
\stext{and}
\Phi_0(x+1) = t(\Phi_1(x,g)) = \Phi_0(x)\cdot g.
\]
Conversely, if we start with a function $\varphi \colon \R \to X$ satisfying $\varphi(x+1) = \varphi(x) \cdot g$, then 
\[
\Phi_1(x,1) := (\varphi(x), g)
\]
defines a loop in $X \mmod G$. We conclude that loops in $X \mmod G$ can be identified with pairs $(\varphi, g)$ such that $\varphi(x+1) = \varphi(x) \cdot g$. Establish the set
\[\label{Pg}
\P_g = \{ \varphi \colon \R \to X \mid \varphi(x+1) = \varphi(x) \cdot g \}.
\]
A morphism $\Lambda \colon \Phi \to \Psi$ in $\L(X \mmod G)$, where $\Phi, \Psi \colon \R \mmod \Z \to X \mmod G$, is a function 
\[
\Lambda \colon \R \to X \times G,
\stext{where}
\Lambda(\R) \subset X \times \{k\}
\]
for some $k \in G$. In fact, for $\Lambda(x)$ to be a morphism in $X \mmod G$ from $\Phi_0(x)$ to $\Psi_0(x)$, we must have $\Lambda(x) = (\Phi_0(x), k)$. This means that $\Psi_0(x) = \Phi_0(x) \cdot k$ for all $x \in \R$. If $\Phi_0 \in \P_g$ and $\Psi_0 \in \P_h$, then the naturality of $\Lambda$ implies that
\[
g \cdot k = k \cdot h 
\quad \implies \quad 
h = k^{-1}gh 
\quad \implies \quad 
\Psi_1(x,1) = (\Phi_0(x) \cdot k, k^{-1}gk).
\]
This leads to the following:

\begin{proposition}
The loop groupoid $\L(X \mmod G)$ is Morita equivalent to the groupoid
\[
\left( \bigsqcup\nolimits_g \P_g \right) \mmod G
\]
where the $G$-action is given by
\[
\left( \bigsqcup\nolimits_g \P_g \right) \times G \to \left( \bigsqcup\nolimits_g \P_g \right),
\quad \bigl( (\varphi, g), h \big) \mapsto (\varphi \cdot h, h^{-1} g h),
\]
so that if $\varphi \in \P_g$, then $\varphi \cdot h \in \P_{h^{-1}gh}$.
\end{proposition}

\noindent We can go one step further. If $g$ and $h$ are not in the same conjugacy class, then there are no morphisms between elements of $\P_g$ and elements of $\P_h$. We can therefore write
\[
\left( \bigsqcup\nolimits_g \P_g \right) \mmod G = \bigsqcup_{[g]}\left[ \left( \bigsqcup\nolimits_{h \in [g]} \P_h \right) \mmod G \right],
\]
where the first union on the right-hand side runs over the conjugacy classes $[g]$ of $G$. Acting by $k \in G$ identifies the objects in $\P_g$ with the objects in $\P_{kgk^{-1}}$. If $g = kgk^{-1}$, that is, $k$ is an element of the centraliser $C_g$ of $g$, then $k$ acts as an automorphism of $\P_g$. We hope that for each $h \in [g]$ we can identify $\P_h$ with $\P_g$. Consider the functor
\[
F: \P_g \mmod C_g \to \left( \bigsqcup\nolimits_{h \in [g]} \P_h \right) \mmod G
\]
defined by $F_0(\varphi) = (\varphi, g)$ and $F_1(\varphi, h) = (\varphi, g, h)$. The $C_g$-action on the left-hand side is $(\varphi \cdot h)(x) = \varphi(x) \cdot h$. This functor is an equivalence of groupoids, leading to the following result \cite[Corollary 6.1.2]{LU02}.

\begin{proposition}
The loop groupoid $\L(X \mmod G)$ is Morita equivalent to the groupoid
\[
\bigsqcup _{[g]}  \P_g \mmod C_g 
\]
where the disjoint union runs over the conjugacy classes $[g]$ of $G$ and the $C_g$-action on $\P_g$ is given by
\[
\P_g \times C_g \to \P_g, \quad (\varphi, h) \mapsto \varphi \cdot h.
\]
\end{proposition}

There is a similarly simple formulation of the loops of $X \mmod G$ that are invariant under translation. The action on $\L(X \mmod G)$ by the real numbers corresponds to the following action on $\P_g$:
\begin{align*}
\P_g \times \R \,\,&\to \,\,  \P_g \\
(\varphi, z) \,\,&\mapsto \,\, \varphi^z
\end{align*}
where $\varphi^z(t) = \varphi(t-z)$. The functions invariant under this action satisfy $\varphi = \varphi^z$ for all $z \in \R$ and are hence constant functions. Then, if $x = \varphi(t)$ for all $t$, we have that
\[
x = \varphi(t+1) = \varphi(t)\cdot g = x \cdot g.
\]
Therefore, $(\P_g)^{\R} \cong X^g$. This leads to \cite[\textsection 6.2.1]{LU02},
\[
\L(X \mmod \G)^{\R} \cong \bigsqcup_{[g]} X^g \mmod C_g.
\]
This matches the description of the inertia groupoid of $X \mmod G$, after we choose representatives in each conjugacy class.

\subsubsection{Sketch of the proof of Proposition \ref{prop:trivialcoverXmmodG}} \label{proofsketch}

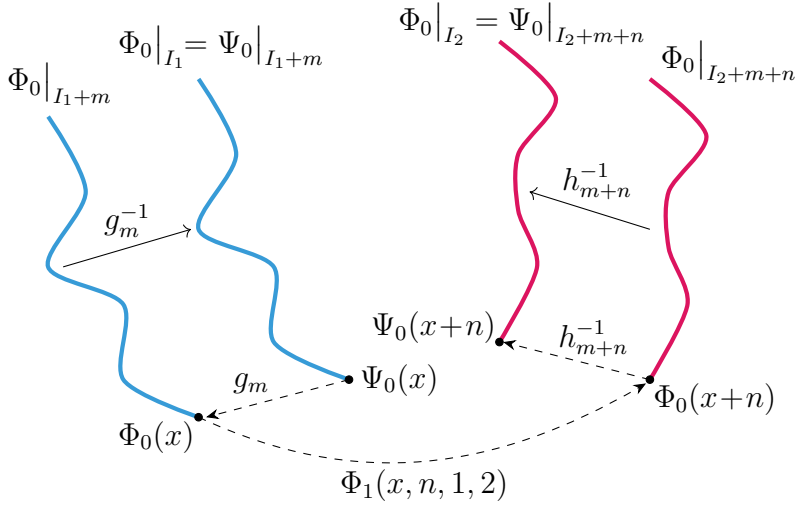
\begin{figure}
    \centering
    \begin{tikzpicture}

    \def \pt {0.05};
    \def \xs {2}
    \def \ys {-0.5}
    
    \def \myarrow {-{Stealth[scale=1,angle'=45]}};
    
    \draw [ultra thick, myblue] plot [smooth] coordinates {(0, -2) (-1,-1.5) (-1,-0.5) (-2,0) (-1.5,1) (-2,2)};
    \draw [ultra thick, mypink] plot [smooth] coordinates {(2, -1.5) (2.5,-0.5) (2.25,0.25) (2.25,1) (2.75,1.75) (2,2.5)};
    
    \begin{scope}[shift={((-2,-0.5)}]
        \draw [ultra thick, myblue] plot [smooth] coordinates {(0, -2) (-1,-1.5) (-1,-0.5) (-2,0) (-1.5,1) (-2,2)};
    \end{scope}
    
    \begin{scope}[shift={((2,-0.5)}]
        \draw [ultra thick, mypink] plot [smooth] coordinates {(2, -1.5) (2.5,-0.5) (2.25,0.25) (2.25,1) (2.75,1.75) (2,2.5)};
    \end{scope}
        
    \node [above, xshift=3mm] at (-2,2) {$\Phi_0\big|_{I_1} \!\! = \Psi_0\big|_{I_1+m}$};
    \node [above left] at (-3,1.5) {$\Phi_0\big|_{I_1+m}$};
    
    \node [above, xshift=3mm] at (2,2.3) {$\Phi_0\big|_{I_2} = \Psi_0\big|_{I_2+m+n}$};
    \node [above right] at (4,1.8) {$\Phi_0\big|_{I_2+m+n}$};
    
    \draw [<-] (-2.1,0) -- (-3.8,-0.5) node [midway, above] {$g_m^{-1}$};
    \draw [<-] (2.4,0.5) -- (4,0) node [midway, above, xshift=1mm] {$h_{m+n}^{-1}$};
    
    \draw [dashed, \myarrow] (0,-2) -- (-1.9, -2.5);
    \draw [dashed, \myarrow] (-2,-2.5) to[out=-25, in=210] (3.95,-2.05);
    \draw [dashed, \myarrow] (4,-2) -- (2.05,-1.5);
    
    \draw [fill] (-2,-2.5) circle [radius=\pt] node [below left, shift={(0.1,0.1)}] {$\Phi_0(x)$};
    \draw [fill] (0,-2) circle [radius=\pt] node [right] {$\Psi_0(x)$};
    \draw [fill] (4,-2) circle [radius=\pt] node [below right, shift={(-0.1,0.1)}] {$\Phi_0(x\!+\!n)$};
    \draw [fill] (2,-1.5) circle [radius=\pt] node [left, yshift=2mm, xshift=1mm] {$\Psi_0(x\!+\!n)$};
    
    \node at (-1.3, -2.05) {$g_m$};
    \node at (1,-3.4) {$\Phi_1(x, n, 1,2)$};
    \node at (3.25, -1.45) {$h^{-1}_{m+n}$};
    
\end{tikzpicture}
    \caption{The loops $\Phi$ and $\Psi$.}
    \label{fig:phiandpsi}
\end{figure}

\begin{figure}
    \centering
    \begin{tikzpicture}

    \def \pt {0.05};

    \begin{scope}[shift = {(4,3)}]
    
        \draw [ultra thick, mypink, (-)] (-0.5,0) -- (1.8,0) node [black, midway, above] {$I_2$};
        \draw [ultra thick, myblue, (-)] (-2,0) -- (0,0) node [black, midway, above] {$I_1$};
    
        \draw [fill] (-1.75,0) circle [radius=\pt] node [below, yshift=-1mm] {$\alpha_0$};
        \draw [fill] (-0.25,0) circle [radius=\pt] node [below, yshift=-1mm] {$\alpha_1$};
        \draw [fill] (1.55,0) circle [radius=\pt] node [below] {$\alpha_0+1$};
    
    \end{scope}
    
    \begin{scope}[shift = {(12,3)}]
    
        \draw [ultra thick, myblue] plot [smooth] coordinates {(0, -2) (-1,-1.5) (-1,-0.5) (-2,0) (-1.5,1) (-2,2)};
        \draw [ultra thick, mypink] plot [smooth] coordinates {(2, -1.5) (2.5,-0.5) (2.25,0.25) (2.25,1) (2.75,1.75) (2,2.5)};
        
        \draw [dashed, -{Stealth[scale=1,angle'=45]}] (0,-2) to[out=-30, in=-120] node [midway, below] {$k_1$} (2,-1.5);
        \draw [dashed, -{Stealth[scale=1,angle'=45]}] (2, 2.5) to[out=150, in=60] node [midway, above] {$k_2$} (-2,2);
        
        \node at (-1.5,-0.7) {$\Psi_0\big|_{I_1}$};
        \node at (2.8,0.5) {$\Psi_0\big|_{I_2}$};
    
    \end{scope}

    \draw [thick, ->] (7,3) -- (9,3) node [midway, above] {$\Psi$};

\end{tikzpicture}
    \caption{A loop in $X \mmod G$ associated to an admissible cover with two vertices.}
    \label{fig:loopinXmmodG}
\end{figure}
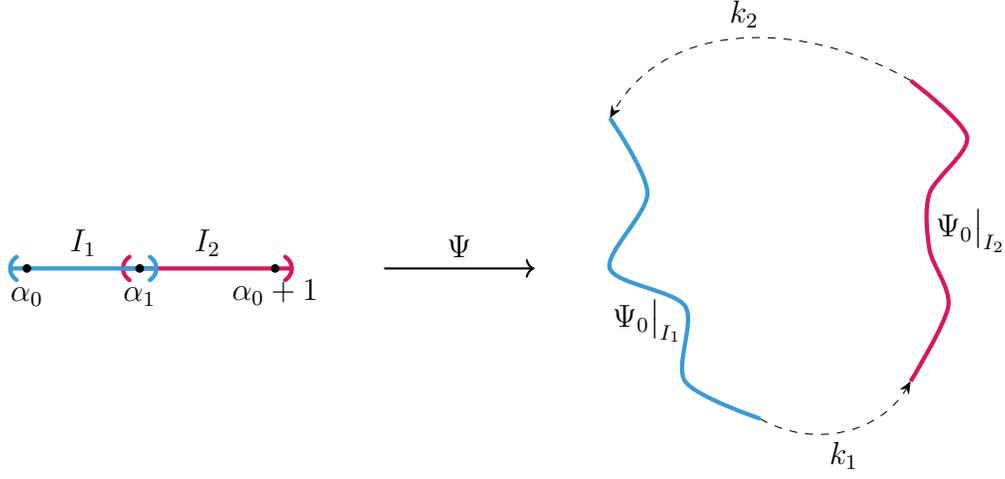

For the remainder of the chapter, we detail a sketch that is intended to complement Lupercio and Uribe's complete proof in \cite{LU02}. Consider a loop $\Phi$ in $X \mmod G$ associated to an admissible cover $\U$ with two vertices
\[
0 < \alpha_0 < \alpha_1 \leq 1.
\]
On objects, $\Phi$ is a collection of paths in $X$, namely for each $m \in \Z$ we have paths 
\[
\Phi_0 |_{I_1 + m} \colon I_1 + m \to X
\quad \text{and} \quad
\Phi_0 |_{I_2 + m} \colon I_2 + m \to X.
\]
Consider the former maps first. For each $m \in \Z$ the paths $\Phi_0|_{I_1}$ and $\Phi_0|_{I_1 + m}$ are connected by morphisms in $\Phi_1(I_1 \times \{m\})$. As $I_1 \times \{m\}$ is connected, $\Phi_1$ is continuous, and $G$ has the discrete topology, we know that
\[
\Phi_1(I_1 \times \{m\}) \subseteq X \times \{g_m\}
\]
for some $g_m \in G$. This means that
\[
\Phi_0(x) \cdot g_m = \Phi_0(x+m)
\quad \text{for all} \quad
x \in I_1.
\]
Similarly, there exists $h_m \in G$ such that $\Phi_0(x) \cdot h_m = \Phi_0(x+m)$ for all $x \in I_2$. We will apply $g_m^{-1}$ and $h_m^{-1}$ appropriately to define a loop $\Psi$ which satisfies
\begin{equation}\label{eq:whatwewant}
\Psi_0 |_{I_i + m} = \Psi_0 |_{I_i} = \Phi_0|_{I_0}.
\end{equation}
In this case, $\Psi_0$ is just a finite collection of paths in $X$ and $\Psi_1$ connects the ends of these paths together. Refer to Figure \ref{fig:phiandpsi} for a schematic description of how $\Psi$ is defined. 
On objects $\Psi$ is defined as follows:
\[
\Psi_0(x) = 
\begin{cases}
\Phi_0(x) \cdot g_m^{-1},    & x \in I_1 + m \\
\Phi_0(x) \cdot h_m^{-1},    & x \in I_2 + m.
\end{cases}
\]
On morphisms we have four cases:
\begin{align*}
    \Psi_1(x, n, 1,1) &:= g_m \cdot \Phi_1(x,n,1,1) \cdot g_{m+n}^{-1}, \quad  x \in I_1 + m, \\
    \Psi_1(x,n,1,2) &:= g_m \cdot \Phi_1(x,n,1,2) \cdot h_{m+n}^{-1}, \quad  x \in I_1 + m,\\
    \Psi_1(x,n,2,1) &:= h_m \cdot \Phi_1(x,n,2,1) \cdot g_{m+n}^{-1}, \quad x \in I_2 + m, \\
    \Psi_1(x,n,2,2) &:= h_m \cdot \Phi_1(x,n,2,2) \cdot h_{m+n}^{-1}, \quad x \in I_2 + m.
\end{align*}
There is an abuse of notation here: $\Phi_1(x,n,i,j) \in X \times G$ but $g_m,h_m \in G$. The first line, for instance, would be correctly written as
\[
\Psi_1(x, n, 1,1) := (\Psi_0(x), g_m) \cdot \Phi_1(x,n,1,1) \cdot (\Phi_0(x+n), g_{m+n}^{-1}).
\]
Now, if $x+m \in I_1 + m$, then
\[
\Psi_0 (x+m) = \Phi_0(x+m) \cdot g_m^{-1} = \Phi_0(x) \cdot g_m \cdot g_m^{-1} = \Phi_0(x) = \Psi_0(x). 
\]
Therefore, condition \eqref{eq:whatwewant} is satisfied. $\Psi$ is drawn again in Figure \ref{fig:loopinXmmodG}. Note that this is the same way we've previously draw loops, because we omit the ``copies''. Now there is nothing to omit because all of the copies have been identified.

There is a natural transformation $\eta \colon \Phi \to \Psi$ defined by
\[
\eta \colon (\S^1_{\U})_0 \to X \times G, \quad \eta(x) = 
\begin{cases}
    (\Phi_0(x), g_m^{-1}),  & x \in I_1 + m, \\
    (\Phi_0(x), h_m^{-1}),  & x \in I_2 + m.
\end{cases}
\]
Indeed, if $x \in I_1 + m$, then 
\[
\Phi_0(x) \cdot g_m^{-1} = \Psi_0(x),
\]
by the definition of $\Psi$. The same holds for $x \in I_2 + m$. If we consider $(x,n,0,1) \in (\S^1_{\U})_1$, with $x \in I_1 + m$, then the following diagram commutes:
\[
\begin{tikzcd}[column sep=1.7cm]
 \Phi_0(x) \arrow[r, "{\Phi_1(x,n,0,1)}"] \arrow[d, "\eta(x)", swap] &  \Phi_{0}(x+n) \arrow[d, "\eta(x+n)"] \\
\Psi_0(x) \arrow[r, "{\Psi_1(x,n,0,1)}"]           & \Psi_0(x+n)         
\end{tikzcd}
\]
This is because
\[
\eta(x) \cdot \Psi_1(x,n,0,1) = g_m^{-1} \cdot g_m \cdot \Phi_1(x,n,0,1) \cdot h_{m+n}^{-1} = \Phi_1(x,n,0,1) \cdot \eta(x+n).
\]
Naturality for the remaining morphisms can be checked by the reader. The point is that $\Phi$ is isomorphic in $\L\G$ to a loop which doesn't contain infinitely many paths - the only ones we need to consider are the restrictions to $I_1$ and $I_2$. This is the first step in producing a loop associated to the trivial cover that is isomorphic to $\Phi$.

The next step is to perform a similar operation to $\Psi$ to obtain a loop associated to the trivial cover. This amounts to joining the two disjoint paths at $\Psi_1(I_1 \cap I_2)$ so that they form a single path. First note that 
\[
\Psi_1(I_1 \cap I_2 \times \{0\}) \subseteq X \times \{k_1\} \stext{and} \Psi_1(I_2 \cap (I_1 + 1) \times \{-1\}) \subseteq X \times \{k_2\}
\]
for some $k_0, k_1 \in G$. The $-1$ is present because $\Psi_1(x, -1, 2, 1)$ is a morphism from $\Psi_0(x,2)$ to $\Psi_0(x-1,1)$ where $x \in I_2 \cap (I_1+1)$ and hence $x-1 \in I_1$. 

\begin{figure}
    \centering
    \begin{tikzpicture}

    \def \pt {0.05};

    \draw [ultra thick, myblue] plot [smooth] coordinates {(0, -2) (-1,-1.5) (-1,-0.5) (-2,0) (-1.5,1) (-2,2)};
    \draw [ultra thick, mypink] plot [smooth] coordinates {(2, -1.5) (2.5,-0.5) (2.25,0.25) (2.25,1) (2.75,1.75) (2,2.5)};
    \draw [ultra thick, mypink] plot [smooth] coordinates {(0, -2) (0.5,-1) (0.25,-0.25) (0.25,0.5) (0.75,1.25) (0,2)};

    \draw [dashed, -{Stealth[scale=1,angle'=45]}] (2, 2.5) to[out=150, in=60] node [midway, above] {$k_2$} (-1.9,2.1);
    \draw [dashed, -{Stealth[scale=1,angle'=45]}] (0, 2) -- (1.9,2.4) node [midway, above] {$k_1$};
    \draw [dashed, -{Stealth[scale=1,angle'=45]}] (0, 2) -- (-1.8,2) node [midway, above] {$k_1k_2$};
    \draw [dashed, -{Stealth[scale=1,angle'=45]}] (0,-2) to[out=-30, in=-120] node [midway, below] {$k_1$} (2,-1.5);
    
    \draw [->] (2,0.5) -- (0.5,0) node [midway, above] {$k_1^{-1}$};
        
    \node at (-1.5,-0.7) {$\Psi_0\big|_{I_1}$};
    \node at (2.8,0.5) {$\Psi_0\big|_{I_2}$};
    \node at (-0.3,0.5) {$\Theta_0\big|_{I_2}$};

\end{tikzpicture}
    \caption{Modify $\Psi$ so that it is associated to a trivial cover.}
    \label{fig:newloopinXmmodG}
\end{figure}
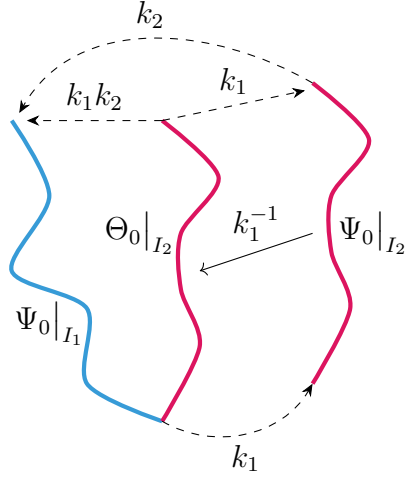

Define a new loop $\Theta \colon \S^1_{\U} \to X \mmod G$ as follows:
\[
\Theta_0(x) = 
\begin{cases}
\Psi_0(x), & x \in I_1, \\
\Psi_0(x) \cdot k_1^{-1}, & x \in I_2, \\
\end{cases} \quad \text{and},
\]\[
\Theta_1(x) = 
\begin{cases}
(\Psi_0(x), 1),           & x \in I_1 \cap I_2 \text{  or  } x \in (I_1+1) \cap I_2, \\
(\Psi_0(x) \cdot k_1^{-1}, k_1k_2),   & x \in I_2 \cap I_1 \text{  or  } x \in I_2 \cap I_1.
\end{cases}
\]
Here we are applying $k_1^{-1}$ to the path $\Psi_0|_{I_2}$ so that it connects to $\Psi_0|_{I_1}$ at $\alpha_1$ as in Figure \ref{fig:newloopinXmmodG}. The morphisms are chosen to be compatible with this change.

Now $\Theta$ is a loop that is associated to the trivial cover: if $x \in I_1 \cap I_2$, then $\Theta_1(x)$ is the identity at $\Theta_0(x)$. We can therefore write $\Theta_0$ as a map
\[
\Theta_0 \colon \R \to X \stext{such that} \Theta_0(x+1) = \Theta_0(x)\cdot (k_1k_2).
\]
We recognise this as a loop in $X \mmod G$ associated to the trivial cover, namely $\Theta_0 \in \P_{k_1k_2}$.

There is a natural transformation $\xi \colon \Psi \to \Theta$ given by
\[
\xi \colon (\S^1_{\U})_0 \to X \times G,  
\quad
\xi(x) =
\begin{cases}
    (\Phi_0(x), 1), & x \in I_1, \\
    (\Phi_0(x), k_1^{-1}), & x \in I_2.
\end{cases}
\]
Therefore, $\Phi$, $\Psi$ and $\Theta$ are isomorphic in $\L(X\mmod G)$. In particular $\Phi$ is isomorphic to a loop associated to the trivial cover.

To summarise, we started with a loop $\Psi$ and first obtained an isomorphic loop $\Psi$ such that $\Psi|_{I_i} \!= \Psi|_{I_i + n}$ for each $i$. In the next step, we obtained a loop $\Theta$ such that the disjoint paths were joined together. In this example, $\Phi$ was associated to a very simple admissible cover - only two intervals. Obtaining $\Psi$ in the general case is a matter of applying the same technique to a greater number of $I_i$. Then, $\Theta$ can be obtained from $\Psi$ be repeatedly joining each disjoint path to the next; in our case this only had to be done once.

\chapter{Structures on Groupoids}\label{StructuresOnGroupoids}


\section{Vector Bundles}

There are two equivalent definitions of a vector bundle on a groupoid and each extends to a definition of vector bundles on orbifolds. The first definition is seen for instance in \cite[\textsection 5.4]{Morbifolds}. Unless otherwise stated we assume that all groupoids are topological groupoids.

\begin{definition}
A \emph{vector bundle} on a groupoid $\G$ is a (right) $\G$-space $E$ such that the anchor map $\pi \colon E \to \G_0$ is a vector bundle and the $\G$-action is fibre-wise linear. Explicitly, this consists of the following data:
\begin{enumerate}[label=(\roman*)]
    \item A vector bundle $\pi \colon E \to \G_0$.
    \item A right $\G$-action,
    \[
    \mu \colon E \mathbin{_\pi \times_t} \G_1 \to E, \quad (e,g) \mapsto e\cdot g,
    \]
    which satisfies, for $z \xrightarrow{g} y \xrightarrow{h} x$ in $\G$,
    \begin{itemize}
        \item $\pi(e \cdot g) = s(g)$,
        \item $e \cdot 1_x = e$,  
        \item $(e \cdot h) \cdot g = e \cdot (gh)$, and
        \item $g \colon E_{t(g)} \to E_{s(g)}$ is a linear isomorphism for each $g \in \G_1$.
    \end{itemize}
    
\end{enumerate}
\end{definition}

\noindent An alternative definition follows a general framework for defining structures over groupoids.

\begin{definition}\label{def:vectorbundledef2}
A \emph{vector bundle} on a groupoid $\G$ consists of the following data:
\begin{enumerate}[label=(\roman*)]
    \item A vector bundle $\pi \colon E \to \G_0$.
    \item An isomorphism $\mu \colon t^* E \to s^*E$ of vector bundles over $\G_1$ which satisfies
    \begin{itemize}
        \item $u^* \mu = \id_{E}$ where $u \colon \G_0 \to \G_1$ is the unit map of the groupoid.
        \item $d_2^* \mu \circ d_0^*\mu = m^* \mu$ where $d_0, m, d_2 \colon \G_2 \to \G_1$ are the maps given by 
        \[d_0(g,h) = h, \quad  m(g,h) = gh,  \quad\text{and} \quad  d_2(g,h) = g.\]
        This is called the associativity condition.
    \end{itemize}
\end{enumerate}
\end{definition}

\begin{remark}
The second condition on the isomorphism $\mu$ is encapsulated by the following commutative diagram:

\[
\begin{tikzcd}
    & m^*t^*E \arrow[r, "m^*\mu"] & m^*s^*E \arrow[rd, equals] &   \\
d_0^*t^*E \arrow[ru, equals] \arrow[rd, "d_0^*\mu"] &      & & d_2^*s^*E \\
    & d_0^*s^*E \arrow[r, equals] & d_2^*t^*E \arrow[ru, "d_2^*\mu"] &  
\end{tikzcd}
\]
To be even more explicit, since $\mu$ is a vector bundle isomorphism we can write
\[
\mu(e, g) = (\mu_1(e,g), g)
\]
for some function $\mu_1 \colon t^*E \to E$. Then the conditions on $\mu$ become
\[
\mu_1(e, 1_{\pi(e)}) = e, \quad \text{and} \quad \mu_1(\mu_1(e,h), g) = \mu(e, gh).
\]
This indicates the relationship between the two definitions; since $t^*E = E\mathbin{_\pi \times_t} \G_1$ we see that $\mu_1$ defines an action on $E$ as per the first definition. In the other direction, if we have a $\G$-action on $E$, then we can define
\[
\mu \colon t^*E \to s^*E, \quad \mu(e, g) = (e\cdot g, g).
\]
The reader can check that this shows that the two definitions are equivalent. We will shift between the two definitions depending on which is convenient at the time.
\end{remark}

\begin{definition}
Let $E$ and $F$ be vector bundles on a groupoid $\G$. An \emph{isomorphism} from $E$ to $F$ is a vector bundle isomorphism $f \colon E \to F$ that is compatible with the $\G$-action i.e. the following diagram commutes:
\[
\begin{tikzcd}
t^*E \arrow[r] \arrow[d, "t^*f", swap] & s^*E \arrow[d, "s^*f"] \\
t^*F \arrow[r]           & s^*F          
\end{tikzcd}
\]
Equivalently, $f(e \cdot g) = f(e) \cdot g$.
\end{definition} 

\begin{definition}
Let $\phi \colon \G \to \H$ be a homomorphism of Lie groupoids and let $E$ be a vector bundle on $\H$. The \emph{pullback bundle} is defined as $\phi^*E = E \mathbin{_\pi \times_{\phi_0}} \G_0$ with the action 
\[
(e, x) \cdot g = (e \cdot \phi_1(g), y), \quad g\colon y \to x.
\]

\end{definition}

Let $\vect{\G}$ denote the category of vector bundles on $\G$.

\begin{proposition}
If $\G$ and $\H$ are Morita equivalent, then there is an equivalence of categories between $\vect{\G}$ and $\vect{\H}$.
\end{proposition}
\begin{proof}
This is best proved in the context of bibundles. The key idea is that a vector bundle on $\G$ is equivalent to a $\GL_n(\C)$-bundle on $\G$, which in turn is equivalent to a bibundle from $\G$ to $\mathbb{B}\!\GL_n(\C)$. If $\G$ and $\H$ are Morita equivalent, then they are isomorphic in the 2-category of Lie groupoids with bibundles. This implies an equivalence of categories between the Hom-categories $\Hom(\G, \mathbb{B}\!\GL_n(\C))$ and $\Hom(\H, \mathbb{B}\!\GL_n(\C))$, which are equivalent to $\vect \G$ and $\vect \H$ respectively. A discussion is found in \cite[\textsection 2.6]{amenta}.
\end{proof}

The notion of $\GL_n(\C)$-valued 1-cocycles are important when discussing vector bundles; every vector bundle is classified up to isomorphism by its associated 1-cocycle. Before further investigating vector bundles on groupoids, we'd like an appropriate notion of cocycles on groupoids. We begin by framing the definition of a cocycle in terms of groupoids. Let $M$ be a space with an open cover $\U = \{U_i\}$ and let $\co = (\co_{ij})$ be a 1-cocycle associated to $\U$. Recall the open cover groupoid, 
\[
M[\U] = \biggl( \bigsqcup\nolimits_{i,j} U_i \cap U_j \rightrightarrows \bigsqcup\nolimits_{i} U_i \biggr).
\]
The 1-cocycle $\co$ can be reinterpreted as the map 
\[
\co \colon M[\U]_1 \to \GL_n(\C), \quad \co(x,i,j) = \co_{ij}(x).
\]
The cocycle condition is a condition on triple intersections; if $x \in U_i \cap U_j \cap U_k$, then $\co_{ij}(x) \co_{jk}(x) = \co_{ik}(x)$. This can be expressed as a condition on
\[
M[\U]_2 = \bigsqcup\nolimits_{i,j,k} U_i \cap U_j \cap U_k.
\]
Precisely, this condition is that if $(x,i,j,k) \in M[\U]_2$, then
\[
\co(x,i,j) \cdot \co(x,j,k) = \co(x,i,k).
\]
In other words, 
\begin{equation}\label{eq:1cocycle}
(\co \circ d_0) \cdot (\co \circ d_2) = (\co \circ d_1).
\end{equation}
To summarise, a $\GL_n(\C)$ valued 1-cocycle $(\co_{ij})$ associated to the open cover $\U$ is equivalent to a map $\co \colon M[\U]_1 \to \GL_n(\C)$ satisfying \eqref{eq:1cocycle}. This formalism can be extended to general groupoids:

\begin{definition}
A function $\co \colon \G_p \to \GL_n(\C)$ is a $p$-cocycle if $\co$ satisfies the following cocycle condition over $\G_{p+1}$:
\[
\prod_{i=0}^p (\co \circ d_i)^{(-1)^i} = 1.
\]
\end{definition}
For example, a 2-cocycle would satisfy
\[
(\co \circ d_0) \cdot (\co \circ d_1)^{-1} \cdot (\co \circ d_2) \cdot (\co \circ d_3)^{-1} = 1.
\]
This definition is a specific case of the notion of \v{C}ech cocycles on groupoids, which appear in the next chapter. We return to our investigation of vector bundles on groupoids with some examples.

\begin{example}\label{ch3:cocyclevectorbundle1}
Consider a 1-cocycle $\ct \colon \G_1 \to \GL_n(\C)$ on $\G$. This means that $\ct$ satisfies
\[
\ct(g_1) \cdot \ct(g_2) = \ct(g_1 g_2)
\]
for composable arrows $g_1$ and $g_2$ in $\G$. We can use $\ct$ to construct a vector bundle on $\G$ as follows. Let $E = \G_0 \times \C^n$ be the trivial vector bundle. Then $t^*E = s^*E = \G_1 \times \C^n$ and $\ct$ induces an isomorphism
\[
\G_1 \times \C^n \to \G_1 \times \C^n,
\quad
(g, v) \mapsto (g, \ct(g)\cdot z).
\]
The cocycle condition on $\ct$ implies that this defines a vector bundle on $\G$. In fact, any vector bundle on an orbifold can be given by a 1-cocycle as above for a suitable choice of representing groupoid. We elaborate on this in the next example.
\end{example}

\begin{example}\label{ch3:cocyclevectorbundle2}
Let $E$ be a vector bundle on $\G$ and consider an open cover $\U = \{U_i\}$ of $\G_0$ such that $E$ is trivial over each $U_i$. We have an open cover groupoid,
\[
\G[\U] = 
\left(
\begin{array}{c}
     \bigsqcup_{i,j} t^{-1}(U_i) \cap s^{-1}(U_j) \\
     \downarrow \downarrow \\
     \bigsqcup_i U_i
\end{array}
\right),
\]
with an equivalence $\varepsilon \colon \G[\U] \to \G$. Pulling back $E$ along $\varepsilon$ gives the trivial vector bundle on $\G[\U]_0$. The isomorphism $t^*E \cong s^*E$ is pulled back along $\varepsilon_1$ to give an isomorphism $\varepsilon_1^*t^*E \cong \varepsilon_1^*s^*E$ of bundles on $\G[\U]_1$. By the functoriality of $\varepsilon$, this is the same as an isomorphism $t^*\varepsilon_0^*E \cong s^*\varepsilon_0^*E$ of trivial bundles. This produces a 1-cocycle
\[
\ct \colon \bigsqcup_{i,j} t^{-1}(U_i) \cap s^{-1}(U_j) \to \GL_n(\C).
\]
Thus, we have a vector bundle on $\G[\U]$ as in the previous example. Since $\G$ and $\G[\U]$ represent the same orbifold we see that any vector bundle on an orbifold is given by a choice of representing groupoid $\G$ and a 1-cocycle $\ct \colon \G_1 \to \GL_n(\C)$.
\end{example}

\begin{example}\label{ch3:equivariantvectorbundle}
Let $G$ be a finite or Lie group. A $G$-equivariant vector bundle on a $G$-space $X$ is a vector bundle $\pi \colon E \to X$ such that
\begin{enumerate}[label=(\roman*)]
    \item $E$ is a $G$-space and $\pi(e \cdot g) = \pi(e) \cdot g$.
    \item For any $g \in G$, the action $g \colon E_x \to E_{x\cdot g}$ is an isomorphism of vector spaces.
\end{enumerate}
Suppose we have a $G$-equivariant vector bundle on $E$. Then we can define an $X \mmod G$ action on $E$ by 
\[
e \cdot (x, g) = e \cdot g^{-1}.\]
Conversely, if we already have a $G$-action on $E$, then we obtain an $X \mmod G$ action on $E$ via
\[
e \cdot g = e \cdot (\pi(e) \cdot g, g^{-1}).
\]
This makes sense because $t(\pi(e) \cdot g, g^{-1}) = \pi(e)$. Going through the details shows that a $G$-equivariant vector bundle on $X$ is the same as a vector bundle on $X\mmod G$.
\end{example}

\begin{example}\label{ch3:representation}
Following from the previous example, a vector bundle on $\mathbb{B}G = * \mmod G$ is a linear representation of $G$: a vector bundle on $\mathbb{B}G_0 = \ast$ is a choice of vector space $V$ and adding an action by $\mathbb{B}G$ gives produces a representation.

\end{example}


\section{Line Bundles with Connection}

To warm up for gerbes with connection, we consider the definition of line bundles with connection on groupoids. Here our groupoids are assumed to be Lie groupoids. First recall the following definition:

\begin{definition}
An \emph{isomorphism} $F \colon (L_1, \nabla_1) \to (L_2, \nabla_2)$ between two line bundles with connection on a smooth manifold $X$ is a line bundle isomorphism $F \colon L_1 \to L_2$ that is compatible with the connections in the sense that
\[
F(\nabla_1(s)) = \nabla_2(F(s))
\]
for all sections $s$ of $L_1$.
\end{definition}

Now we can essentially rewrite the definition of a line bundle on a groupoid so that it includes connections:

\begin{definition}
A \emph{line bundle with connection} on a Lie groupoid $\G$ consists of the following data:
\begin{enumerate}[label=(\roman*)]
    \item A line bundle $L$ with connection $\nabla$ on $\G_0$.
    \item An isomorphism $\mu \colon (t^* L, t^*\nabla) \to (s^*L, s^*\nabla)$ of line bundles with connection that satisfies $u^* \mu = \id_{L}$ and $d_2^* \mu \circ d_0^*\mu = m^* \mu$.
\end{enumerate}
\end{definition}

\begin{example}\label{ch3:linebundleconnectionlocal}
Consider a groupoid $\G$ in which each $\G_n$ is a disjoint union of contractible spaces. Then, a line bundle $L$ on $\G_0$ must be trivial, $L \cong \G_0 \times \C$. Moreover, $t^*L$ and  $s^*L $ are isomorphic to  $\G_1 \times \C$. The isomorphism $\mu$ is therefore of the form
\[
\mu \colon \G_1 \times \C \to \G_1 \times \C, \quad (g, z) \to (g, \ct(g)\cdot z)
\]
for some function $\ct \colon \G_1 \to \C^*$. The associativity condition on $\mu$ implies that $\ct$ satisfies a cocycle condition, 
\begin{equation}\label{cocycle}
\ct(g) \cdot \ct(h) = \ct(gh).
\end{equation}
A section of $L$ is a function $\xi \colon \G_0 \to \C$ and a connection on $L$ is of the form
\[
\nabla(\xi) = d\xi + A \otimes \xi
\]
where $A$ is a 1-form on $\G_0$. The pullback connections $t^*\nabla$ and $s^*\nabla$ are given by
\[
t^*\nabla(\xi) = d\xi + (t^*A) \otimes \xi, \quad s^*\nabla (\xi) = d\xi + (s^*A) \otimes \xi.
\]
We have that
\[
\mu(t^*\nabla(\xi)) = \mu( d\xi + (t^*A) \otimes \xi) = \ct\,d\xi + (t^*A)\otimes (\ct\xi),
\]
and
\[
s^*\nabla( \mu (\xi)) = s^*\nabla(\ct\xi) = \ct\,d\xi + \xi \,d\ct + (s^*A) \otimes (\ct\xi). 
\]
Compatibility with the connection implies that these are equal, hence
\begin{equation}\label{connection}
(t^*A - s^*A) \otimes (\ct\xi) = \xi\,d\ct \quad \implies \quad t^*A - s^*A = \frac{d\ct}{\ct} = d\log \ct. 
\end{equation}
In fact, one can show that a line bundle with connection on $\G$ is equivalent (up to isomorphism) to a function $\ct \colon \G_1 \to \C^*$ and a 1-form $A$ satisfying conditions \eqref{cocycle} and \eqref{connection}. 

\end{example}

\section{Gerbes}

Line bundles on $X$ are classified up to isomorphism by their first Chern classes, which are elements of $H^2(X;\Z)$. Gerbes are a generalisation of line bundles in the sense that $p$-gerbes are classified by $H^{p+2}(X;\Z)$. Just like line bundles, gerbes have duals, pullbacks, connections and characteristic classes. Gerbes, namely 1-gerbes, were first introduced by Giraud in \cite{Giraud} and have since taken several forms. One formulation is Hitchin-Chatterjee gerbes, which on a space $M$ consists of the following data \cite{Hitchin:LecturesOnSpecialLagrangianSubmans}: 
\begin{itemize}
    \item An open cover $\{U_i\}$ of $M$.
    \item A collection of line bundles $L_{ij}$ on intersections $U_i \cap U_j$.
    \item Isomorphisms $L_{ij} \cong L_{ji}^{-1}$.
    \item Isomorphisms $L_{ij} \otimes L_{jk} \cong L_{ik}$ of line bundles on $U_i \cap U_j \cap U_k$.
    \item The following diagram commutes:
    \[
    \begin{tikzcd}
         L_i \otimes L_j \otimes L_k \arrow[d] \arrow[r] & L_i \otimes L_{jk} \arrow[d] \\
        L_{ij} \otimes L_k \arrow[r]           & L_{ijk}          
    \end{tikzcd}
    \]
    This is called the associativity condition.
\end{itemize}
This definition is analogous to the local definition of line bundles, which  are $\C^*$-valued cocycles on the intersections of open sets in an open cover. The analogy is even stronger because $(-1)$-gerbes are in fact $\C^*$-valued functions. This indicates a hierarchy of gerbes in which $p$-gerbes can be constructed from $(p-1)$-gerbes on intersections of open sets in a cover. 

The definition of a Hitchin-Chatterjee gerbe is generalised to groupoids in the following way \cite[\textsection 6]{LU04}:
\begin{definition}\label{def:gerbeongroupoid}
A \emph{gerbe} $L$ on a groupoid $\G$ is given by:
\begin{enumerate}[label=(\roman*)]
    \item A line bundle $L$ on $\G_1$.
    \item An isomorphism $i^*L \cong L^{-1}$.
    \item An isomorphism $\mu \colon d_2^*L \otimes d_0^*L \to m^*L$ satisfying the following associativity condition on fibres:
    \[
    \begin{tikzcd}
         L_g \otimes L_h \otimes L_k \arrow[d] \arrow[r] & L_g \otimes L_{hk} \arrow[d] \\
        L_{gh} \otimes L_k \arrow[r]           & L_{ghk}          
    \end{tikzcd}
    \]
\end{enumerate}
\end{definition}

Gerbes on a groupoid $\G$ are classified by $H^3(B\G;\Z)$, the degree 3 integral cohomology of its classifying space $B\G$, see \cite[Prop 6.2.2]{LU04}. In particular, the classifying space of  $X \mmod G$ is the Borel construction $X_G$, so gerbes on $X \mmod G$ are classified by $H^3(X_G;\Z)$. 

\begin{example}\label{ch3:gerbeonmanifold}
Let $M$ be a manifold with an open cover $\U = \{U_i\}$ and consider the open cover groupoid $M[\U]$. Since
\[
M[\U]_1 = \bigsqcup_{i,j} U_i \cap U_j,
\]
a line bundle $L$ on $M[\U]_1$ is a collection of line bundles $L_{ij}$ on each $U_i \cap U_j$.  The isomorphism $i^*L \cong L^{-1}$ implies that $L_{ij} \cong L_{ji}^{-1}$. The isomorphism $d_2^*L \otimes d_0^*L \cong m^*L$ implies that $L_{ij} \otimes L_{jk} \cong L_{ik}$ and that this satisfies the required associativity condition. Therefore, the line bundle $L$, together with these isomorphisms, defines a gerbe on $M$.
\end{example}

\begin{example}\label{ch3:multiplicativelinebunde}
Let $G$ be a Lie group. A gerbe on $\B G$ is the same as a multiplicative line bundle on $G$, which is a line bundle $L$ on $\B G_1 \cong G$ with a ``line bundle multiplication'':
\[
L_g \otimes L_h \cong L_{gh}.
\]
This satisfies the associativity condition given in the definition of a gerbe. Multiplicative line bundles on $G$ are in bijection with isomorphism classes of central extensions of $G$ by $\C^*$ and the degree 2 group cohomology of $G$.
\end{example}

\begin{example}\label{ch3:bundlegerbe}
\cite[\textsection 7.3]{LU04} A bundle gerbe on a space $M$ is a smooth manifold $Y$ with a surjective submersion $\pi \colon Y \to M$ and a gerbe on the groupoid $Y \!\!\fp{\pi}{\pi}\! Y \rightrightarrows Y$. In this groupoid, $s(y_1,y_2) = y_1$ and $t(y_1,y_2) = y_2$. 
\end{example}

\begin{example}\label{ch3:gerbefromacocycle}
Let $\co \colon \G_2 \to \C^*$ be a 2-cocycle, which means that $\co$ satisfies
\[
\co(a,b) \co(a,b c)^{-1} \co(a b, c) \co(b,c)^{-1} = 1 
\quad \text{for } (a,b,c) \in \G_3.
\]
Let $L = \G_1 \times \C$ be the trivial line bundle on $\G_1$. Then 
\[
d_2^*L \otimes d_0^*L \cong \G_2 \times \C \cong m^*L
\]
and we can define an isomorphism
\[
\G_2 \times \C \to \G_2 \times \C, \quad ((a,b), z) \mapsto ((a,b), \co(a,b)z).
\]
The cocycle condition on $\co$ ensures that this isomorphism satisfies the associativity condition, so we have a gerbe on $\G$. If we start with a gerbe in which the line bundle $L$ is trivial, we can go in the other direction, obtaining a 2-cocycle on $\G$. One can show that any gerbe on an orbifold can be represented by a 2-cocycle on a ``good'' representing groupoid.
\end{example}

\begin{example}\label{ex:gerbeasgroupoid}
A gerbe on $\G$ can be viewed as a groupoid $\G_L$ over $\G$. Given a gerbe $L$ on $\G$, define the groupoid
\[
\G_L =  ( L \rightrightarrows \G_0),
\]
where an element $\ell \in L$ is a morphism from $s( \pi(\ell))$ to $t(\pi(\ell))$. In other words, the source and target maps are chosen so that there is a functor $F \colon \G_L \to \G$ such that $F_1 = \pi$ and $F_0 = \id_{\G_0}$. The multiplication map is defined as 
\[
(\G_L)_2 = L \fp{t\pi}{s\pi} L \to L, \quad (\ell_1, \ell_2) \mapsto \mu(\ell_1 \otimes \ell_2),
\]
where $\mu$ is in the definition of a gerbe. The gerbe $\G_L$ will be called an \emph{$L$-twisted gerbe} and will return when we discuss twisted K-theory.
\end{example}

\begin{remark}
Definition \ref{def:gerbeongroupoid} is a local definition of bundle gerbes on groupoids. If we were starting from scratch, we may want our definition to be as follows: A bundle gerbe on a groupoid $\G$ is a bundle gerbe $L$ on $\G_0$ together with a stable isomorphism $t^*L \cong s^*L$ satisfying certain unit and associativity conditions. This approach is discussed in Appendix \ref{appendix:bundlegerbes}, where we discuss bundle gerbes and see how this leads to Definition \ref{def:gerbeongroupoid}.
\end{remark}

Let $M$ be a smooth manifold with a good open cover $\{U_i\}$, that is, each $U_i$ and finite intersection of $U_i$ is contractible. A gerbe with connection on $M$ is described by the following data:
\begin{itemize}
    \item A 2-cocycle $\co = (\co_{ijk})$.
    \item A collection $A = (A_{ij})$ of 1-forms $A_{ij} \in \Omega^1(U_i \cap U_j)$.
    \item A collection $B = (B_{i})$ of 2-forms $B_i \in \Omega^2(U_i)$.
\end{itemize}
These must satisfy the following conditions:
\begin{itemize}
    \item $A_{jk} - A_{ik} + A_{ik} = d\log g_{ijk}$ on $U_i \cap U_j \cap U_k$.
    \item $B_{j} - B_{i} = -dA_{ij}$ on $U_i \cap U_j$.
\end{itemize}
Here the gerbe itself is given by the 2-cocycle and the differential forms give the connective structure. In Example \ref{ch3:gerbefromacocycle} we saw that a gerbe on a groupoid can also be given by a 2-cocycle on $\G$. We generalise the connective structure in an analogous way \cite[\textsection 6.3]{LU04}:

\begin{definition}
A \emph{gerbe with connection} on a groupoid $\G$ is the following data:
\[
\co \colon \G_2 \to \C^*,
\quad 
A \in \Omega^1(\G_1),
\quad
B \in \Omega^2(\G_2).
\]
These must satisfy:
\begin{itemize}
    \item $\co(a,b) \co(a,b c)^{-1} \co(a b, c) \co(b,c)^{-1} = 1 $ for $(a,b,c) \in \G_3$.
    \item $d_0^*A - m^*A + d_2^*A = d\log \co$.
    \item $t^*B - s^*B = -dA$.
\end{itemize}
\end{definition}
We will see in the next section that this is precisely the description of a gerbe with connection as a Deligne 2-cocycle.

\begin{example}
If $M$ is a smooth manifold and $\U = \{U_i\}$ is an open cover then a gerbe with connection on $M[\U]$ is precisely a gerbe with connection on $M$ that is associated to the open cover $\U$.
\end{example}

\begin{remark}
This definition can be compared to the local definition of a vector bundle with connection, which consists of a 1-cocycle and connection 1-forms associated to a fixed open cover. In this situation, one does not necessarily obtain all of the vector bundles on a space by only considering a single open cover. For this to be true, one needs to consider a good open cover. It is the same for this definition of gerbes with connection. If one considers a fixed groupoid representing an orbifold, then this definition may not give all of the possible gerbes with connection unless the groupoid is a ``good'' groupoid. This ``good'' groupoid is called a Leray groupoid in the next chapter.
\end{remark}


\chapter{Deligne Cohomology and Transgression}
\label{DeligneCohomologyandTransgression}

\section{Deligne Cohomology}

We know that $p$-gerbes are classified by degree $(p+2)$ integral cohomology or equivalently degree $(p+1)$ \v{C}ech cocycles. What if we want to include connective structure in this classification? The answer is that line bundles, gerbes and higher gerbes with connective structure are classified by Deligne cohomology. This theory was first introduced in complex analytic geometry by Deligne \cite{deligne}, though we refer the reader to \cite{brylinski} for an introduction. We are motivated to study Deligne cohomology by the following facts:
\begin{itemize}
    \item There is an isomorphism between degree 1 Deligne cohomology and isomorphism classes of line bundles with connection. 
    \item The holonomy of a line bundle with connection can be expressed as a transgression map on Deligne cohomology.
\end{itemize}
These statements may be generalised to higher gerbes with connective structure, for example bundle gerbes \cite[\textsection 4]{MSstabisom}. Deligne cohomology for orbifolds is defined in Lupercio and Uribe's work \cite{LUdeligne}, and is our main reference for this section. We give a somewhat simpler formulation that is suitable for our needs.

When defining \v Cech cohomology for a topological space, it is convenient to work with a good open cover, sometimes called a Leray cover. This is an open cover in which open sets and finite intersection of open sets in the cover are contractible. Thinking of a groupoid representing an orbifold as an atlas or open cover of the orbifold, the analogous notion for orbifolds is that of a Leray groupoid:

\begin{definition}
A \emph{Leray groupoid} $\G$ is a topological groupoid such that each $\G_n$ is homeomorphic to a disjoint union of contractible spaces.
\end{definition}

\begin{proposition}
Every orbifold can be represented by a Leray groupoid.
\end{proposition}
\begin{proof}
The relevant result is \cite[Cor 1.2.5]{MPsimplicialcohomologyoforbifolds}, which says that we can represent any orbifold with an atlas which is closed under finite intersections and for which every chart in the atlas is contractible. Reviewing how one uses an orbifold atlas to build a representing groupoid, for instance in \cite[Ex 2.1.8]{amenta}, one sees that the groupoid obtained is Leray.
\end{proof}

\begin{example}\label{ex:Leraygroupoidmanifold}
If $\U$ is a good open cover of a manifold $M$ then $M[\U]$ is a Leray groupoid.
\end{example}

\begin{example}\label{ex:LerayGroupoidOfCirclemodC2}
Consider the cyclic group $C_2 = \langle r \mid r^2=1 \rangle$ acting on $S^1 \subseteq \C$ by $z \cdot r = -z$. The groupoid $S^1 \mmod C_2$ is not a Leray groupoid but an appropriate choice of open cover will produce one. Let $\U = \{U_1, U_2, U_3, U_4\}$ be the open cover in Figure \ref{fig:opencoverforexample}. Note that $U_1 \cdot r = U_3$ and $U_2 \cdot r = U_4$. Consider the groupoid
\[
\bigsqcup_{i,j} U_i \cap U_j \times C_2 \rightrightarrows \bigsqcup_{i} U_i
\]
where for example $r$ induces a morphism from $x \in U_1$ to $x\cdot r \in U_3$. This is a Leray groupoid that is equivalent to $S^1 \mmod C_2$.
\begin{figure}
    \centering
    \begin{tikzpicture}

    
    \draw (0,0) circle [radius=2];
    
    \draw [ultra thick, mypink, -)] (0:2) arc (0:55:2);
    \draw [ultra thick, myblue, (-)] (45-10:2) arc (35:145:2);
    \draw [ultra thick, mypink, (-)] (135-10:2) arc (125:235:2);
    \draw [ultra thick, myblue, (-)] (225-10:2) arc (215:325:2);
    \draw [ultra thick, mypink, (-] (315-10:2) arc (305:360:2);
    
    \node [right] at (0:2) {$U_1$};
    \node [above] at (90:2) {$U_2$};
    \node [left] at (180:2) {$U_3$};
    \node [below] at (270:2) {$U_4$};
    
    \draw [thick, ->] (0:1) arc (0:180:1) node [midway, below] {$r$};
    
\end{tikzpicture}
    \caption{The open cover in Example \ref{ex:LerayGroupoidOfCirclemodC2}.}
    \label{fig:opencoverforexample}
\end{figure}
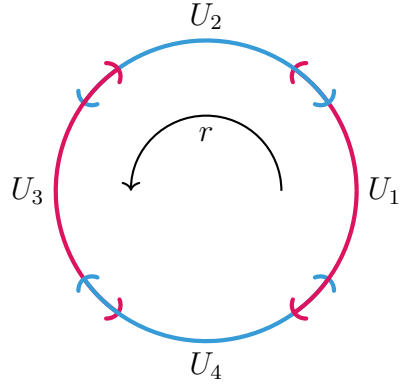
\end{example}

Generally speaking, if $\U$ is a good open cover of $\G_0$, then $\G[\U]$ need not be a Leray groupoid - the cover also needs to be compatible with the morphism structure of $\G$. For $X \mmod G$, the solution is to choose a good open cover $\{U_i\}$ of $X$ with the additional property that for any $i$ and $g \in G$ there exists $j$ such that $U_i \cdot g = U_j$. A Leray groupoid can then be constructed as in Example \ref{ex:LerayGroupoidOfCirclemodC2} \cite[12]{LU:introtogerbesonorbifolds}.

We will define the Deligne cohomology of a groupoid. The Deligne cohomology of an orbifold will then be the Deligne cohomology of a Leray groupoid representing it. The Deligne complex is the following complex of sheaves:
\begin{equation}
\underline{\C}^* \xrightarrow{\,\, d\log \,\,} \Omega^1 \xrightarrow{\,\,\,\, d\,\,\,\,} \cdots \xrightarrow{\,\,\,\, d\,\,\,\,} \Omega^{p}
\end{equation}
This is to be recognised as the truncated de Rham complex with 0-forms replaced with $\C^*$-valued functions. 

For a groupoid $\G$, consider the double complex 
\[
K^{k,\ell} =
\begin{cases}
C^\infty(\G_\ell, \C^*), & k = 0, \ell \geq 0, \\
\Omega^k(\G_\ell), & 0 \leq k \leq p, \ell \geq 0, \\
0, & \text{otherwise}.
\end{cases}
\]
Here $C^\infty(\G_n, \C^*)$ denotes smooth $\C^*$-valued maps on $\G_n$ and $\Omega^k(\G_n)$ denotes complex valued $k$-forms on $\G_n$. The coboundary maps are as denoted the following diagram:
\begin{center}
\begin{tikzcd}
\vdots          & \vdots         &         & \vdots   \\
C^\infty(\G_2, \C^*) \arrow[r, "d\log"] \arrow[u, "\delta"] & \Omega^1(\G_2) \arrow[r, "d"] \arrow[u, "\delta"] & \cdots \arrow[r, "d"]  & \Omega^p(\G_2) \arrow[u, "\delta"] \\
C^\infty(\G_1, \C^*) \arrow[r, "d\log"] \arrow[u, "\delta"] & \Omega^1(\G_1) \arrow[r, "d"] \arrow[u, "\delta"] & \cdots \arrow[r, "d"]  & \Omega^p(\G_1) \arrow[u, "\delta"] \\
C^\infty(\G_0, \C^*) \arrow[r, "d\log"] \arrow[u, "\delta"] & \Omega^1(\G_0) \arrow[r, "d"] \arrow[u, "\delta"] & \cdots \arrow[r, "d"]  & \Omega^p(\G_0) \arrow[u, "\delta"]
\end{tikzcd}
\end{center}
 The map $d\log$ takes $\co \colon \G_n \to \C^*$ to the 1-form $d\log \co = d\co/\co$. The remaining horizontal maps are given by the exterior derivative of the de Rham complex. The vertical maps are given by
\[
\delta \colon \Omega^k(\G_n) \to \Omega^k(\G_{n+1}), 
\quad
\omega \mapsto \sum_{i=0}^n (-1)^i d_0^*\omega.
\]
In the first column we must write this multiplicatively;
\[
\delta \colon C^\infty(\G_n, \C^*) \to C^\infty(\G_{n+1}, \C^*),
\quad
\co \mapsto \prod_{i=0}^n (\co\circ d_i)^{(-1)^i}.
\]
This complex is called the \v Cech de Rham double complex because the columns are \v Cech complexes and the rows are truncated de Rham complexes with 0-forms replaces with $\C^*$ valued functions. From here we can obtain a cochain complex $K^\bullet$ by:
\[
K^n = \bigoplus_{k + \ell = n} K^{k,\ell}
\]
The coboundary map $D \colon K^n \to K^{n+1}$ is given by $D = d + (-1)^p \delta$. For example, an element of $K^1$ is a pair $(\co, \omega)$ with $\co \in  C^\infty(\G_1, \C^*)$ and $\omega \in \Omega^1(\G_0)$. This is mapped into $K^2$ by
\[
D(\co, \omega) = (\delta \co, d\log \co - \delta \omega, d\omega).
\]
Then because $d^2 = 0$, $\delta^2 = 0$ and $d\delta = \delta d$ we have
\begin{align*}
D^2(\co, \omega) = (\delta^2\co, d\log \delta \co - \delta d\log \co + \delta^2 \omega, d^2\log \co - d\delta \omega + \delta d \omega , d^2 \omega) = 0.
\end{align*} 
In general, we have $D^2 =0$ as we would hope. This is a common way of obtaining a chain complex from a double complex and is often called the \emph{total complex} of the double complex \cite[\textsection 1.2]{brylinski}.

\begin{definition}
\emph{Deligne cohomology} of a groupoid $\G$, denoted by
\[
H^*(\G, \underline{\C}^* \to \Omega^1 \to \cdots \to \Omega^{p}),
\]
is the cohomology of the cochain complex $K^\bullet$ defined above. Deligne cohomology of an orbifold is the Deligne cohomology of a Leray groupoid representing it.
\end{definition}

\begin{example}\label{ch3:deligneformanifolds}
Let $\U = \{U_i\}$ is a Leray cover of a smooth manifold $M$. The Deligne cohomology of $M[\U]$ matches the Deligne cohomology of the manifold $M$. For example, a Deligne 1-cocycle is $(\co, A)$ where $\co = (\co_{ij})$ is a 1-cocycle and $A = (A_i)$ is a collection of 1-forms on each $U_i$. This defines a line bundle with connection on $M$.
\end{example}

\begin{proposition}
Let $\G$ be a Leray groupoid. There is an isomorphism between the group of line bundles with connection on $\G$ and $H^1(\G, \C^* \to \Omega^1)$.
\end{proposition}
\begin{proof}
This is a modified version of the proof for smooth manifolds in \cite[Thm 2.2.11]{brylinski}. Let $L$ be a line bundle with connection on $\G$. As $\G_0$ is a disjoint union of contractible spaces, $L$ has is a non-vanishing section $\sigma$. Let $A \in \Omega^1(\G_0)$ be the 1-form satisfying $\nabla(\sigma) = A \otimes \sigma$. Letting $\mu$ be the isomorphism $t^*L \to s^*L$, the maps $\mu \circ \sigma \circ t$ and $\sigma \circ s$ are two non-vanishing sections of $s^*L$. This means there is a map $\co \colon \G_1 \to \C^*$ such that $\mu \circ \sigma \circ t = \co \cdot (\sigma \circ s)$. This implies that
\[
t^*A = \frac{d(\mu \circ \sigma \circ t)}{\mu \circ \sigma \circ t} = \frac{d(\co (\sigma \circ s))}{\co (\sigma \circ s)} = \frac{d\co (\sigma \circ s) + \co d(\sigma \circ s)}{\sigma \circ s} = \frac{d\co}{\co} + s^*A.
\]
In other words, $t^*A - s^*A = d\log \co$, and so $(\co,A)$ is a Deligne 1-cocycle. Suppose we made a different choice of section n $\sigma' = f\cdot \sigma$ and cocycle $\co'$ such that $\mu \circ \sigma' \circ t = \co'\cdot (\sigma' \circ s)$. Then
\[
\mu \circ \sigma' \circ t = (f \circ t)\cdot (\mu \circ \sigma \circ t) = (f \circ t)\cdot \co \cdot (\sigma \circ s) = \frac{f \circ t}{f \circ s} \cdot \co \cdot (\sigma' \circ s),
\]
which implies that $\co' \co^{-1} = (f\circ t) (f \circ s)^{-1}$. Moreover, letting $A' = d\sigma' / \sigma'$, we have that
\[
A'  = \frac{d\sigma'}{\sigma'} = \frac{d(f\sigma)}{f\sigma} = \frac{df}{f} + \frac{d\sigma}{\sigma} = d \log f + A.
\]
Therefore, $(\co', A')-(\co, A) = (\delta(f), d\log f)$ and the two Deligne 1-cocycles are cohomologous. Thus, we have a well-defined homomorphism from the isomorphism classes of line bundles with connection to $H^1(\G, \C^* \to \Omega^1)$.

For injectivity, assume that $(\co, A)$ is trivial. In this case, we can choose a section $\sigma$ such that $\co = 1$ and $A = 0$. This means our line bundle with connection is trivial.

For surjectivity, start with a Deligne 1-cocycle $(\co, A)$. Choose $L = \G_0 \times \C^*$, so that $t^*L$ and $s^*L$ are both equal to $\G_1 \times \C^*$. Then $\co$ determines an isomorphism $t^*L \to s^*L$ and $A$ defines a connection on $L$ by $\nabla(\sigma) = A \otimes \sigma$.
\end{proof}

\begin{proposition}
Let $\G$ be a Leray groupoid. There is an isomorphism between the group of gerbes with connection on $\G$ and $H^2(\G, \C^* \to \Omega^1 \to \Omega^2)$.
\end{proposition}
\begin{proof}
We've somewhat cheated by defining gerbes with connection as Deligne 1-cocycles. For bundle gerbes on manifolds the result is proved in \cite[\textsection 4]{MSstabisom} and this can be generalised to the groupoid case in a similar way to the previous proof.
\end{proof}

We want to define 2-gerbes with connection on groupoids. At this point, the simplest way to do this is via Deligne cohomology.

\begin{definition}
A \emph{2-gerbe with connection} on a groupoid $\G$ is a Deligne 3-cocycle. In other words, a 2-gerbe consists of
\[
\co \colon \G_3 \to \C^*, \quad \omega^1 \in \Omega^1(\G_2), \quad \omega^2 \in \Omega^2(\G_1), \quad \omega^4 \in \Omega^3(\G_0).
\]
These must satisfy the following conditions:
\begin{gather*}
    \co(b, c, d) \cdot  \co(ab, c, d)^{-1} \cdot \co(a, bc, d) \cdot \co(a, b, cd)^{-1} \cdot \co(a, b, c) = 1, \\
    d_0^*\omega^1 - d_1^*\omega^1 + d_2^*\omega^1 - d_3^*\omega^1 = d\log \co, \\
    d_0^*\omega^2 - d_1^*\omega^2 + d_2^*\omega^2 = -d\omega^1, \\
    t^* \omega^3 - s^*\omega^3 = d\omega^2.
\end{gather*}
\end{definition}

\noindent Without any connective structure a 2-gerbe on $\G$ is just a 3-cocycle $\co \colon \G_3 \to \C^*$.

A 3-cocycle $\co \colon \G_3 \to \C^*$ is normalised if $\co(1, a, b) = \co(a, 1, b) = \co(a,b,1)$ for all suitable $a,b \in \G_1$. It is usually convenient to work with normalised cocycles, and the following implies that we always can.

\begin{lemma}\label{lem:normalisedcocycles}
Every 3-cocycle $\co \colon \G_3 \to \C^*$ is cohomologous to a normalised cocycle.
\end{lemma}

A clear proof is provided in \cite[Lemma 2.2]{normalised3cocycleresult}.

\section{Transgression}

Holonomy of a line bundle with connection on $X$ is a map $LX \to \C^*$ that associates a complex number to each loop in $X$. In other words, holonomy is a homomorphism
\[
\{ \text{Line bundles with connection on $X$} \} \longrightarrow \{ \text{Maps $LX \to \C^*$} \}. 
\]
The left-hand side can be interpreted as Deligne 1-cocycles on $X$ and the right-hand side as Deligne $0$-cocycles on $LX$. Holonomy can thus be interpreted as a transgression map on Deligne cohomology \cite[\textsection 6.1]{brylinski}:
\[
H^1(X, \C^* \to \Omega^1) \longrightarrow H^0(LX, \C^*).
\]
Moving one step higher in the gerbe hierarchy, the holonomy of gerbes with connective structure is a homomorphism:
\[
\{ \text{Gerbes with connection on $X$} \} \longrightarrow \{ \text{Line bundles with connection on $LX$} \}. 
\]
This may also be interpreted as a transgression map on Deligne cohomology, namely:
\[
H^2(X, \C^* \to \Omega^1 \to \Omega^2) \longrightarrow H^1(LX, \C^* \to \Omega^1).
\]
One starts to see the pattern. Holonomy of $p$-gerbes with connective structure is a transgression map
\[
H^p(X, \C^* \to \Omega^1 \to \cdots \to \Omega^p) \longrightarrow H^{p-1}(LX, \C^* \to \Omega^1 \to \cdots \to \Omega^{p-1}).
\]
This means that we start with a $p$-gerbe with connective structure on $X$ and obtain a $(p-1)$-gerbe with connective structure on the loop space. Gawedzki \cite{Gawedkzi} was the first to define the transgression map in the gerbe case, though in the context of quantum field theory. An expository reference is \cite[\textsection 6.5]{brylinski}, where Brylinski additionally defines the transgression map for general $p$-gerbes with connective structure. Gomi and Terashima \cite{gomiterashima} define a general transgression map which includes $p$-gerbes as a particular case.

We are interested in the holonomy of 2-gerbes with connection on groupoids. In this case, we need to consider the loop groupoid rather than the loop space, considering a transgression map
\[
\tau \colon  H^3(\G, \C^* \to \Omega^1 \to \Omega^2 \to \Omega^3) \longrightarrow H^{2}(\L\G, \C^* \to \Omega^1 \to \Omega^{2}),
\]
for a Leray groupoid $\G$ and its loop groupoid $\L\G$. We owe this generalisation to Lupercio and Uribe \cite{LUhol} who explicitly define holonomy for line bundles and gerbes on orbifolds in addition to providing a formulation of the general case. In this section, we will use their formulation to write down the transgression map in the 2-gerbe case and show that it is the map we're after. 

\subsection{Preamble}

Before defining the transgression map, some notation needs to be established. Recall that a loop in $\G$ is represented by a smooth functor $\Phi \colon \S^1_{\U} \to \G$ where $\U$ is an admissible cover of $S^1$ with vertices, say,
\[
0 < \alpha_0 < \alpha_1 < \dotsb < \alpha_n \leq 1
\]
and a sufficiently small choice of $\epsilon$. Choosing a smaller $\epsilon$ results in a refinement of $\U$ and the restriction of $\Phi$ to this refinement gives an equivalent loop in $\L\G$. Define
\[
I_i = [ \alpha_{i-1}, \alpha_i], \quad i \in \{1, \dotsc, n+1\},
\]
where $\alpha_{n+1} = 1 + \alpha_0$. Next, define the groupoid
\[
\S^1_{\alpha} = 
\left(
\begin{array}{c}
    \bigsqcup_{i,j} (W_i \cap W_j)  \times \Z \\
     \downdownarrows \\
     \bigsqcup_i W_i
\end{array}
\right)
\stext{where}
W_i = \bigsqcup_{n \in \Z} (I_i + n).
\]
One should compare $\S^1_\alpha$ to the groupoid $\S^1_{\U}$ - we think of $\S^1_\alpha$ as being the limiting case of the usual $\S^1_{\U}$ as we take finer and finer refinements of $\U$ by letting $\epsilon$ approach 0. As an element of $\L\G$, which is a colimit of all such refinements, $\Phi$ can be viewed as a functor $\Phi \colon \S^1_{\alpha} \to \G$.

Define, for $i \in \{1, ..., n+1\}$, the paths
\begin{equation}\label{eq:Phii}
\varphi_i \colon I_i \to \G_0, \quad \varphi_i  := \Phi_0|_{I_i}.
\end{equation}
Also define
\begin{equation}\label{eq:Phi}
\varphi \colon \{\alpha_0, \dotsc , \alpha_n\} \to \G_1, \quad 
\varphi(\alpha_i) = 
\begin{cases}
    \Phi_1(\alpha_i, 0),        &i \in \{1, ... ,n\}, \\
    \Phi_1(\alpha_{n+1}, -1),   & i=n+1.
\end{cases}
\end{equation}
Then we have
\[
s(\varphi(\alpha_i)) = \varphi_i(\alpha_i) \stext{and} t(\varphi(\alpha_i)) = \varphi_{i+1}(\alpha_i)
\]
with $t(\varphi(\alpha_{n+1})) = \varphi_1(\alpha_0)$. The reasoning for these definitions is as follows. We want to restrict our attention from $\S^1_\alpha$ to the sub-groupoid
$\mathcal{I}$ defined by
\[
\mathcal{I}_0 = \bigsqcup_i I_i, \quad \text{and,}
\] 
\[
\mathcal{I}_1 = (I_1+1) \cap I_{n+1} \times \!\{-1\} \,\sqcup\, \bigsqcup_{i,j} \,(I_i \cap I_j) \times \!\{0\}.
\]
The inclusion $\mathcal{I} \hookrightarrow \S^1_\alpha$ is an equivalence and so the restriction of $\Phi$ from $\S^1_{\alpha}$ to $\mathcal{I}$, which is given by the $\varphi_i$ and $\varphi$ defined above, gives an equivalent generalised map via the diagram
\begin{center}
\begin{tikzpicture}


\node (L) at (0,1)  {$\S^1$};
\node (R) at (5,1) {$\G$};
\node (U) at (2.5, 1.75) {$\S^1_\alpha$};
\node (D) at (2.5, 0.25) {$\mathcal{I}$};

\draw [<-] (L) -- (U);
\draw [->] (U) -- (R) node [midway, yshift=7pt] {\small$\Phi$};
\draw [<-] (L) -- (D);
\draw [->] (D) -- (R) node [below right, xshift=-5pt, yshift=2pt, midway] {\small$\varphi$};
\draw [->] (D) -- (U);

\end{tikzpicture}
\end{center}
In this section, whenever there is a loop $\Phi$ associated to a collection of vertices $\alpha_i$, we will assume that the intervals $I_i = [\alpha_{i-1}, \alpha_i]$ and the functions defined in \eqref{eq:Phii} and \eqref{eq:Phi} are implicitly defined. We will continue using the convention that the relevant lower-case Greek letter will be used for the maps obtained from loops in $\L\G$, which are denoted by upper-case Greek letters. Morphisms $\Lambda \colon \Phi \to \Psi$ will be viewed as a collection of maps $ \Lambda_i \colon I_i \to \G_1$ such that $\Lambda_i(x)$ is a morphism from $\varphi_i(x)$ to $\psi_i(x)$ and the following diagram commutes:
\begin{center}
    \begin{tikzcd}[column sep =15mm]
        \varphi_i(\alpha_i) \arrow[r, "\Lambda_i(\alpha_i)"] \arrow[d, "\varphi(\alpha_i)", swap] & \psi_i(\alpha_i) \arrow[d, "\psi(\alpha_i)"]\\
        \varphi_{i+1}(\alpha_i) \arrow[r, "\Lambda_{i+1}(\alpha_i)"] & \psi_{i+1}(\alpha_i) 
    \end{tikzcd}
\end{center}
For convenience, $\Lambda_{n+2}(\alpha_{n+1}):= \Lambda_1(\alpha_0)$. 

\subsection{The Transgression Map}

The transgression map is first defined as a map on Deligne cochains:
\[
\tau \colon (\co, \omega^1, \omega^2, \omega^3) \longmapsto (h, \theta^1, \theta^2).
\]
Here we have
\[
\co \colon \G_3 \to \C^*, \quad \omega^1 \in \Omega^1(\G_2), \quad \omega^2 \in \Omega^2(\G_1), \quad \omega^4 \in \Omega^3(\G_0),
\]
and we obtain
\[
\ct \colon \L\G_2 \to \C^*, \quad \theta^1 \in \Omega^1(\L\G_1), \quad \theta_2 \in \Omega^2(\L\G_0). 
\]
We do not care about the connective structure on the gerbe that we obtain, so we will ignore $\theta^1$ and $\theta^2$. It turns out that only $\co$ and $\omega^1$ are required to define $\ct$. Moreover, when the obtained gerbe is restricted to the inertia groupoid, we will see that the dependence on $\omega^1$ vanishes. 

Let $\U$ be an admissible cover of $S^1$ and let $(\Lambda, \Omega) \in \L\G(\U)_2$ where
\[
\Phi \xrightarrow{\,\, \Lambda \,\,} \Psi \xrightarrow{\,\, \Omega \,\,} \Xi.
\]
As in the previous section, we obtain maps 
\[
\varphi_i, \psi_i, \xi_i \colon I_i \to \G_0, 
\quad
\phi, \psi, \xi \colon \{\alpha_0, ..., \alpha_{n+1}\} \to \G_1
\stext{and}
\Lambda_i, \Omega_i \colon I_i \to \G_1
\]
such that the following diagram commutes in $\G$ for each $\alpha_i$:
\[
\begin{tikzcd}
\varphi_i(\alpha_i) \arrow[d, "\varphi(\alpha_i)"] \arrow[r, "\Lambda_i(\alpha_i)"] & \psi_i(\alpha_i) \arrow[d, "\psi(\alpha_i)"] \arrow[r, "\Omega_i(\alpha_i)"] & \xi_i(\alpha_i) \arrow[d, "\xi(\alpha_i)"] \\
\varphi_{i+1}(\alpha_i) \arrow[r, "\Lambda_{i+1}(\alpha_i)"]           & \psi_{i+1}(\alpha_i) \arrow[r, "\Omega_{i+1}(\alpha_i)"]           & \xi_{i+1}(\alpha_i)          
\end{tikzcd}
\]
Now define:
\[
\ct(\Lambda, \Omega) 
    =\exp\left(  \sum_{i=1}^n \int_{I_i} (\Lambda_i, \Omega_i)^*\omega^1 \right)   \cdot
        \prod_{i=1}^n \frac{ \co( \varphi(\alpha_i), \Lambda_{i+1}(\alpha_i), \Omega_{i+1}(\alpha_i)) \cdot \co(\Lambda_i(\alpha_i), \Omega_i(\alpha_i), \xi(\alpha_i)) }{\co( \Lambda_i(\alpha_i), \psi(\alpha_i), \Omega_{i+1}(\alpha_i))}
\]
We need $\ct$ to satisfy three properties:
\begin{enumerate}
    \item $\ct$ doesn't depend on the choice of refinement $\U$.
    \item If $\co$ and $\omega^1$ are components of a Deligne 3-cocycle, then $\ct$ is a 2-cocycle.
    \item If $\co$ and $\omega^1$ are components of a Deligne 2-coboundary, then $\ct$ is a 1-coboundary.
\end{enumerate}
These three results together imply that $\tau$ gives us a well-defined gerbe on $\L\G$. Before starting the proofs, we will introduce the notation
\[
\ct^i_1(\Lambda, \Omega) = \exp\left(\int_{I_i} (\Lambda_i, \Omega_i)^*\omega^1 \right)
\]
and
\[
\ct^i_2(\Lambda, \Omega) = 
\frac{ \co( \varphi(\alpha_i), \Lambda_{i+1}(\alpha_i), \Omega_{i+1}(\alpha_i)) \cdot \co(\Lambda_i(\alpha_i), \Omega_i(\alpha_i), \xi(\alpha_i)) }{\co( \Lambda_i(\alpha_i), \psi(\alpha_i), \Omega_{i+1}(\alpha_i))}
\]
so that 
\[
\ct(\Lambda, \Omega) = \prod_{i=1}^n \ct_1^i(\Lambda, \Omega) \cdot \ct_2^i(\Lambda, \Omega).
\]

\begin{lemma}
If $\co$ is a normalised cocycle, then $\ct$ does not depend on the choice of admissible cover $\U$ used to represent $\Lambda$ and $\Omega$.
\end{lemma}

\begin{proof}

Let $\U$ and $\widetilde \U$ be two admissible covers of the circle. By considering a common refinement of $\U$ and $\widetilde \U$ we can assume that $\widetilde \U$ is a refinement of $\U$. Moreover, we may assume that $\widetilde \U$ is obtained from $\U$ by adding a vertex $\widetilde \alpha$ to the interval $I_{i}$, dividing the interval into two new intervals $J_1$ and $J_2$. This is because $\widetilde \U$ can always be obtained from $\U$ by adding a finite number of vertices. 

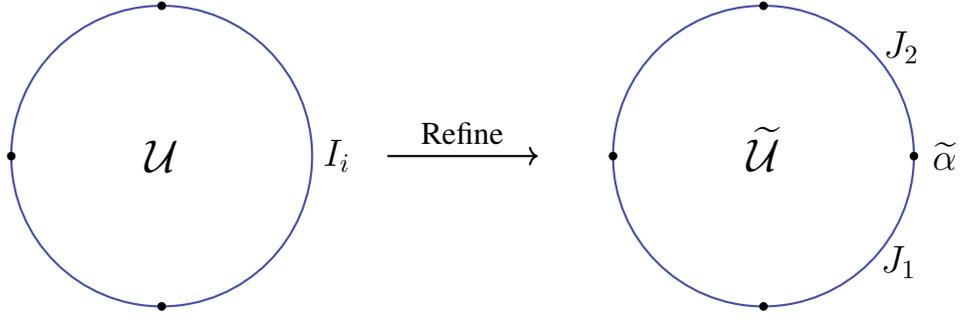
\begin{figure}
    \centering
    \begin{tikzpicture}

    \def \pt {0.05};

    \begin{scope}[shift={(4,2)}]
    
        \draw [mydarkblue, thick] (0,0) circle [radius=2];
        \draw [fill] (90:2) circle [radius=\pt];
        \draw [fill] (180:2) circle [radius=\pt];
        \draw [fill] (270:2) circle [radius=\pt];
        
        \node [right] at (2,0) {\large $I_i$};
        
        \node at (0,0) {\Large $\U$};
    
    \end{scope}
    
    \begin{scope}[shift={(12,2)}]
    
        \draw [thick, mydarkblue] (0,0) circle [radius=2];
        \draw [fill] (90:2) circle [radius=\pt];
        \draw [fill] (180:2) circle [radius=\pt];
        \draw [fill] (270:2) circle [radius=\pt];
        \draw [fill] (0:2) circle [radius=\pt] node [right, xshift=1mm] {\large $\widetilde{\alpha}$};
        
        \node [right] at (45:2.05) {\large$J_2$};
        \node [right] at (-45:2) {\large$J_1$};
        
        \node at (0,0.1) {\Large $\widetilde{\U}$};
    
    \end{scope}
    
    \draw [thick, ->] (7,2) -- (9,2) node [midway, above] {Refine};

\end{tikzpicture}
    \caption{$\U$ is refined by adding an additional vertex $\widetilde{\alpha}$ to the interval $I_i$}
    \label{fig:my_label}
\end{figure}

Let $\rho \colon S^1_{\widetilde{\U}} \to \S^1_\U$ be the natural inclusion map. Let $(\Lambda, \Omega) \in \L\G(\U)_2$, where
\[
\varphi \xrightarrow{\,\, \Lambda \,\,} \psi \xrightarrow{\,\, \Omega \,\,} \xi.
\]
The equivalent maps in $L\G(\widetilde \U)$, which we notate using tildes, are obtained by composing with $\rho$, 
\[
\widetilde \Lambda = \Lambda \circ \rho_1, \quad \widetilde \varphi_i = \varphi_i \circ \rho_0,
\]
and similarly for the remaining maps. These equivalent maps in $\L\G(\widetilde{\U})$ are equal to their counterpart in $\L\G(\U)$ everywhere except on the divided intervals $J_1$ and $J_2$. Use the indices $j_1$ and $j_2$ to represent the components of the maps on $J_1$ and $J_2$ respectively. We have
\[
\widetilde \Lambda_{j_1} = \Lambda_i |_{J_1} \quad \text{and} \quad \widetilde \Omega_{j_1} = \Omega_i |_{J_1}
\]
and similarly for $J_2$. Therefore,
\[
\int_{J_1} (\widetilde \Lambda_{j_1}, \widetilde \Omega_{j_1})^*\omega^1 + \int_{J_2} (\widetilde \Lambda_{j_2}, \widetilde\Omega_{j_2})^*\omega^1 = \int_{I_i} (\Lambda, \Omega)^* \omega^1
\]
and so the exponential part of $\ct$ remains the same. The remaining part of the equation is changed by the addition of the new vertex $\widetilde{\alpha}$. The morphisms $\widetilde{\varphi_i}(\widetilde{\alpha}), \widetilde{\psi_i}(\widetilde{\alpha})$ and $\widetilde{\xi}(\widetilde{\alpha})$ are the identity morphisms of $\varphi_i(\widetilde{\alpha}), \psi_i(\widetilde{\alpha})$ and $\xi_i(\widetilde{\alpha})$, respectively. Since $\co$ is a normalised cocycle, this implies that the $\widetilde \alpha$ terms are equal to 1. Therefore, $\ct$ has the same value whether the loops are represented by $\U$ or $\widetilde{\U}$.
\end{proof}

\begin{lemma}
Suppose that $\co$ and $\omega^1$ satisfy, for $(a,b,c,d) \in \G_4$\textnormal{:}
\begin{gather*}
    \co(b, c, d) \cdot  \co(ab, c, d)^{-1} \cdot \co(a, bc, d) \cdot \co(a, b, cd)^{-1} \cdot \co(a, b, c) = 1, \\
    d_0^*\omega^1 - d_1^*\omega^1 + d_2^*\omega^1 - d_3^*\omega^1 = d\log \co.
\end{gather*}
Then $\ct$ is a 2-cocycle, that is, if $(\Lambda, \Omega, \Theta) \in L\G(\U)_3$, then
\[
\frac{ \ct(\Omega, \Theta) \cdot \ct(\Lambda, \Omega \cdot \Theta) }{ \ct(\Lambda\cdot\Omega, \Theta) \cdot \ct(\Lambda, \Omega)} = 1.
\]
\end{lemma}
\begin{proof}
This is more or less a direct calculation. Assume that
\[
\varphi_1 \xrightarrow{\,\, \Lambda \,\,} \varphi_2 \xrightarrow{\,\, \Omega \,\,} \varphi_3 \xrightarrow{\,\, \Theta \,\,} \varphi_4.
\]
First, we see that
\begingroup
\addtolength{\jot}{7pt}
\setlength{\abovedisplayskip}{0pt}
\setlength{\abovedisplayshortskip}{0pt}
\begin{multline*}
\ct_1^i(\Omega, \Theta)\cdot \ct_1^i(\Lambda\cdot\Omega, \Theta)^{-1} \cdot \ct_1^i(\Lambda, \Omega \cdot \Theta)  \cdot \ct_1^i(\Lambda, \Omega)^{-1} \\
\begin{aligned}
    &= \exp\left(  \int_{I_i} (\Omega_i, \Theta_i)^*\omega^1 - ((\Lambda \cdot \Omega)_i, \Theta_i)^*\omega^1 + (\Lambda_i, (\Omega \cdot \Theta)_i)^*\omega^1 - (\Lambda_i, \Omega_i)^*\omega^1\right) \\
    &= \exp\left( \int_{I_i} (\Lambda_i, \Omega_i, \Theta_i)^* ( d_0^*\omega^1 - d_1^*\omega^1 + d_2^*\omega^1 - d_3^*\omega^1 ) \right) \\
    &= \exp\left( \int_{I_i} (\Lambda_i, \Omega_i, \Theta_i)^* d\log \co  \right) \\
    &= \frac{
            \co(\Lambda_i(\alpha_i), \Omega_i(\alpha_i), \Theta_i(\alpha_i))
            }{
            \co(\Lambda_i(\alpha_{i-1}), \Omega_i(\alpha_{i-1}), \Theta_i(\alpha_{i-1}))}.
\end{aligned}
\end{multline*}
\endgroup
Now calculate what happens to the $\ct_2^i$ part of the equation - the steps will be explained after the calculation. To simplify notation, we omit writing the arguments of the inputs of $\co$.
\begingroup
\addtolength{\jot}{7pt}
\setlength{\abovedisplayskip}{0pt}
\setlength{\abovedisplayshortskip}{0pt}
\begin{multline*}
\ct_2^i(\Omega, \Theta)\cdot \ct_2^i(\Lambda\cdot\Omega, \Theta)^{-1} \cdot \ct_2^i(\Lambda, \Omega \cdot \Theta)  \cdot \ct_2^i(\Lambda, \Omega)^{-1} \\
\begin{aligned}
        &= \frac{\co(\Lambda_{i+1}, \Omega_{i+1}, \Theta_{i+1})}{\co(\varphi_1 \cdot \Lambda_{i+1}, \Omega_{i+1}, \Theta_{i+1})} \times
            \frac{\co(\Lambda_i, \Omega_i, \Theta_i\cdot \varphi_4)}{\co(\Lambda_i, \Omega_i, \Theta_i)}  \times
            \frac{\co(\varphi_2, \Omega_{i+1}, \Theta_{i+1})}{\co(\Omega_i, \varphi_3, \Theta_{i+1})}\\
    &\phantom{=~}
            \times \frac{\co(\Lambda_i, \varphi_2, \Omega_{i+1})}{\co(\Lambda_i, \Omega_i, \varphi_3)}
            \times \frac{\co((\Lambda \cdot \Omega)_i, \varphi_3, \Theta_{i+1})}{\co(\Lambda_i, \varphi_2, (\Omega\cdot\Theta)_{i+1})} \\
        &= \frac{ \co(\Lambda_{i+1}, \Omega_{i+1}, \Theta_{i+1}) }{ \co(\Lambda_i, \Omega_i, \Theta_i)} 
        \times \frac{\co(\Lambda_i, \varphi_2\cdot\Omega_{i+1}, \Theta_{i+1})}{\co(\Lambda_i, \Omega_I \cdot \varphi_3, \Theta_{i+1})} \\
        &= \frac{\co(\Lambda_{i+1}, \Omega_{i+1}, \Theta_{i+1})}{\co(\Lambda_i, \Omega_i, \Theta_i)}.
\end{aligned}
\end{multline*}
\endgroup
In the first step we use the cocycle condition on the tuples
\[
(\varphi_1, \Lambda_{i+1}, \Omega_{i+1}, \Theta_{i+1})
\quad \text{and} \quad
(\Lambda_i, \Omega_i, \Theta_i, \varphi_4).
\]
In the next step we use the cocycle condition on
\[
(\Lambda_i, \varphi_2, \Omega_{i+1}, \Theta_{i+1})
\quad \text{and} \quad
(\Lambda_i, \Omega_i, \varphi_3, \Theta_{i+1}).
\]
Using that $\varphi_2 \cdot \Omega_{i+1} = \Omega_i \cdot \varphi_3$ we simplify to get the final quotient. Putting the two calculations together we conclude that
\begin{multline*}
    \frac{ \ct(\Omega, \Theta) \cdot \ct(\Lambda, \Omega \cdot \Theta) }{ \ct(\Lambda\cdot\Omega, \Theta) \cdot \ct(\Lambda, \Omega)} \\
    = \prod_{i=1}^n \frac{
            \co(\Lambda_i(\alpha_i), \Omega_i(\alpha_i), \Theta_i(\alpha_i))
            }{
            \co(\Lambda_i(\alpha_{i-1}), \Omega_i(\alpha_{i-1}), \Theta_i(\alpha_{i-1}))}
            \cdot
            \frac{\co(\Lambda_{i+1}(\alpha_i), \Omega_{i+1}(\alpha_i), \Theta_{i+1}(\alpha_i))}{\co(\Lambda_i(\alpha_i), \Omega_i(\alpha_i), \Theta_i(\alpha_i))}
            = 1. \raisebox{-2em}{\qedhere}
\end{multline*}
\end{proof}

\begin{lemma}
If there exists $\cth \colon \G_2 \to \C^*$ and $\mu \in \Omega^1(\G_1)$ such that
\[
\co = \delta \cth
\quad \text{and} \quad
\omega^1 = d\log \cth - \delta \mu
\]
then there exists $F \colon \L\G(\U)_1 \to \C^*$ such that $\ct = \delta F$.
\end{lemma}
\begin{proof}
The 1-cocycle $F$ is obtained from the transgression map of \cite[\textsection 4]{LUhol}, which sends a Deligne 2-cocycle on $\G$ to a Deligne 1-cocycle on $\L\G$. This is one step below the transgression map we've been considering. Letting $\tau_1$ and $\tau_2$ be the gerbe and 2-gerbe transgression maps, we obtain the following commutative diagram: 
\[
\begin{tikzcd}
(\cth, \mu) \arrow[r, maps to, "D"] \arrow[d, maps to, "\tau_1", swap] & (\co, \omega^1) \arrow[d, maps to, "\tau_2"] \\
F \arrow[r, maps to, "D"]                    & \ct       
\end{tikzcd}
\]
It is left to the reader to check that this diagram commutes.
\end{proof}

\subsection{Manifolds}

Let $M$ be a smooth manifold with open cover $\U = \{U_i\}$. We will consider the transgression map in the case where $X = M[\U]$ just as the authors of \cite[\textsection 4.1]{LUhol} do for gerbes. Here we start with the following data:
\begin{itemize}
    \item A 3-cocycle $\co \colon M[\U]_3 \to \C^*$,  which is equivalent to a collection of maps $\co = (\co_{ijkl})$ where
    \[
    \co_{ijkl}  \colon U_i \cap U_j \cap U_k \cap U_\ell \to \C^*.
    \]
    \item A 1-form $\omega^1 \in \Omega^1(M[\U]_2)$,  which is equivalent to a collection of 1-forms $\omega^1 = (\omega^1_{ijk})$ with
    \[
    \omega^1_{ijk} \in \Omega^1(U_i \cap U_j \cap U_k).
    \]
\end{itemize}
These satisfy the familiar conditions. Given an admissible cover $\U$, a loop $\varphi$ in $L\G(\U)$ is a collection of maps 
\[
\varphi_i \colon I_i \to M[\U]_0.
\]
These maps are continuous and $M[\U]_0$ is a disjoint union of open sets, so there exists $\kappa_i$ such that $\varphi_i(I_i) \subset U_{\kappa_i}$. Given another loop $\psi$ with $\psi_i(I_i) \subset U_{\lambda_i}$, a morphism $\Lambda \colon \varphi \to \psi$ is a collection of maps $\Lambda_i \colon I_i \to M[\U]_1$ with $\Lambda_i(t) = \varphi_i(t) = \psi_i(t) \in M$. This implies that
$
\Lambda_i(I_i) \subset U_{\kappa_i} \cap U_{\lambda_i}.
$

Consider a pair of arrows $\varphi \xrightarrow{\,\,\Lambda\,\,} \psi \xrightarrow{\,\, \Omega\,\,} \xi$, where 
\[
\varphi_i(I_i) \subset U_{\kappa_i}, 
\quad
\psi_i(I_i) \subset U_{\lambda_i},
\quad \text{and} \quad
\xi_i(I_i) \subset U_{\mu_i}.
\]
We have that $\Lambda_i(I_i) \subset U_{\kappa_i} \cap U_{\lambda_i}$ and $\Omega_i(I_i) \subset U_{\lambda_i} \cap U_{\mu_i}$. The transgression map in this context becomes
\[
h(\Lambda, \Omega)
    = \exp\left( \sum_{i=1}^n \int_{I_i} (\Lambda_i, \Omega_i)^* \omega^1_{\kappa_i \lambda_i \mu_i} \right)
    \cdot \prod_{i=1}^n
    \frac{
        \co_{\kappa_i \kappa_{i+1} \lambda_{i+1} \mu_{i+1}}(\Lambda_i(\alpha_i)) \cdot 
        \co_{\kappa_i \lambda_i \mu_i \mu_{i+1}}(\Lambda_i(\alpha_i))
        }{ 
        \co_{\kappa_i \lambda_i \lambda_{i+1} \mu_{i+1}}(\Lambda_i(\alpha_i))
        }.
\]
We have recovered the transgression map in \cite{brylinski} and \cite{gomiterashima}, the latter only after restricting to the gerbe case.

\subsection{Gerbes on the Inertia Groupoid}

Recall that we have an inclusion functor $\iota \colon \Lambda \G \to \L\G$. If we restrict our attention to the inertia groupoid, then we only need to consider loops associated to the trivial cover. For an arrow $(g,h)$ in the inertia groupoid, the corresponding morphism $\Lambda = \iota_1(g,h)$ is the constant map
\[
\Lambda \colon [0,1] \to \G_1, \quad  \Lambda(t) = h.
\]
This map should be written with a subscript, namely $\Lambda_1$. Since we only have a single interval, $I_1 = [0,1]$, this is not necessary.
A pair of morphisms $(\Lambda, \Omega)$ is also constant, and so
\[
(\Lambda, \Omega)^*\omega^1 = 0 
\quad \implies \quad
\exp\left(\int_0^1  (\Lambda, \Omega)^*\omega^1 \right) = 1.
\]
Therefore, if $\Lambda = \iota_1(g,h)$ and $\Omega = \iota_1(h^{-1}gh, k)$ are the arrows 
\[
\varphi := \iota_0(g) \xrightarrow{\,\, \Lambda\,\,} \psi := \iota_0(h^{-1}gh) \xrightarrow{\,\,\Omega\,\,} \xi := \iota_0(k^{-1}h^{-1}ghk),
\]
then the transgression map becomes
\[
\ct(\Lambda, \Omega) = 
\frac{ \co(g, h, k) \cdot \co(h, k,k^{-1}h^{-1}ghk) }{ \co(h, h^{-1}gh, k) }.
\]

\begin{remark}
Restricting the transgression map to the action groupoid completely removed the contribution of the connection on the 2-gerbe we started with. We now have a transgression from 2-gerbes on $\G$ to gerbes on $\Lambda \G$, without any talk of connective structure.
\end{remark}

\begin{remark}
This formula appears in several places in the literature, for instance in the context of the twisted Drinfeld double \cite{Willerton:Drinfeld}, gauge theory \cite{Dijkgraaf:Witten:TopologicalGaugeTheories} and, importantly, Moonshine \cite{Mason:unpublished}.
\end{remark}

\subsection{Transgression for \texorpdfstring{$X \mmod G$}{X/G}}

Consider a 2-gerbe on $X \mmod G$ which is given by a 3-cocycle 
\[
\alpha \colon X \times G^3 \to \C^*
\]
with trivial connective structure. The only relevant loops in $X \mmod G$ are those associated to the trivial cover. Choosing three loops $\R \mmod \Z \to X \mmod G$ leads to maps $\varphi, \psi, \xi \colon [0,1] \to X$ and elements $g_1, g_2, g_3 \in G$ such that
\[
\varphi(1) = \varphi(0)\cdot g_1, \quad \psi(1)=\psi(0)\cdot g_2, \stext{and} \xi(1)=\xi(0)\cdot g_3.
\]
Morphisms $\varphi \xrightarrow{\,\, \Lambda \,\,} \psi \xrightarrow{\,\, \Omega \,\,} \xi$ are maps $\Lambda, \Omega \colon [0,1] \rightarrow X \times G$ connecting $\varphi$ to $\psi$ and $\psi$ to $\xi$, respectively. At this point we don't need to talk about $\psi$ and $\xi$ anymore as they're completely determined by $\varphi, g_1$ and $g_2$. There exists $h_1, h_2 \in G$ such that
\[
\Lambda\bigl([0,1]\bigr) \subseteq X \times \{h_1 \} \stext{and} \Omega([0,1]) \subseteq X \times \{h_2\}.
\]
Explicitly this means that $\varphi(x) \cdot h_1 = \psi(x)$ and $\psi(x)\cdot h_2 = \xi(x)$. Naturality implies that $g_2 = h_1^{-1}g_1 h_1$ and $g_3 = h_2^{-1}g_2h_2 = h_2^{-1} h_1^{-1} g_1 h_1h_2$, so we don't need to talk about $g_2$ and $g_3$ anymore either.

The transgression map gives us a gerbe $\theta$ on $\L(X \mmod G) = \left( \bigsqcup\nolimits_g \P_g \right) \mmod G$,
\[
\theta \colon \left( \bigsqcup\nolimits_g \P_g \right) \times G^2 \to \C^*, 
\]\[
\theta(\varphi, h_1, h_2) = 
\frac{
    \alpha(\varphi(0), g_1, h_1, h_2) \cdot \alpha(\varphi(0), h_1, h_2, h_2^{-1} h_1^{-1} g_1 h_1h_2)
}{
    \alpha(\varphi(0), h_1, h_1^{-1} g_1 h_1, h_2)
}.
\]
Restricting to the case where $h_1,h_2 \in C_{g_1}$, this simplifies to
\[
\theta(\varphi, h_1, h_2) = 
\frac{
    \alpha(\varphi(0), g_1, h_1, h_2) \cdot \alpha(\varphi(0), h_1, h_2, g_1)
}{
    \alpha(\varphi(0), h_1, g_1, h_2)
}.
\]
This formula defines a gerbe on 
\[
\bigsqcup_{[g]} \P_g \mmod C_g,
\]
which is another description of $\L(X \mmod G)$. If we restrict to the inertia groupoid, then we obtain formulas that are equivalent to those in the previous section.

\chapter{Twisted Equivariant Tate K-Theory}
\label{TETKTheory}

\section{Twisted K-Theory}

\subsection{Twisted K-Theory of Spaces}

Ordinary K-theory is twisted by an element of degree 3 integral cohomology, or equivalently a $\C^*$-valued \v{C}ech 2-cocycle. When this is a torsion class, the twisted K-theory can be described geometrically in terms of twisted vector bundles.

\begin{definition}\label{def:twistedvb}
Let $\{U_i\}$ be an open cover of a space $X$ and let $\alpha = (\alpha_{ijk})$ be a $\C^*$-valued 2-cocycle associated to this open cover. An \emph{$\alpha$-twisted vector bundle} of rank $n$ is a collection $(\co_{ij})$ of continuous maps $\co_{ij} \colon U_{i} \cap U_j \to  \GL_n(\C)$ such that 
\[
    \co_{ii} = 1,                 \quad 
    \co_{ij} = \co_{ji}^{-1},       \quad  \text{and} \quad 
    \co_{ij} \cdot \co_{jk} = \alpha_{ijk} \cdot \co_{ik}.
\]
\end{definition}

\noindent We see straight away that if $\alpha$ is the trivial cocycle, then this becomes the definition of an ordinary vector bundle.

\begin{definition}
Let $\co = (\co_{ij})$ and $\ct = (\ct_{ij})$ be two $\alpha$-twisted vector bundles. An \emph{isomorphism} from $\co$ to $\ct$ is a collection $f = (f_i)$ of maps $f_i \colon U_i \to  \GL_n(\C)$ such that
\[
f_i \cdot \co_{ij} = \ct_{ij} \cdot f_j.
\]
\end{definition}

Now we can define $\alpha$-twisted K-theory.

\begin{definition}
Let $\alpha$ be a 2-cocycle associated to a good open cover $\{U_i\}$ of $X$. Then \emph{$\alpha$-twisted K-theory}, ${}^{\alpha}K(X)$, is the Grothendieck group of $\alpha$-twisted vector bundles on $X$. 
\end{definition}

If $\co = (\co_{ij})$ is a twisted line bundle then the condition
\[
\co_{ij} \cdot \co_{ik}^{-1} \cdot  \co_{jk} = \alpha_{ijk},
\]
means that $\alpha$ is a coboundary, hence is equal 1 in \v{C}ech cohomology. Moreover, taking the $n$th exterior product of a rank $n$ twisted vector bundle results in an $\alpha^n$-twisted line bundle, which indicates that $\alpha^n$ is trivial. To elaborate, a matrix $A\in  \GL_n(\C)$ is a linear map $\C^n \to \C^n$, and after applying the $n$th exterior product functor we obtain,
\[
\begin{tikzcd}
\C^n \arrow[r, "A"] \arrow[d] & \R^n \arrow[d] \\
 \C \arrow[r, "\Lambda^n(A)"]           & \C          
\end{tikzcd}
\]
since $\Lambda^n(\C^n) \cong \C$. If $A = \lambda I$, then $\Lambda^n(A) = \lambda^n I$. Applying this to $\co$ results in an $\alpha^n$-twisted cocycle $\tilde \co_{ij}:= \Lambda^n(\co_{ij})$, which takes values in $\C^*$. Functoriality implies that the twisted cocycle condition becomes
\[
\tilde \co_{ij} \tilde \co_{ik}^{-1} \tilde \co_{jk} = \alpha_{ijk}^n.
\]
Therefore, $\alpha^n$ is a coboundary. This means that if $\co$ is an $\alpha$-twisted vector bundle of rank $n$, then $\alpha^n = 1$. In particular, if $\alpha$ is of infinite order, then there are no $\alpha$-twisted vector bundles and the twisted K-theory is meaningless. To twist K-theory by an arbitrary cohomology class, one must employ some infinite dimensional analogue of a twisted vector bundles, for instance projective bundles of infinite dimensional Hilbert spaces \cite{AStwisted}. 

\subsection{Twisted Orbifold K-Theory}

In the previous section, a 2-cocycle was used to define a notion of twisted vector bundles. This is generalised to twisted vector bundles on orbifolds \cite[\textsection 7]{LU04}:

\begin{definition}
Let $\G$ be a topological groupoid and let $L$ be a gerbe on $\G$. An \emph{$L$-twisted vector bundle} on $\G$ consists of the following data:
\begin{enumerate}[label=(\roman*)]
    \item A vector bundle $\pi \colon E \to \G_0$.
    \item An isomorphism $\mu \colon L \otimes t^*E \to s^*E$ of vector bundles over $\G_1$ that is compatible with the gerbe multiplication, that is, for $z \xrightarrow{g} y \xrightarrow{h} x$, the following commutes:
    \[
    \begin{tikzcd}
         L_g \otimes L_h \otimes E_{x} \arrow[d] \arrow[r] & L_g \otimes E_y \arrow[d] \\
        L_{gh} \otimes E_x \arrow[r]           & E_z
    \end{tikzcd}
    \]
\end{enumerate}
\end{definition}

If the gerbe is trivial than we recover the definition of a vector bundle on a groupoid.  Moreover, if $\G$ represents a manifold, this definition reduces to the definition of a twisted vector bundle introduced in the previous section. We elaborate on this in the following example.

\begin{example}
Let $M$ be a manifold with a good open cover $\U = \{U_i\}$. A gerbe $L$ on $M[\U]$ is uniquely determined, up to isomorphism, by a 2-cocycle $\alpha \colon M[\U]_2 \to \C^*$. Let $E$ be a vector bundle on $M[\U]$. As $\U$ is a good open cover, there exists a non-vanishing section $\sigma \colon M[\U]_0 \to E$. Then, $\mu \circ (\xi \otimes (\sigma \circ t))$ and $\sigma \circ s$ are non-vanishing sections of $s^*E$, where $\xi$ is a non-vanishing section of $L$. These two sections differ by a map $\co \colon M[\U]_1 \to \C^*$. The reader can check that the compatibility condition implies that $\co$ is an $\alpha$-twisted vector bundle on $M$ with respect to the open cover $\U$.
\end{example}

\begin{example}\label{ch4:twistedvectorbundleonBG}
As we've seen, a gerbe on $\mathbb{B}G$ is a multiplicative line bundle on $G$ or equivalently a group 2-cocycle $\theta \colon G \times G \to \C^*$. This gives rise to a central extension
\[
1 \to \C^* \to \widetilde G \to G \to 1
\]
where $\widetilde G =  G \times \C^*$ with multiplication
\[
(g,z) \cdot (h,w) = (gh, \theta(g,h) \cdot z \cdot w).
\]
This is how multiplicative line bundles on $G$ classify central extensions of $G$ by $\C^*$. Let $L$ be the gerbe on $\mathbb{B}G$ associated to the cocycle $\theta \colon G \times G \to \C^*$. An $L$-twisted vector bundle is:
\begin{itemize}
\item A vector bundle on $\mathbb{B}G_0 = *$. This is just a vector space $V$.
\item An isomorphism of vector bundles over $G$:
\[
L \otimes t^*V \cong G \times V \to G \times V \cong s^*V.
\]
This is given by a function $\rho \colon G \to  \GL(V)$. To be compatible with the gerbe multiplication we must have
\[
\rho(g) \cdot \rho(h) = \theta(g,h) \rho(gh).
\]
This is called a \emph{twisted representation} of $G$ and is in particular a projective representation.
\end{itemize}
We obtain linear representation of $\widetilde G$ defined as
\[
\tilde \rho \colon \widetilde G \to  \GL(V), \quad (g, z) \mapsto z \cdot \rho(g).
\]
This works because
\[
\tilde\rho( (g,z) \cdot (h,w) ) 
= \tilde\rho(gh , \theta(g,h) \cdot z \cdot w ) 
= \theta(g,h) \cdot z \cdot w \cdot \tilde\rho(gh) 
= z\cdot w \cdot \tilde\rho(g) \tilde\rho(h) 
= \tilde\rho(g,z) \cdot \tilde\rho(h,w).
\]
In this way we find that an $L$-twisted vector bundle on $\mathbb{B}G$ is a representation of $\widetilde G$ in which $\C^* \subseteq \widetilde G$ acts by complex multiplication.
\end{example}

\begin{example}\label{ch4:twistedvbviacocycle}
Suppose that our gerbe $L$ is given by a 2-cocycle $\alpha \colon \G_2 \to \C^*$. A twisted vector bundle in this case is a vector bundle $\pi \colon E \to \G_0$ together with an isomorphism $\mu \colon t^*E \to s^*E$ satisfying
\[
d_2^* \mu \circ d_0^* \mu = \theta \cdot m^* \mu.
\]
Compare this to the twisted cocycle condition in Definition \ref{def:twistedvb}
\end{example}

\begin{definition}
Let $E$ and $F$ be $L$-twisted vector bundles on a groupoid $\G$. An \emph{isomorphism} from $E$ to $F$ is a vector bundle isomorphism $f \colon E \to F$ that commutes with the $L$-action, that is, the following diagram commutes:
\[
    \begin{tikzcd}
         L \otimes t^*E \arrow[d, "\id \otimes \,t^*f", swap] \arrow[r] & s^*E \arrow[d, "s^*f"] \\
        L \otimes t^*E  \arrow[r]           & s^*F
    \end{tikzcd}
\]
\end{definition}

\noindent Given two $L$-twisted bundles $E$ and $F$ we obtain an $L$-twisted bundle $E \oplus F$ with action given by the isomorphisms
\[
L \otimes t^*(E \oplus F) \cong (L \otimes t^*E) \oplus (L \otimes t^*E) \cong s^*E \oplus s^*F \cong s^*(E \oplus F).
\]
Therefore, the group of isomorphism classes of $L$-twisted vector bundles on $\G$ is well-defined and we can define the corresponding K-theory.

\begin{definition}
Let $\G$ be a topological groupoid and let $L$ be a gerbe on $\G$. Then the \emph{$L$-twisted $K$-theory} of $\G$, ${}^L K(\G)$, is the Grothendieck group of $L$-twisted vector bundles on $\G$.
\end{definition}

\begin{example}\label{ch4:twistedktheoryofBG}
In Example \ref{ch4:twistedvectorbundleonBG} we saw that an $L$-twisted vector bundle on $\mathbb{B}G$ is the same as a representation of the corresponding central extension $\widetilde G$ in which the central $\C^*$ acts by complex multiplication. Therefore, ${}^L K(\mathbb{B}G)$ is the sub-ring of $R(\widetilde G)$ generated by these representations.
\end{example}

In defining equivariant structures on $X$, we built the equivariant structure into a groupoid $X \mmod G$. We hope to achieve the same thing for twisted structures; in Example \ref{ex:gerbeasgroupoid} we defined a twisted gerbe $\G_L$. The question is whether or not a twisted vector bundle on $\G$ is the same as an ordinary vector bundle on $\G_L$. The answer is no, but it is almost true. Recall that a vector bundle on $\G_L$ is a vector bundle $p \colon E \to (\G_L)_0 = \G_0$ together with an action of $(G_L)_1$, which is given by a map $E \fp{p}{t\pi} L \to E$ that is linear on the fibres of $E$. This is equivalent to a map $t^*E \oplus L \to s^*E$, but this is not a map of vector bundles unless it is bilinear on fibres. We say that the $L$-action is bilinear on fibres if the corresponding map $t^*E \oplus L \to s^*E$ is bilinear on fibres.

\begin{proposition}\label{prop:vbontwistedgrp}
An $L$-twisted vector bundle on $\G$ is equivalent to a vector bundle on $\G_L$ such that the $L$-action is bilinear on fibres.
\end{proposition}

\begin{proof}
Let $\pi \colon E \to \G_0$ be a vector bundle. We show that an isomorphism $\mu \colon L \otimes t^*E \to s^*E$ of vector bundles over $\G_1$ satisfying the compatibility conditions is equivalent to an action of $L$ on $E$. Indeed, the universal property of tensor products implies that $\mu$ is equivalent to a map
\[
L \oplus t^*E \to s^*E
\]
that is bilinear on fibres. There are homeomorphisms
\[
L \oplus t^*E \cong t^*E \oplus L \cong E \fp{p}{\pi t} L.
\]
The required action is the composition
\[
 E \fp{p}{\pi t}L \cong L \oplus t^*E  \to s^*E \to E,
\]
where $s^*E \to E$ is the projection onto $E$. Since $L \oplus t^*E \to E$ is bilinear on fibres, so is this composition. The reader can check that this satisfies the required properties and that the reverse construction works as well.
\end{proof}

\begin{corollary}
The twisted K-theory ${}^L K(\G)$ is isomorphic to the sub-ring of $K(\G_L)$ generated by vector bundles on $\G_L$ in which the $\G_L$-action is bilinear on fibres.
\end{corollary}

\section{Tate K-theory}

As a cohomology theory, Tate K-theory is obtained by trying to make sense of the functor
\begin{equation}\label{ch4:tatektheory}
X \mapsto K^*_{S^1}(LX).
\end{equation}
As it stands, this is not a cohomology theory. The reason is that loop spaces are not compatible with the Mayer-Vietoris principle. Specifically,
\[
X = U \cup V \quad \notimplies \quad LX = LU \cup LV,
\]
and so (\ref{ch4:tatektheory}) does not inherit all of the properties of $S^1$-equivariant K-theory.

The way of mending this issue is by restricting to the subspace $X \subseteq LX$ of constant loops. Consider the following completed version of $S^1$-equivariant K-theory from \cite{kitchloo:morava:thomprospectra}:
\[\label{KS1hat}
\widehat K_{S^1}(X) := K_{S^1}(X^{S^1}) \otimes_{Z[q^{\pm}]} \Z \ls{q}.
\]
This is obtained by completing the coefficient ring 
\[
K_{S^1}(*) = R(S^1) \cong \Z[q^{\pm}]
\]
at positive powers of $q$, which denotes the standard representation of $S^1$:
\[
q \colon S^1 \to U(1), \quad t \mapsto e^{2\pi it}.
\]
We have the following localisation property \cite[Thm 5.1, Cor 5.3]{kitchloo:morava:thomprospectra}:

\begin{theorem}
Let $X$ be an $S^1$-CW-complex. The inclusion $j \colon X^{S^1} \hookrightarrow X$ of fixed points defines an isomorphism
\[
j^* \colon \widehat K_{S^1}(X) \to \widehat K_{S^1}(X^{S^1}).
\]
\end{theorem}

We will not concern ourselves too much with what an $S^1$-CW-complex is. The important thing is that $LX$ is an $S^1$-CW-complex, so we can apply the theorem to obtain the isomorphism
\begin{equation}\label{ch4:loopspacektheory}
\widehat{K}_{S^1}(LX) \cong \widehat{K}_{S^1}(X).
\end{equation}
The true definition of Tate K-theory is obtained:

\begin{definition}
\emph{Tate $K$-theory} is the cohomology theory defined by
\[
K_{\Tate}(X) := \widehat{K}_{S^1}(LX) \cong \widehat{K}_{S^1}(X).
\]
\end{definition}

This definition can be simplified a bit. The constant loops are rotation invariant, so we know that \cite[Prop 2.2]{segal1968}
\[
K_{S^1}(X) \cong K(X) \otimes_{\Z} R(S^1) \cong K(X)[q^{\pm}].
\]
If we assume that $K_{S^1}(X)$ is finitely generated as an $R(S^1)$-module then we have that
\[
\widehat{K}_{S^1}(X) \cong K(X)[q^{\pm }] \otimes_{\Z[q^{\pm}]} Z \ls{q} \cong K(X) \ls{q}.
\]
One might not think that this is a particularly interesting cohomology theory. The importance of Tate K-theory comes when we view it is an \emph{elliptic cohomology theory}.

\begin{definition}
\cite[Def 1.2]{Lurie:SurveyofEllipticCohomology} An \emph{elliptic cohomology theory} consists of the following data:
\begin{enumerate}[label=(\roman*)]
    \item A commutative ring $R$.
    \item An elliptic curve $E$ over $R$.
    \item A multiplicative cohomology theory $A$ which is even and weakly periodic.
    \item Isomorphisms $A(*) \cong R$ and $\widehat E \cong \operatorname{Spf} A(\C P^\infty)$. Here $\widehat E$ is the formal completion of the elliptic curve $E$ along its identity section and $\operatorname{Spf}$ denotes the formal spectrum.
\end{enumerate}
\end{definition}

\noindent The relevant elliptic curve for Tate K-theory is the Tate curve over $R=\Z\ls{q}$. As a sanity check, we can note that
\[
K_{\Tate}(*) = K(*)\ls{q} = \Z\ls{q}
\]
and because K-theory is even and weakly periodic, so is Tate K-theory. We will consider the relationship between Tate K-theory and the Tate curve to be outside the scope of this thesis.

\section{Equivariant Tate K-theory}

To build equivariant Tate K-theory we consider the equivariant analogues of all the structures considered in the previous section:
\begin{align*}
    X  \quad &\leadsto \quad X \mmod G \\
    LX \quad &\leadsto \quad \L(X \mmod G) \\
    X \subseteq LX \quad &\leadsto \quad \Lambda(X \mmod G) \subseteq \L(X \mmod G).
\end{align*}
We will first investigate how $S^1$ acts on the loop groupoid $\L(X \mmod G)$. An object of $\L(X \mmod G)$ is a map $\varphi \colon \R \to X$ such that $\varphi(t+1) = \varphi(t)\cdot g$. This is acted upon by $h \in C_g$ via
\[
(\varphi \cdot h)(t) = \varphi(t) \cdot h.
\]
The real numbers also act via rotations; if $z \in \R$, then
\[
(\varphi \cdot z)(t) = \varphi(t+z).
\]
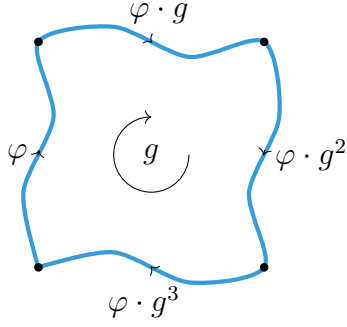
\begin{figure}
    \centering
    \begin{tikzpicture}

    \def \pt {0.05};

    \draw [->] (4,3.5) arc (0:-270:0.5);
    \node at (3.5,3.5) {$g$};

    \draw [ultra thick, myblue] plot [smooth] coordinates
            {(2,2) (1.8,3) (2.2,4) (2,5) (3, 5.2) (4,4.8) (5,5) (5.2, 4) (4.8,3) (5,2) (4, 1.8) (3, 2.2) (2,2)};
            
    \draw [fill] (2,2) circle [radius=\pt];
    \draw [fill] (2,5) circle [radius=\pt];
    \draw [fill] (5,5) circle [radius=\pt];
    \draw [fill] (5,2) circle [radius=\pt];  
            
    \def \e {0.01};
    \draw [->] (2-\e, 3.5-\e) -- (2+\e, 3.5+\e+0.05);
    \draw [->] (3.5-\e, 5+\e) -- (3.5+\e, 5-\e);
    \draw [->] (5+\e, 3.5+\e+0.05) -- (5-\e, 3.5-\e);
    \draw [->] (3.5+\e, 2-\e) -- (3.5-\e, 2+\e);

    \node [left] at (2,3.5) {$\varphi$};
    \node [above] at (3.5+0.1,5+0.1) {$\varphi\cdot g$};
    \node [right] at (5,3.5) {$\varphi\cdot g^2$};
    \node [below] at (3.5-0.1,2-0.1) {$\varphi\cdot g^3$};

\end{tikzpicture}
    \caption{The element $g \in C_g$ acts by rotating loops in $\P_g$.}
    \label{ch4:loopinpgrotatedbyg}
\end{figure}
In fact, $g \in C_g$ acts in the same way as $1 \in \R$:
\[
(\varphi \cdot g)(t) = \varphi(t)\cdot g = \varphi(t+1).
\]
In other words, acting by $g$ and then by $-1$ is the same as not doing anything at all. This motivates the definition of the following group:
\[
\Lambda_g  := \frac{ C_g \times \R }{\langle (g, -1) \rangle }.
\]
Modify the loop groupoid to incorporate the $S^1$-action as follows:
\[
\L_{S^1}(X \mmod G) = \bigsqcup_{[g]} \mathcal{P}_g \mmod \Lambda_g, \quad 
\Lambda_{S^1}(X \mmod G) = \bigsqcup_{[g]} X^g \mmod \Lambda_g.
\]
The K-theory of these groupoids will be analogous to the $S^1$-equivariant K-theory of $LX$ and $X$ considered in the previous section.

We could now follow the formalism of the previous section and define equivariant Tate K-theory as
\begin{equation}\label{ch4:equivtatektheory1}
K(\Lambda_{S^1}(X \mmod G)) \otimes_{\Z[q^{\pm}]} Z\ls{q}.
\end{equation}
However, we wish for a more concrete definition. First note that
\begin{equation}\label{ch4:quasielliptic}
K(\Lambda_{S^1}(X \mmod G)) = \bigoplus_{[g]} K( X^g \mmod \Lambda_g) = \bigoplus_{[g]} K_{\Lambda_g}(X^g).
\end{equation}
The quotient map $C_g \times \R / |g|\Z \to \Lambda_g$ induces
\begin{equation}\label{eq:quotient}
K_{\Lambda_g}(X^g) \to K_{C_g \times \R / |g|\Z}(X^g).
\end{equation}
The image of this map is generated by vector bundles on $X^g$ with both a $C_g$ and $S^1$-action, where the action of $g \in C_g$ agrees with the action of $1 \in S^1$. The circle acts trivially on $X^g$, so
\[
K_{C_g \times \R / |g|\Z}(X^g) \cong  K_{C_g}(X^g) \otimes_{\Z} R(\R / \abs{g}\Z) \cong K_{C_g}(X^g)[q^{\pm \frac{1}{\abs{g}}}],
\]
where $q^{\frac{1}{\abs{g}}} \colon t \mapsto e^{2\pi i t/ \abs{g}}$ is the standard representation of the ``long circle'' $\R / \abs{g}\Z$. Using \eqref{eq:quotient}, we therefore identify $K_{\Lambda_g}(X^g)$ with the sub-ring of 
\[
K_{C_g}(X^g) \otimes_{\Z} R(S^1) \cong K_{C_g}(X^g)[q^{\pm \frac{1}{|g|}}]
\]
consisting of by polynomials whose $q^{j / \abs{g}}$ coefficient is a virtual $C_g$-equivariant vector bundle on $X^g$ such that $g$ acts as multiplication by $\zeta^j_{\abs{g}}$. Here $\zeta_{\abs{g}} = e^{2\pi i / \abs{g}}$ is the primitive $\abs{g}$th root of unity. Adopting the language of \cite{ganter2013}, we will call a condition of this form the \emph{rotation condition} with respect to $g$. Under the assumption that $K_{C_g}(X^g)$ is finite over $R(C_g)$, we have
\[
K_{C_g}(X^g)[q^{\pm \frac{1}{|g|}}] \otimes_{\Z[q^{\pm]}} Z\ls{q} \cong K_{C_g}\ls{q^{\frac{1}{|g|}}}.
\]
We are led to the following definition \cite{ganterstringy}.

\begin{definition}
Let $X$ be a space acted on by a finite group $G$. \emph{Equivariant Tate K-theory}, $K_{\Tate}(X \mmod G)$, is the sub-ring of 
\[
\bigoplus\nolimits_{[g]} K_{C_g}(X^g)\ls{q^{\frac{1}{|g|}}},
\]
in which the $[g]$-summand is the Grothendieck group of Laurent series
\[
\sum_{j > n} V_j q^{\frac{j}{|g|}}, \quad n \in \Z,
\]
satisfying the following rotation condition with respect to $g$: for each $j$ the coefficient $V_j$ is a $C_g$-equivariant vector bundle on $X^g$ such that $g$ acts as multiplication by $\zeta^j_{|g|}$.

\end{definition}

\begin{example}
As we would hope, if $G$ is the trivial group, then we recover the definition of non-equivariant Tate K-theory, $K(X)\ls{q}$.
\end{example}

\begin{example}
The other extreme is when $X$ is a point. In this case we obtain a sub-ring of
\[
\bigoplus_{[g]} R(C_g)\ls{q^{\frac{1}{|g|}}}.
\]
Similarly, if $X$ is a trivial $G$-space, then we have
\[
\bigoplus_{[g]} \left[K(X) \otimes_{\Z} R(C_g) \right]\ls{ q^{\frac{1}{|g|}}}.
\]
\end{example}

\begin{remark}\label{rem:quasiellipticcohomology}
The K-theory in (\ref{ch4:quasielliptic}) is in fact the definition of Zhen Huan's quasi-elliptic cohomology $\text{QEll}_G(X)$ for finite $G$ \cite{Huan:thesis}. This is an intermediate cohomology theory which, although is not itself elliptic, is more computationally convenient then Tate K-theory. Moreover, it is defined for compact Lie groups, not just finite groups. Its relationship with equivariant Tate K-theory is
\[
K_{\Tate}(X \mmod G) \cong \text{QEll}_G (X) \otimes_{\Z[q^{\pm}]} Z\ls{q}.
\]
\end{remark}

\section{Twisted Equivariant Tate K-Theory}

\subsection{Twisting by a 3-Cocycle}

In this section we finally arrive at a twisted version of equivariant Tate K-theory. As before we begin with the groupoid $X \mmod G$ but now we also need a 2-gerbe to do the twisting. Consider a 2-gerbe on $X \mmod G$ given by a 3-cocycle
\[
\alpha \colon X \times G^3 \to \C^*.
\]
As discussed in Chapter \ref{StructuresOnGroupoids}, not every 2-gerbe on $X \mmod G$ can be described in this form. Our work in Chapter \ref{DeligneCohomologyandTransgression} allows us to use $\alpha$ to build a gerbe $\theta$ on $\Lambda(X \mmod G)$. When restricted to $X^g \mmod C_g$ this is the 2-cocycle $\theta \colon X^g \times C_g \times C_g \to \C^*$ with
\[
\theta(x, h, k) = \frac{\alpha(x, g, h, k) \cdot \alpha(x, h, k, g)}{\alpha(x, h, g, k)}.
\]
This gerbe is used to define a twisted inertia groupoid:
\[
\Lambda_{\theta}(X \mmod G) := \bigsqcup_{[g]} \, (X^g \mmod C_g)_\theta = \bigsqcup_{[g]} 
  \left( X^g \times C_g \times \C^* \rightrightarrows X^g \right).
\]
Here $(x, h, z)$ is a morphism from $x$ to $x \cdot h$, just as in the untwisted inertia groupoid. The twist comes into play when morphisms are composed; composition in $\Lambda_{\theta}(X \mmod G)$ is defined to be
\[
(x, h, z) \cdot (x\cdot h, k, w) = (x, hk, \theta(h,k)zw).
\]
This groupoid was previously denoted by $\G_L$, defined in Example \ref{ex:gerbeasgroupoid}.

\begin{example}
If $X$ is a point then $\theta$ is a collection of group 2-cocycles, one on each centraliser $C_g$. This gives us a collection of central extensions $\widetilde{C}_g$. Explicitly, $\widetilde{C}_g = C_g \times \C^*$ with multiplication 
\[
(h, z) \cdot (k, w) = (hk, \theta(h,k)zw).
\] 
We therefore find that
\[
\Lambda_{\theta}(\mathbb{B}G) = \bigsqcup_{[g]} \mathbb{B}\widetilde{C}_g.
\]
\end{example}

In the untwisted case we had powers of $q^{\frac{1}{|g|}}$ in our definition because $g$ acted on loops by a rotation that matched the action of $1 \in \R / |g|\Z$. For the twisted case we will instead consider $\widetilde g = (g, 1) \in C_g \times \C^*$ acting on loops in $X \mmod G$ and impose that this will be the same as acting by $1 \in \R / \abs{\widetilde{g}} \Z$. In other words, we consider Laurent series satisfying a rotation condition with respect to $\widetilde{g}$:

\begin{definition}\label{def:TETKtheoryalpha}
Let $X$ be a space acted on by a finite group $G$ and let $\alpha \colon X \times G^3 \to \C^*$ be a 3-cocycle. We obtain a 2-cocycle $\theta$ on $\Lambda(X \mmod G)$ via transgression. The \emph{$\alpha$-twisted $G$-equivariant Tate K-theory} of $X$, denoted by ${}^\alpha K_{\Tate}(X \mmod G)$, is defined to be the sub-ring of
\[
\bigoplus_{[g]} K\bigl((X^g \mmod C_g)_\theta\bigr) \ls{ q^{\frac{1}{\abs{\widetilde{g}}}} }
\]
in which the $[g]$-summand is the Grothendieck group of Laurent series
\[
\sum_{j>n} V_j q^{\frac{j}{\abs{\widetilde{g}}}}, \quad n \in \Z,
\]
satisfying the following rotation condition with respect to $\widetilde{g}$: for each $j$ the coefficient $V_j$ is a vector bundle on $(X^g \mmod C_g)_\theta$ such that $\widetilde{g}$ acts as multiplication by $\zeta^j_{\abs{\widetilde{g}}}$. 
\end{definition}

By Proposition \ref{prop:vbontwistedgrp} we have the inclusion
\[
\bigoplus_{[g]} {}^{\theta}K_{C_g}(X^g) \ls{ q^{\frac{1}{h|g|}} }
\hookrightarrow
\bigoplus_{[g]} K\bigl((X^g \mmod C_g)_\theta\bigr) \ls{ q^{\frac{1}{h|g|}} },
\]
so twisted equivariant Tate K-theory in particular contains Laurent series with coefficients in the twisted $C_g$-equivariant K-theory of $X^g$. Restricting to this sub-ring amounts to imposing that the morphisms $(1,z) \in C_g \times \C^*$ act as multiplication by $z$.

\begin{remark}
To investigate the order of $\widetilde{g}$, note that
\[
\widetilde{g}^k = (g,1)^k = (g^k, \theta(x, g, g) \theta(x, g^2, g) \dotsm \theta(x, g^{k-1}, g) ).
\]
Letting $h$ be the order of $\alpha$ restricted to $X \times \langle g \rangle^3$ we have that $\abs{\widetilde{g}}$ must divide $h|g|$. Therefore, we could instead consider twisted $C_g$-equivariant vector bundles in which $\widetilde{g}$ acts by some power of $\zeta_{h|g|}$. This $h|g|$ arises in Nora Ganter's description of twisted equivariant Tate K-theory of a point \cite{ganter2009}.
\end{remark}

In light of Remark \ref{rem:quasiellipticcohomology}, one might ask about a version of twisted quasi-elliptic cohomology for finite $G$. Define
\[
\Lambda_{\theta, S^1}(X \mmod G) = \Lambda_{\theta}(X \mmod G) \times \mathbb{B}\R / \sim
\]
where $\sim$ is the equivalence relation generated by identifying $\widetilde g \in C_g \times \C^*$ and $1 \in \R$. This is analogous to $\Lambda_{S^1}(X \mmod G)$, which is obtained from $\Lambda(X \mmod G) \times \B\R$ by identifying $g \in C_g$ and $1 \in \R$. Twisted quasi-elliptic cohomology could then be defined as
\[
{}^\alpha \textnormal{QEll}_G(X) := K\bigl( \Lambda_{\theta, S^1}(X \mmod G) \bigr)
\]

\subsection{Twisted Equivariant Tate K-Theory in General}

Let $\L$ be a gerbe on $\Lambda(X \mmod G)$ that is obtained from a 2-gerbe on $X \mmod G$ via a transgression procedure. Define $\Lambda_L(X \mmod G)$ to be the $L$-twisted groupoid of $\Lambda(X \mmod G)$. Explicitly, 
\[
\Lambda_L(X \mmod G) = \bigsqcup_{[g]} \, (X^g \mmod C_g)_{L^g} = \bigsqcup_{[g]}\, ( L_g \rightrightarrows X^g),
\]
where each $L^g$ is a line bundle on $X^g \times C_g$. The element $g \in C_g$ defines an element in the center of $X^g \mmod C_g$, see Example \ref{ex:centerofXgCg}. In the previous section, we lifted $g \in C_g$ to an element $\widetilde g$ in the center of $(X^g \mmod C_g)_\theta$ and considered Laurent series satisfying a rotation condition with respect to $\widetilde g$. In the general case, let $\xi_g$ be a lift of $g \in C_g$ to the center of $(X^g \mmod C_g)_{L^g}$. We consider Laurent series satisfying a rotation condition with respect to $\xi_g$:

\begin{definition}
Let $X$ be a space acted on by a finite group $G$ and let $L$ be a gerbe on $\Lambda(X \mmod G)$. The \emph{$L$-twisted $G$-equivariant Tate K-theory} of $X$, denoted ${}^LK_{\Tate}(X \mmod G)$, is defined to be the sub-ring of 
\[
\bigoplus_{[g]} K\bigl( (X^g \mmod C_g)_{L^g} \bigr) \ls{q^{\frac{1}{\abs{\xi_g}}}}
\]
in which the $[g]$-summand is the Grothendieck group of Laurent series
\[
\sum_{j>n} V_j q^{\frac{j}{\abs{\xi_g}}}
\]
satisfying the rotation condition with respect to $\xi_g$: for each $j$ the coefficient $V_j$ is a vector bundle on $(X^g \mmod C_g)_{\theta}$ such that $\xi$ acts as multiplication by $\zeta^j_{\abs{\xi_g}}$.
\end{definition}

This definition is less concrete than Definition \ref{def:TETKtheoryalpha}. In practice, before using the transgression formula one must pass from $X \mmod G$ to a Morita equivalent groupoid over which $L$ can be represented as a Deligne cocycle. Thus, the obtained gerbe on $\Lambda(X \mmod G)$ will be in the form of a 2-cocycle on a groupoid that is Morita equivalent to $\Lambda(X \mmod G)$. In theory, there is a gerbe on $\Lambda(X \mmod G)$ which corresponds to the one obtained via transgression, but it is not clear to us how to find this explicitly.

\subsection{Tate K-theory and Moonshine}

Consider the case of $X$ being a point. We twist by a group 2-cocycle $\alpha \colon G^3 \to \C^*$, which produces a collection of 1-cocycles $\theta \colon C_g \times C_g \to \C^*$. These, in turn, define central extensions $\widetilde{C}_g$.  Take an element of twisted equivariant Tate K-theory of a point,
\[
F \in \bigoplus\nolimits_{[g]} R(\widetilde{C}_g) \ls{ q^{\frac{1}{\abs{\widetilde g}} } },
\]
Let $F(g,\widetilde h; q^{\frac{1}{\abs{\widetilde g}}})$ denote the $[g]$-summand, thought of as a function of $\widetilde h \in \widetilde{C}_g$. We may write
\[
F\bigl(g,\widetilde h; q^{\frac{1}{\abs{\widetilde g}}}\bigr) = \sum_j V_j(\widetilde h) q^{\frac{j}{\abs{\widetilde g}}},
\]
where each $V_j$ is a character of $\widetilde{C}_g$. The rotation condition implies that
\[
F\bigl(g,\widetilde{g} \widetilde h ; q^{\frac{1}{\abs{\widetilde g}}}\bigr) 
= \sum_j V_j(\widetilde{g}\widetilde{h}) q^{\frac{j}{\abs{\widetilde g}}} = \sum_j V_j(\widetilde{h}) \zeta^j_{\abs{\widetilde{g}}} q^{\frac{j}{\abs{\widetilde g}}} = F\bigl(g, \widetilde h; \zeta_{\abs{\widetilde{g}}} q^{\frac{1}{\abs{\widetilde g}}}
\bigr).
\]
If we interpret representations of $\widetilde{C}_g$ as twisted representations of $C_g$, as in Example \ref{ch4:twistedvectorbundleonBG}, then $F(g,h; q^{\frac{1}{\abs{\widetilde g}}})$ is a function of $h \in C_g$. We obtain a function $F(g,h;\tau)$ on the upper half plane by setting $q = e^{2\pi i \tau}$. This satisfies
\[
F(g, gh; \tau) 
= \sum_j V_j(gh) q^{\frac{j}{\abs{\widetilde{g}}}}
= \zeta \sum_j V_j(h) \zeta^j_{\abs{\widetilde{g}}} q^{\frac{j}{\abs{\widetilde{g}}}}
=\zeta \cdot F(g, h; \tau+1)
\]
where $\zeta = \theta(g,h)^{-1}$ is a root of unity. For example, consider the Laurent series in the [1]-summand. The cocycle $\theta$, restricted to $\mathbb{B}G$, is trivial because
\[
\theta(h,k) = \frac{\alpha(1, h, k) \cdot \alpha(h,k,1)}{\alpha(h, 1, k)} = 1,
\]
assuming that $\alpha$ is normalised. Therefore, the $[1]$-summand is an element of $R(G)\ls{q}$ and
\begin{equation}\label{eq:1summand}
F(1,g;\tau) = F(1,g;\tau+1).
\end{equation}
We have entered the world of modular functions. Consider Norton's generalised Moonshine conjecture:\footnote{Originally published, without twisting, in the appendix of \cite{Mason:Norton} Twisting arose naturally in later calculations, see \cite{Mason:unpublished}.}

\begin{conjecture}
For each pair $(g,h)$ of commuting elements of the monster group $\M$ there exists a modular function $f(g,h;\tau)$ with the following properties:
\begin{enumerate}[label=(\alph*)]
    \item Up to a constant factor (which will be a root of unity), there is an equality
    $$
        f(g^a h^c, g^bh^d; \tau) = f\left(  g,h; \frac{a\tau + b}{c\tau + d}\right) 
    $$
    whenever $ad-bc = 1$.
    \item For any group element $g$ and nonzero rational number $l$ the coefficient of $q^l = e^{2\pi i l \tau}$ in $f(g,h; \tau)$ is, as a function of $h$, a character of a central extension of $C_g$. Note that nonzero characters can occur for non-integral $l$, but that generalised characters are not needed.
    \item Conjugation of $g$ and $h$ leaves the function unchanged.
    \item Unless $f(g,h;\tau)$ is a constant function, its invariance group will be a modular group of genus zero, commensurable with the standard modular group $SL_2(\Z)$.
\end{enumerate}
\end{conjecture}

A lot of this looks familiar to us; parts (b) and (c) imply that the family of modular functions in the conjecture is an element of the ring
\[
\bigoplus\nolimits_{[g]} R(\widetilde{C}_g) \ls{q^{\frac{1}{h\abs{g}}}}.
\]
The rotation condition in Tate K-theory is half of condition (a). As a particular example, generalised Moonshine should restrict to ordinary Moonshine when considering the pair $(1,g)$. This is the situation of equation \eqref{eq:1summand} and, in fact, $F(1,1;\tau)$ is expected to be the normalised $j$-function
\[
j(\tau) - 744 = \frac{1}{q}  + 196884\,q + 21493760\,q^2 + 864299970\,q^3 + \dotsc.
\]

\subsection{Connection with Luecke's work}

In his recent paper \cite{Kiranspaper}, Luecke introduces a formulation of twisted equivariant Tate K-theory for when $G$ is a compact connected Lie group. The loop space considered is the \emph{thickened inertia groupoid} $\widetilde{\Lambda}(X \mmod G)$, which carries a strict $S^1$-action that is not present in the ordinary inertia groupoid when $G$ is a Lie group. In the case that $G$ is finite, $\widetilde{\Lambda}(X \mmod G)$ is the same as $\Lambda(X \mmod G)$.

The twists come in Luecke's construction come in the form of a graded $S^1$-central extension, which on a groupoid $\G$ is a graded $S^1$-bundle $L$ over $\G_1$ together with an isomorphism of graded $S^1$-bundles over $\G_2$,
\[
L_g \otimes L_f \to L_{f \cdot g},
\]
satisfying a cocycle condition over $\G_3$. This definition is from \cite[\textsection 2.2]{FHT1}. The data of an $S^1$-graded central extension can be built into a groupoid $\tau = (L \rightrightarrows \G_0)$ together with a functor $P \colon \tau \to \G$ such that $P_1 \colon L \to \G_1$ is a graded $S^1$-bundle. Recalling that $S^1$-bundles are equivalent to Hermitian line bundles, a graded $S^1$-central extension on $\G$ is equivalent to a gerbe with graded structure and a Hermitian metric. The groupoid $\tau$ corresponds to the twisted groupoid $\G_L$ that was defined in Example \ref{ex:gerbeasgroupoid}. The relationship between $S^1$-central extensions and $S^1$-gerbes is made precise in \cite{Behren:Xu:DifferentiableStacksandGerbes}. The major difference between our approach and Luecke's is that we do not start with a gerbe on the loop groupoid; we start with a 2-gerbe on $X \mmod G$ and use transgression techniques to obtain a gerbe on the loop groupoid. 

Given an equivariant graded central extension $\tau \to \widetilde{\Lambda}(X \mmod G)$ Luecke's twisted equivariant Tate K-theory is the completed $S^1$-equivariant K-theory of $\tau$ itself. In the case of finite groups, there weren't any problems with using the completion
\begin{equation}\label{ex:completedktheory}
K_{S^1}(\widetilde\Lambda(X \mmod G)) \otimes_{\Z[q^{\pm}]} Z \ls{q},
\end{equation}
defined in \cite{kitchloo:morava:thomprospectra}. In the case of compact connected Lie groups, Luecke constructs an $S^1$-equivariant K-theory that is ``intrinsically $q$-completed.'' This is required because when $X$ is a point and $G$ is a simply connected Lie group, the twisted $G$-equivariant elliptic cohomology is expected to produce level $k$ positive energy representations of the loop group $LG$, where $k$ is an integer determined by the twist, see \cite{Grojnowski:delocalised}. Such representations are infinite dimensional, hence \eqref{ex:completedktheory} will not meet this expectation; without modification, $K_{S^1}$ can only produce finite dimensional representations. This is not an issue in the finite case, where $LG \cong G$.

\begin{appendices}
    \addtocontents{toc}{\protect\renewcommand{\protect\cftchappresnum}{\appendixname\space}}
    \addtocontents{toc}{\protect\renewcommand{\protect\cftchapnumwidth}{6em}}

    \chapter{Bibundles}\label{appendix:bibundles}

Another approach to generalised maps between Lie groupoids is to view a generalised map $\G \to H$ as a \emph{bibundle} from $\G$ to $\H$. These maps are also called \emph{Hilsum-Skandalis maps} after the authors who first introduced them in the context of groupoids. One can read about this approach in \cite{Lerman}, which is the main reference for this section.

\begin{definition}
Let $\H$ be a Lie groupoid. A right $\H$-space $(P,a)$ is a manifold $P$ with a map $a\colon P \to \H_0$ together with a map 
$$
    P \mathbin{_a\times_t} \H_1  \longrightarrow P, \quad 
    (p, h)  \longmapsto p\cdot h
$$
satisfying 
$$
a(p \cdot g) = s(g), \qquad p \cdot 1_x = p  \quad \text{and} \quad (p \cdot h_1) \cdot h_2 = p \cdot (h_2 h_1)
$$
for $z \xrightarrow{\, h_2 \,} y \xrightarrow{\, h_1 \,} x$ and $a(p) = x$.
\end{definition}

\begin{definition}
A principal $\H$-bundle $(P, a, \pi)$ on a space $B$ is a right $\H$-space $(P,a)$ together with a surjective submersion $\pi\colon P \to B$ such that
\begin{enumerate}[label=(\roman*)]
    \item $\pi$ is $\H$-invariant: $\pi(x \cdot h) = \pi(x)$.
    \item $\H$ acts freely and transitively on the fibres of $\pi$, that is, the map
    $$
    P \mathbin{_a\times_t} \H_1  \longrightarrow P \times_B P, \quad 
    (p, h)  \longmapsto (p, p\cdot h)
    $$
    is a diffeomorphism.
\end{enumerate}
\end{definition}

A principal $\H$-bundle will be depicted diagrammatically as follows:
\begin{center}
    \begin{tikzpicture}
    
    \node (P) at (0,0) {$P$};
    \node (B) at (-1.5,0) {$B$};
    \node (H0) at (1.5,0) {$\H_0$};
    
    \draw [<-] (B) -- (P) node [midway, above] {\footnotesize $\pi$};
    \draw [->] (P) -- (H0) node [midway, above] {\footnotesize $a$};
    
    \node (Pe) at (0, -0.7) {$p$};
    \node (Be) at (-1.5, -0.7) {$\pi(p)$};
    \node (H0e) at (1.5, -0.7) {$a(p)$};
    
    \draw [<-|] (Be) -- (Pe);
    \draw [|->] (Pe) -- (H0e);
    
    \end{tikzpicture}
\end{center}

\begin{example}
$\H_1 \xrightarrow{\,\, t\,\,} \H_0$ is a principal $\H$-bundle with anchor map $s \colon \H_1 \to \H_0$. The action is given by composition of arrows,
\[
\H_1 \, \fp{s}{t} \H_1 \to \H_1, \quad (h_1, h_2) \mapsto h_2 \cdot h_1.
\]
\end{example}

Principal bundles have pullbacks. Let $f \colon A \to B$ be a map and $(P, a, \pi)$ a principal $\H$-bundle on $B$. Then the space $f^*P = A \fp{f}{\pi} P$ has a $\H$-action given by $(a, x)\cdot h = (a, x \cdot h)$. This gives the principal bundle
\begin{center}
    \begin{tikzpicture}
    
    \node (fP) at (0,0) {$A \fp{f}{\pi} P$};
    \node (A) at (-2,0) {$A$};
    \node (H0) at (2,0) {$\H_0$};
    
    \draw [<-] (A) -- (fP);
    \draw [->] (fP) -- (H0);
    
    \node (fPe) at (0, -0.7) {\,\,\,$(x, p)\,\,\,$};
    \node (Ae) at (-2, -0.7) {\,\,\,$x$\,\,\,};
    \node (H0e) at (2, -0.7) {\,\,\,$a(p)$.\,\,\,};
    
    \draw [<-|] (Ae) -- (fPe);
    \draw [|->] (fPe) -- (H0e);
    
    \end{tikzpicture}
\end{center}
A bibundle between two groupoids $\G$ and $\H$ will be a principal $\H$-bundle on $\G_0$ with an additional left action by $\G$. For the left action, we introduce the notion of a left $\G$-space.

\begin{definition}
Let $\G$ be a Lie groupoid. A left $\G$-space $(P,a)$ is a manifold $P$ with a map $a\colon P \to \G_0$ together with a map
$$
    \G_1 \mathbin{_s\times_a} P  \longrightarrow P, \quad 
    (g,p)  \longmapsto g \cdot p
$$
satisfying 
$$
a(g \cdot p) = t(g), \qquad 1_x \cdot p = p  \quad \text{and} \quad g_1 \cdot (g_2 \cdot p) = (g_2 g_1) \cdot p
$$
for $x \xrightarrow{\, g_2 \,} y \xrightarrow{\, g_1 \,} z$ and $a(p) = x$.
\end{definition}

\begin{definition}
Let $\G$ and $\H$ be groupoids. A bibundle from $\G$ to $\H$ is a $\H$-principle bundle $(P, a_R, a_L)$ on $\G_0$ such that $(P, a_L)$ is a left $\G$-space satisfying the following conditions:
\begin{enumerate}[label=(\roman*)]
    \item $a_R$ is $\G$-invariant: $a_R(g \cdot p) = a_R(p)$.
    \item The actions commute: $(g \cdot p) \cdot h = g \cdot (p \cdot h)$.
\end{enumerate}
\end{definition}

A bibundle may be drawn as
\begin{center}
    \begin{tikzpicture}
    
    \node (P) at (0,0) {$P$};
    \node (G) at (-1.5,0) {$\G$};
    \node (H) at (1.5,0) {$\H,$};
    
    \draw [<-] (B) -- (P) node [midway, above, xshift=3pt] {\footnotesize $a_L$};
    \draw [->] (P) -- (H) node [midway, above] {\footnotesize $a_R$};
    
    \end{tikzpicture}
\end{center}
with the understanding the $a_L$ and $a_R$ actually map to $\G_0$ and $\H_0$ respectively.

\begin{example}
A bibundle can be obtained from a homomorphism $f\colon \G \to \H$ by pulling back the principal $\H$-bundle $\H_1 \xrightarrow{\,\, t\,\,} \H_0$ along the map $f_0\colon \G_0 \to \H_0$. Explicitly we obtain
\begin{center}
    \begin{tikzpicture}
    
    \def \x {2}
    \def \l {\hspace{0cm}}
        
    \node (P) at (0,0) {$\G_0 \fp{f}{t} \H_1$};
    \node (G) at (-\x,0) {$\G_0$};
    \node (H) at (\x,0) {$\H_0$};
    
    \draw [<-] (G) -- (P);
    \draw [->] (P) -- (H);
    
    \node (Pe) at (0, -0.7) {\l $(x, h)$ \l};
    \node (Ge) at (-\x, -0.7) {\l $x$ \l};
    \node (He) at (\x, -0.7) {\l $s(h)$. \l};
    
    \draw [<-|] (Ge) -- (Pe);
    \draw [|->] (Pe) -- (He);
    
    \end{tikzpicture}
\end{center}
The left $\G$-action is given by $g \cdot (s(g), h) = (t(g), h \cdot f_1(g))$. This bibundle will be denoted $\langle f \rangle$ and is the first step in connecting generalised maps and bibundles.
\end{example}

\begin{definition}
Two bibundles $P,Q \colon \G \to \H$ are isomorphic if there is a diffeomorphism $\alpha\colon P \to Q$ that is $\G$- and $\H$-equivariant, which means that
$$
\alpha(g \cdot p \cdot h) = g \cdot \alpha(p) \cdot h,
\quad  (g,p,h) \in \G_1 \times_{\G_0} P \times_{\H_0} \H_1.
$$
\end{definition}
Pictorially we may denote an isomorphism of bibundles by
\begin{center}
\begin{tikzpicture}
    \node (L) at (0,1)  {$\G_0$};
    \node (R) at (5,1) {$\H_0$.};
    \node (U) at (2.5, 1.75) {$Q$};
    \node (D) at (2.5, 0.25) {$P$};
    \draw [<-] (L) -- (U);
    \draw [->] (U) -- (R);
    \draw [<-] (L) -- (D);
    \draw [->] (D) -- (R);
    \draw [->] (D) -- (U) node [midway, right] {\footnotesize $\alpha$};
\end{tikzpicture}
\end{center}
The resemblance with generalised maps of groupoids should be clear.

\begin{definition}
A Hilsum-Skandalis map from $\G$ to $\H$ is an isomorphism class of bibundles from $\G$ to $\H$.
\end{definition}

Bibundles can be composed: if $P$ is a bibundle from $\G$ to $\H$ and $Q$ is a bibundle from $\H$ to $\K$, then define
\[
Q \circ P := (P \times_{\H_0} Q) / \H
\]
where the fibre product fits into the diagram
\begin{center}
    \begin{tikzcd}
     P \times_{H_0} Q  \arrow[r] \arrow[d] & Q \arrow[d] \\
     P \arrow[r] & \H_0
    \end{tikzcd}
\end{center}
and has a $\H$-action $(p, q) \cdot h = (p \cdot h, h^{-1} \cdot q)$. Then $Q \circ P$ inherits a left $\G$-action from $P$ and a right $\K$-action from $Q$ and is a bibundle from $\G$ to $\K$. This composition is not strictly associative, but it is up to isomorphism.

\begin{definition}
A bibundle $P \colon \G \to \H$ is invertible if the right anchor map $a_R \colon P \to \H_0$ is a $\G$-principal bundle. In this case, the bibundle $P^{-1} \colon \H \to \G$ is obtained by switching the anchor maps, changing the left $\G$-action to a right action and turning the right $\H$-action into a left action.
\end{definition}

An invertible bibundle is directly analogous to an equivalence of Lie groupoids, as seen in the following \cite[Lemma 3.34]{Lerman}.

\begin{lemma}
A functor $f \colon \G \to \H$ is an equivalence of Lie groupoids if and only if the corresponding bibundle $\langle f \rangle \colon \G \to \H$ is invertible.
\end{lemma}

Let HS denote the category of Lie groupoids with Hilsum-Skandalis and let Gp be the category of Lie groupoids with generalised maps. Given a generalised map $(K, \varepsilon, \phi)$ from $\G$ to $\H$,
\[
\G \xleftarrow{\,\, \varepsilon \,\,} \K \xrightarrow{\,\, \phi \,\,} \H,
\]
the bibundle $\langle \phi \rangle \circ \langle \varepsilon \rangle^{-1} \colon \G \to \H$ is well-defined since $\varepsilon$ is an equivalence. This gives rise to a functor $F \colon \text{Gp} \to \text{Hs}$ which is the identity on objects and maps the generalised map $(\K, \varepsilon, \phi)$ to the Hilsum-Skandalis map $\langle \phi \rangle \circ \langle \varepsilon \rangle^{-1}$. The following result \cite[Prop 3.39]{Lerman}  makes precise the equivalence of the two approached to orbifold maps.

\begin{proposition}\label{prop:equivalenceHSandGp}
The functor $F \colon \textnormal{Gp} \to \textnormal{HS}$ is an equivalence of categories.
\end{proposition}

We can use this equivalence to prove an enlightening fact about generalised maps. First note the following \cite[Lemma 3.19]{Lerman}:

\begin{lemma}\label{lemma:lerman:3.37}
Let $P \colon \G \to \H$ be a bibundle from $\G$ to $\H$. there exists an open cover $\U$ of $\G_0$ and a functor $f \colon \G[\U] \to \H$ such that there is an isomorphism of bibundles,
\[
P \cong \langle f \rangle \circ \langle \varepsilon \rangle^{-1}
\]
where $\varepsilon \colon \G[\U] \to \G$ is the induced equivalence.
\end{lemma}

\begin{proposition}
Any generalised map $\G \to \H$ can be represented by an open cover $\mathcal{U}$ of $\G_0$ and a homomorphism $\phi$ so that we have
$$
\G \xleftarrow{\,\,\varepsilon\,\,} \G[\U] \xrightarrow{\,\,\phi\,\,} \H.
$$
where $\varepsilon$ is defined in Example \ref{ch2:groupoidopencoverequivalence}. Furthermore, another such $(\G[\U'], \varepsilon', \phi')$ represents the same generalised map if and only if $\phi$ and $\phi'$ are related by a natural transformation when restricted to a common refinement of $\mathcal{U}$ and $\mathcal{U}'$.
\end{proposition}
\begin{proof}
Consider a generalised map from $\G$ to $\H$.
$$
\G \xleftarrow{\,\,\varepsilon\,\,} \K \xrightarrow{\,\,\phi\,\,} \H.
$$
Since $\varepsilon$ is an equivalence the bibundle $\langle \varepsilon \rangle$ is invertible. Therefore $\langle \phi \rangle \circ \langle \varepsilon \rangle^{-1}$ is a well-defined bibundle from $\G$ to $\H$. By Lemma \ref{lemma:lerman:3.37} there is an open cover $\U$ of $\G_0$ and a functor $f\colon \G[\U] \to \H$ such that 
$$
\langle \phi \rangle \circ \langle \varepsilon \rangle^{-1}
\cong 
\langle f \rangle \circ \langle \psi \rangle^{-1}
$$
where $\psi$ is the equivalence $\G[\U] \to \G$. By the equivalence of categories in Proposition \ref{prop:equivalenceHSandGp} this means that 
$$
\G \xleftarrow{\,\,\psi\,\,} \G[\U] \xrightarrow{\,\,f\,\,} \H.
$$
represents the same generalised map we started with. This proves the first part of the proposition. Let $(G[\U], \varepsilon, \phi)$ and $(G[\U'], \varepsilon', \phi')$ be equivalent generalised maps where $\varepsilon$ and $\varepsilon'$ are the usual equivalences for open cover groupoids. By definition there is a smooth functor $f \colon G[\U] \to G[\U']$ such that the following commutes up to natural transformation:
\begin{center}
\begin{tikzpicture}
    \node (L) at (0,1)  {$\G$};
    \node (R) at (5,1) {$\H$};
    \node (U) at (2.5, 1.75) {$\G[\U']$};
    \node (D) at (2.5, 0.25) {$\G[\U]$};
    \draw [->] (U) -- (L) node [midway, above] {\small$\varepsilon'$};
    \draw [->] (U) -- (R) node [midway, above] {\small$\phi'$};
    \draw [->] (D) -- (L) node [midway, below] {\small$\varepsilon$};
    \draw [->] (D) -- (R) node [midway, below, yshift=1mm] {\small$\phi$};
    \draw [->] (D) -- (U) node [midway, right] {\small$f$};
\end{tikzpicture}
\end{center}
In particular by looking at the left triangle there is a natural transformation $\eta$ from $\varepsilon$ to $\varepsilon' \circ f$. This is a map  $\eta\colon \G[\U]_0 \to \G$ such that $\eta(x)$ is a morphism in $\G$ from $\varepsilon_0(x)$ to $\varepsilon'_0(f_0(x))$. For $\U \cap \U'$ a common refinement of $\U$ and $\U'$ we want the following to commute up to natural transformation:
\begin{equation}\label{eq:refinement}
\begin{tikzcd}
    \G[\U \cap \U'] \arrow[hook, r, "\iota'"] \arrow[hook, d, "\iota"]          & \G[\U'] \arrow[d, "\phi'"]        \\
    \G[\U] \arrow[ru, "f"]    \arrow[r, "\phi"]                             & \H
\end{tikzcd}
\end{equation}
The maps $\iota, \iota'$ are the natural inclusions. We already know that the lower right triangle commutes up to natural transformation, so it's sufficient to find a natural transformation $\xi$ from $f \circ \iota$ to $\iota'$. This exists because if $x \in \G[\U \cap \U']$, then $\eta(\iota(x))$ is a morphism in $\G$ from $\varepsilon(f(x))$ to $\varepsilon'(x)$ and this morphism exists in $\G[\U']$ as well.

The other direction follows from the definition of a generalised map; if $\phi$ and $\phi'$ are related by a natural transformation when restricted to $\U \cap \U'$, then \eqref{eq:refinement}, without the $f$, commutes up to natural transformation. Therefore, the two generalised map are both equivalent to the same generalised map, which associated to the common refinement.

\end{proof}

It is convenient to talk about bibundles in the context of vector bundles on groupoids: a vector bundle on $\G$ is equivalent to a $\GL_n(\C)$-bundle on $\G$, and a $\GL_n(\C)$-bundle on $\G$ is equivalent to a bibundle from $\G$ to $\mathbb{B}\GL_n(\C)$. Then, if $\G$ and $\H$ are Morita equivalent, then they are isomorphic in HS and there are bijections
\[
\vect \G \cong \text{HS-maps}(\G, \mathbb{B}\GL_n(\C)) \cong \text{HS-maps}(\H, \mathbb{B}\GL_n(\C)) \cong \vect \H.
\]
These are in fact equivalences of categories, but we will not get into this discussion. We will just state the result that tells us that vector bundles on groupoids are Morita invariant, and refer the reader to \cite[\textsection 2.6]{amenta} for a more detailed discussion.

\begin{proposition}
If $\G$ and $\H$ are Morita equivalent then there is an equivalence of categories between $\vect{\G}$ and $\vect{\H}$.
\end{proposition}

    \chapter{Bundle Gerbes on Orbifolds}\label{appendix:bundlegerbes}

In this section, we recall the definition of a bundle gerbe on both smooth manifolds and groupoids. This will lead to an explanation of the definition used in Chapter 3. These were first introduced by Murray in \cite{murraybg} though our main reference is his expository paper \cite{Murrayintotobg}. 

Given a map $\pi \colon Y \to M$ we write $Y^{[n]}$ for the $n$-fold fiber product
\[
Y^{[n]} = Y \fp{\pi}{\pi} Y \fp{\pi}{\pi} \dotsm \fp{\pi}{\pi} Y.
\]
We now define a bundle gerbe:

\begin{definition}
A bundle gerbe on a smooth manifold $M$ is a pair $(P,Y)$ where $Y \to M$ is a surjective submersion and $P$ is a line bundle on $Y^{[2]}$. These come equipped with a bundle gerbe multiplication, which is an isomorphism
\[
\pi_3^* P \otimes \pi_1^*P \to \pi_2^*P
\]
satisfying the following associativity condition on fibers:
\[
\begin{tikzcd}
P_{(y_1, y_2)} \otimes P_{(y_2, y_3)} \otimes P_{(y_3, y_4)} \arrow[r] \arrow[d] & P_{(y_1, y_2)} \otimes P_{(y_2, y_4)} \arrow[d] \\
P_{(y_1, y_3)} \otimes P_{(y_3, y_4)} \arrow[r]           & P_{(y_1, y_4)}          
\end{tikzcd}
\]
\end{definition}

\begin{example}\label{appendix:hitchinchatterjeegerbes}
Let $\{ U_i \}$ be an open cover of $M$ and let 
\[
Y = \bigsqcup\nolimits_{i} U_i.
\]
Then the natural inclusion $Y \to M$ is a surjective submersion and a line bundle
\[
P \to Y^{[2]} = \bigsqcup\nolimits_{i,j} U_i \cap U_j 
\]
is a collection of line bundles $P_{ij} \to U_i \cap U_j$. The bundle gerbe multiplication is then an isomorphism
\[
P_{ij} \otimes P_{jk} \cong P_{ik}.
\]
In this way we obtain a Hitchin-Chatterjee gerbe.
\end{example}

The original purpose of bundle gerbes was to provide a geometric interpretation of degree three integral cohomology. The correct notion of isomorphism which allows one to identify isomorphism classes of bundle gerbes with integral cohomology is the following.

\begin{definition}
A stable isomorphism between bundle gerbes $(P_1, Y_1)$ and $(P_2, Y_2)$ is a line bundle $L$ on $Z = Y_1 \times_M Y_2$ together with an isomorphism
\[
p_1^* P_1 \otimes \zeta_2^* L \to \zeta_1^* L \otimes p_2^*P_2
\]
of line bundles over $Z^{[2]}$ where $p_1, p_2, \zeta_1, \zeta_2$ are the relevant projection maps. The isomorphism must be compatible with the respective bundle gerbe multiplications. 
\end{definition}

If we only care about bundle gerbes up to stable isomorphism, then the only bundle gerbes we need to care about are Hitchin-Chatterjee gerbes of example \ref{appendix:hitchinchatterjeegerbes}.

\begin{theorem}[{\cite[Prop 5.5]{Murrayintotobg}}]
Every bundle gerbe is stably isomorphic to a Hitchin-Chatterjee gerbe.
\end{theorem}

Suppose that $(P, Y_1)$ and $(Q, Y_2)$ are Hitchin-Chatterjee gerbes with
\[
Y_1 = \bigsqcup\nolimits_{i} U_i, \quad Y_2 = \bigsqcup\nolimits_{\alpha} V_\alpha,
\]
where $\{U_i\}$ and $\{V_\alpha\}$ are two open covers. A stable isomorphism between these gerbes is a line bundle
\[
A \to Y_1 \times_M Y_2 = \bigsqcup\nolimits_{i, \alpha} U_i \cap V_\alpha.
\]
In other words, a collection of line bundles $A_{i\alpha} \to U_i \cap V_\alpha$. These have isomorphisms
\[
P_{ij} \otimes A_{j, \beta} \cong A_{i\alpha} \otimes Q_{\alpha\beta}
\]
that are compatible with the isomorphisms $P_{ij} \otimes P_{jk} \cong P_{ik}$ and $Q_{\alpha\beta} \otimes Q_{\beta\gamma} \cong Q_{\alpha\gamma}$.

\begin{definition}
Let $P = (P_{ij})$ and $Q = (Q_{\alpha\beta})$ be Hitchin-Chatterjee gerbes associated to open covers $\{U_i\}$ and $\{ V_\alpha \}$ respectively. A stable isomorphism from $P$ to $Q$ is a collection of line bundles $A_{i\alpha} \to U_i \cap V_\alpha$ together with an isomorphism
\[
P_{ij} \otimes A_{j\beta} \cong A_{i\alpha} \otimes Q_{\alpha\beta}
\]
compatible with the respective gerbe multiplication maps.
\end{definition}

Consider a map $f \colon M \to N$ and a Hitchin-Chatterjee gerbe $P = (P_{ij})$ on $N$ associated to the open cover $\{U_i\}$. Pulling back along $f$ we can obtain an open cover $\{f^{-1}(U_i)\}$ of $M$ with line bundles
\[
f^*P_{ij} \to f^{-1}(U_i \cap U_j).
\]
Thus, we have a pull-back Hitchin-Chatterjee gerbe $f^*P$ on $M$.

Using the usual formalism for structures over groupoids, we define a bundle gerbe over a groupoid as follows.

\begin{definition}\label{def:bundlegerbeongroupoid}
A bundle gerbe on a groupoid $\G$ consists of the following data:
\begin{enumerate}[label=(\roman*)]
    \item A bundle gerbe $(P,Y)$ on $\G_0$.
    \item A stable isomorphism $A \colon t^* (P,Y) \to s^*(P,Y)$ of bundle gerbes over $\G_1$ which satisfies:
    \begin{itemize}
        \item $u^* A = \id_{(P,Y)}$.
        \item There is an isomorphism $d_2^*A \otimes d_0^*A \to m^*A$ of line bundles over $\G_2$ such that the following diagram commutes:
        \[
        \begin{tikzcd}
         A_g \otimes A_h \otimes A_k   \ar[r] \ar[d]    & A_g \otimes A_{hk} \ar[d] \\
         A_{gh} \otimes A_k     \ar[r]             & A_{ghk}
        \end{tikzcd}
        \]
    \end{itemize}
\end{enumerate}
\end{definition}

Note that there is an extra compatibility condition that isn't included in, for instance, the definition of vector bundles on groupoids (Definition \ref{def:vectorbundledef2}). The reason for this the fact that bundle gerbes form a 2-category - there is a notion of morphisms between stable isomorphisms. This extra condition can be thought of as a compatibility condition on the 2-morphisms.

We investigate this definition for Hitchin-Chatterjee gerbes. Start with a bundle gerbe $(P,Y)$ on $\G_0$ with
\[
Y = \bigsqcup\nolimits_i U_i
\quad \text{and} \quad
P = \bigsqcup\nolimits_{ij} P_{ij} \to Y^{[2]},
\]
where $\U = \{U_i\}$ is an open cover of $\G_0$. A stable isomorphism from $t^*P$ to $s^*P$ is a collection of line bundles $A_{ij} \to t^{-1}(U_i) \cap s^{-1}(U_j)$ with a collection of isomorphisms detailed above. These line bundles must satisfy
\[
A_{ij} \otimes A_{jk} \cong A_{ik}.
\]
By assuming that $\U$ is suitably refined, we may assume that $P$ is a trivial line bundle. In which case our Hitchin Chatterjee gerbe is given by a line bundle
\[
A = \bigsqcup\nolimits_{i,j} A_{ij} \to \bigsqcup\nolimits_{i,j} t^{-1}(U_i) \cap s^{-1}(U_j)
\]
with isomorphisms $A_{i j} \otimes A_{jk} \cong A_{ik}$. In other words, we have a gerbe on $\G[\U]$ as defined in Definition \ref{def:gerbeongroupoid}. This means that both definitions of gerbes on orbifolds are equivalent up to the choice of representing groupoid.

Let us now introduce connective structure on bundle gerbes. We will repeat the previous theory with connective structure added.

\begin{definition}\label{def:bundlegerbeswithconnection}
Let $(P,Y)$ be a bundle gerbe on $M$. A connection on $(P,Y)$ is a connection $\nabla$ on $P$ that is compatible with the bundle gerbe multiplication, which means that 
\[
\pi_3^* P \otimes \pi_1^*P \to \pi_2P
\]
is an isomorphism of line bundles with connection.
\end{definition}

\begin{definition}\label{def:stableisomofbgwithconn}
A stable isomorphism between bundle gerbes $(P_1, Y_1)$ and $(P_2, Y_2)$ with connection $\nabla_1$ and $\nabla_2$ respectively is a line bundle $L$ with connection $\nabla$ on $Z = Y_1 \times_M Y_2$ together with an isomorphism 
\[
p_1^*P_1 \otimes \zeta_2^*L \to \zeta_1^*L \otimes p_2^*P_2
\]
of line bundles with connection over $Z^{[2]}$ where $p_1, p_2, \zeta_1, \zeta_2$ are the relevant projection maps. There are two conditions that must be satisfied,
\begin{enumerate}[label=(\roman*)]
    \item The connection on $L$ satisfies $\nabla = p_2^* \nabla_2 - p_1^* \nabla_1$.
    \item The isomorphism above is compatible with the respective bundle gerbe multiplications.
\end{enumerate}
\end{definition}

A bundle gerbe with connection on a groupoid is obtained by repeating the definition of a bundle gerbe and adding the words ``with connection'' throughout. 

Suppose that $\G$ is a Leray groupoid. A bundle gerbe $(P,Y)$ on $\G_0$ is stably isomorphic to the trivial bundle gerbe because $\G_0$ is a disjoint union of contractible spaces. Therefore, assume that $Y = \G_0$ and $P = \G_0 \times \C$. Suppose that $P$ has curvature $B \in \Omega^2(\G_0)$. A stable isomorphism between $t^*P$ and $s^*P$, both trivial, is a line bundle $L$ on $\G_1$ with curvature $t^*B - s^*B$. Again, $L$ must be the trivial line bundle because $\G$ is Leray and hence there is a 1-form $A \in \Omega^1(\G_1)$ such that $-dA = t^*B - s^*B$. The associativity conditions on this stable isomorphism imply that there is an isomorphism
\[
d_2^*L \otimes d_0^*L \cong \G_2 \times \C \rightarrow \G_2 \times \C \cong m^*L
\]
of line bundles with connection. This is given by a 2-cocycle $g \colon \G_2 \to \C^*$. Compatibility with the connective structure shows, via a similar calculation to Example \ref{ch3:linebundleconnectionlocal}, that
\[
d_0^*A - m^*A + d_2^*A = d\log f.
\]
We have obtained a 2-cocycle $g$ with differential forms $A \in \Omega^1(\G_1)$ and $B \in \Omega^2(\Omega^2)$ satisfying
\[
d_0^*A - m^*A + d_2^*A = d\log f \stext{and}  t^*B - s^*B = -dA.
\]
This is the reasoning behind the definition of a gerbe with connection in Chapter 4.

    \clearpage
\end{appendices}

\chapter*{Notation Index}
\addcontentsline{toc}{chapter}{Notation Index}

\newcommand{\ite}[1]{\item[#1 \hfill]}

\begin{list}{}
   {\setlength{\labelwidth}{3cm} 
    \setlength{\labelsep}{0cm}
    \setlength{\leftmargin}{3.5cm}
    \setlength{\rightmargin}{0cm}
    \setlength{\itemsep}{0cm} 
    \setlength{\parsep}{0.05cm} 
    \setlength{\itemindent}{0cm}
    \setlength{\listparindent}{0cm}}
\ite{$\G,\H,\K$}    {groupoids},
\ite{$s,t,u,i,m$}   {structure maps of a groupoid},
\ite{$\G_n$}        {the space of sequences of $n$ composable arrows in $\G$, see Page \pageref{Gn}},
\ite{$d_j$}         {maps in the nerve of $\G$, see Page \pageref{Gn}},
\ite{$\U, \V$}      {open covers},
\ite{$X \mmod G$}   {global quotient orbifold, Example \ref{ch2:actiongroupoid}},
\ite{$X_G$}        {the Borel construction of $X$},
\ite{$\mathbb{B}G$} {the groupoid $* \mmod G$, Example \ref{ch2:group}},
\ite{$\G{[}\U{]}$}      {open cover groupoid, Example \ref{ch2:opencovergroupoid}},
\ite{$\Lambda \G$}  {the inertia groupoid, Example \ref{ch2:inertiagroupoid}},
\ite{$X^g$}         {elements $x \in X$ such that $x \cdot g = x$},
\ite{$\S^1_{\U}$}   {see Page \pageref{S1U}},
\ite{$\L\G(\U)$}    {the groupoid of loops in $\G$ associated to the open cover $\U$, Definition \ref{def:LG(U)}}
\ite{$\L\G$}        {the loop groupoid, Definition \ref{LoopGroupoids}}
\ite{$\P_g$}        {the space of maps $\varphi\colon \R \to X$ such that $\varphi(x+1) = \varphi(x)\cdot g$, Page \pageref{Pg}}
\ite{$C_g$}         {the centraliser of $g$ in $G$},
\ite{$\vect{\G}$}   {the category of vector bundles on $\G$},
\ite{$\co, \ct$}    {typically used to denote cocycles},
\ite{$\widetilde{G}$} {a central extension of $G$ by $\C^*$},
\ite{$\Omega^k$}      {complex valued differential $k$-forms},
\ite{${}^\alpha K$} {twisted K-theory},
\ite{$\widehat{K}_{S^1}$}   {completed $S_1$-equivariant K-theory, Page \pageref{KS1hat}},
\ite{$\zeta_j$} {the primitive $j$th root of unity $e^{2\pi i / j}$},
\ite{$K_{\Tate}$}   {Tate K-theory, Chapter \ref{TETKTheory}},
\ite{$^{\alpha}K_{\Tate}$}   {Twisted Tate K-theory, Definition \ref{def:TETKtheoryalpha}}.
\end{list}

\clearpage


\addcontentsline{toc}{chapter}{References}  
\begin{singlespace}  
	\setlength\bibitemsep{\baselineskip}  
	\printbibliography[title={References}]
\end{singlespace}

\end{document}